\let\und=\underline

\def\MM{{\overline{M}}}
\def\MMM{{\overline{M}_{g,n}(X,\beta)}}

\def\M0{{\overline{M}_{0,n}(\mathbb{P}^{r},d)}}

\def\mm{{\overline{\mathcal{M}}}}
\def\m0{{\overline{\mathcal{M}}_{0,n}(\mathbb{P}^{r},d)}}
\def\mznXb{{\overline{\mathcal{M}}_{0,n}(X,\beta)}}
\def\mmm{{\overline{\mathcal{M}}_{g,n}(X,\beta)}}
\def\mg{{\overline{\mathcal{M}}_{g,n}(\mathbb{P}^{r},d)}}

\def\nobarm0{{\mathcal{M}_{0,n}(\mathbb{P}^{r},d)}}
\def\nobarmm{{\mathcal{M}}}

\def\eM{{\overline{\mathbf M}_{g,n}(X,\beta)}}

\def\D{{\Delta}}
\def\G{{\Gamma}}
\def\FF{{\mathbb F}}
\def\GG{{\mathbb G}}

\def\La{{\Lambda}}
\def\O{{\mathcal O}}

\def\X{{\mathcal X}}
\def\AA{{\mathbb A}}
\def\CC{{\mathbb C}}

\def\N{{\mathcal N}}
\def\NN{{\mathbb N}}
\def\PP{{\mathbb P}}
\def\QQ{{\mathbb Q}}
\def\R{{\mathcal R}}

\def\ZZ{{\mathbb Z}}
\def\C{{\mathcal C}}
\def\PPr{{{\mathbb P}}^{r}}

\def\a{\alpha}
\def\b{\beta}
\def\f{\phi}
\def\la{\lambda}
\def\d{\delta}
\def\dd{\partial}
\def\e{\epsilon}
\def\g{\gamma}
\def\n{{\eta}}
\def\r{\rho}
\def\s{\sigma}
\def\t{\tau}
\def\w{\omega}

\def\pr{^{\prime}}
\def\dpr{^{\prime\prime}}

\def\x{\times}
\def\*{\otimes}
\def\iso{\simeq}
\def\sub{\subset}

\def\+{\oplus}

\def\bra{\langle}
\def\ket{\rangle}

\def\ra{\rightarrow}
\def\lra{\longrightarrow}

\def\pj{\operatorname{pr}}
 
\def\ker{\operatorname{ker}} 

\def\Spec{\operatorname{Spec}}

\def\Ext{\operatorname{Ext}}

\def\gr{\operatorname{gr}}

\def\id{\operatorname{id}}
\def\ev{\operatorname{ev}}
\def\dim{\operatorname{dim}}

\def\deg{\operatorname{deg}}

\def\Bl{\operatorname{B\ell}}
\def\eqeul{\operatorname{Euler_T}}
\def\ser{\operatorname{Serre}}

\def\PGL{\operatorname{PGL}}
\def\Div{\operatorname{Div}}
\def\div{\operatorname{div}}

\def\val{\operatorname{val}}
\def\GL{\operatorname{GL}}
\def\SL{\operatorname{SL}}
\def\point{\operatorname{point}}
\def\Sym{\operatorname{Sym}}
\def\Aut{\operatorname{Aut}}

\def\Trace{\operatorname{Tr}}
\def\st{\operatorname{st}}
\def\nor{\operatorname{nor}}

\def\one{\textrm{\makebox[0.02in][l]{1}1}}

\documentclass[12pt]{report}

\usepackage{amsmath,amssymb}
\usepackage{latexsym,amsfonts}
\usepackage{amsthm,xypic}
\usepackage{pstricks,pst-node,pst-tree}
\usepackage{graphicx}
\usepackage{rotating}

\begin{document}

\newtheorem{thm}{Theorem}
\newtheorem{cor}{Corollary}
\newtheorem{Def}{Definition}
\newtheorem{eg}{Example}
\newtheorem{prop}{Proposition}
\newtheorem{rmk}{Remark}
\newtheorem{lem}{Lemma}
\newtheorem{conj}{Conjecture}

\newcommand{\chews}[2]{\genfrac{[}{]}{0pt}{}{#1}{#2}}
\newcommand{\lchoose}[2]{{#1\choose#2}}



\setlength{\oddsidemargin}{.5in}
\setlength{\evensidemargin}{.5in}
\setlength{\textwidth}{4in}
\setlength{\topmargin}{1.062in}
\setlength{\textheight}{6.9375in}
\renewcommand{\baselinestretch}{1}

\pagenumbering{roman}

\title{
\vspace{-2.75in}
A PRESENTATION FOR\\
\mbox{} \\
THE CHOW RING\\
\mbox{} \\
$A^*(\mm_{0,2}(\PP^1,2))$}

\author{
\vspace{-.5in}\\
By\\
\\
JONATHAN ANDREW COX\\
\\
Bachelor of Science\\ 
\\
Wisconsin Lutheran College\\ 
\\
1997
}

\date{
\mbox{}\\ 
Submitted to the Faculty of the\\ 
Graduate College of the\\ 
Oklahoma State University\\ 
in partial fulfillment of\\ 
the requirements for\\ 
the degree of\\
DOCTOR OF PHILOSOPHY\\
July 2004}

\maketitle


\setlength{\oddsidemargin}{.53125in}
\setlength{\evensidemargin}{.53125in}
\setlength{\textwidth}{5.96875in}
\setlength{\topmargin}{.03125in}
\setlength{\textheight}{8.96875in}
\setlength{\footskip}{2in}


\begin{center}
{\Large
\mbox{}

\vspace{.25in}

ACKNOWLEDGMENTS}
\end{center}

I dedicate this dissertation with love to my wife Tedi. Thank you for enduring
the insanity, stress, dysfunctionality, and long hours associated with
this doctoral quest, while at the same time taking on extra responsibilities
and supporting me in so many ways. I am indebted to you for at least the
next seven years (but likely forever)!

I am also deeply grateful to my dissertation adviser, Dr. Sheldon Katz, for 
presenting me with the invaluable opportunity to work with him as well as
with many other opportunities, for plenteous financial support, for patience
with my sometimes slow learning, for encouragement when I was ready to give
up hope, for inspiration during the many times 
when I was at an impasse, and for
his contagious love and excitement for
mathematics and dedication to his vocation.

I am further beholden to Dr. William Jaco and Dr. 
Alan Adolphson for providing 
additional financial 
support during my graduate school years, and to Dr. Adolphson for his willing
and helpful service as the chair of my committee.
I appreciate the generous support given to me by the University of Illinois 
mathematics department during my two and a half years as a visiting graduate
student. I also acknowledge 


\setlength{\textheight}{8.84375in}
\setlength{\footskip}{-3in}

\noindent
with humble gratitude the
Oklahoma State University mathematics department for continued 
faith in me
and support during this period of absence.

I am thankful for my high school mathematics teacher
Mr. Matthew Schlawin and my undergraduate adviser Dr. Mel Friske for
shaping my mathematical skills and 
inspiring me to pursue advanced study in mathematics.

Finally, thanks go to my parents, J. Gary and Joan Cox, for instilling in me
a love of learning and encouraging me to play with geometric objects at
a young age.



\setlength{\textheight}{8.5in}

\vspace{-.5in}

\tableofcontents

\pagebreak


\setlength{\textheight}{8.84375in}

\vspace{-.5in}

\listoftables


\pagebreak



\setlength{\footskip}{0in}

\setlength{\oddsidemargin}{0.53125in}
\setlength{\evensidemargin}{0.53125in}
\setlength{\textwidth}{5.96875in}
\setlength{\topmargin}{.03125in}
\setlength{\headheight}{0in}
\setlength{\headsep}{0in}
\setlength{\topskip}{0in}
\setlength{\textheight}{8.96875in}

\chapter{Introduction}
\label{sec:intro}

\pagenumbering{arabic}
\setlength{\footskip}{.5in}

Let $\mmm$ be the moduli space of stable maps from $n$-pointed, 
genus $g$ curves to $X$ of class $\b$, to be defined in Section 
\ref{sec:modintro}. If $X=\PP^r$, we can identify the homology class
$\b$ with an integer $d=d[line]$, the {\em degree} of the stable map.
In this case we write the moduli space of stable maps as $\mg$.
In this dissertation we will study the intersection theory of the moduli
spaces $\mm_{0,2}(\PP^r,2)$.

Moduli spaces of stable maps have proven useful in studying
both superstring theory and enumerative geometry. Here we will only mention
briefly the physics side of the story.

In the last fifteen years, superstring theory emerged
as a serious contender
for the ``Theory of Everything," {\em i.e.}, 
a fundamental physical description of the 
universe.  At the
heart of the extremely complex string theory revolution lies a simple 
concept. Elementary ``particles,'' the basic constituents from which 
everything in our universe is built, are actually tiny vibrating loops.
More recently it has become apparent that string
theory is not quite the final answer. Higher dimensional
versions of strings, called $D$-branes, are added in the more comprehensive
$M$-theory. However, the theory of strings continues to
play a key role in this bigger picture. See \cite{G} for a non-technical
introduction to string theory and $M$-theory.

As a string propagates through time, it traces out a {\em world sheet}, which
can be viewed 
mathematically as a Riemann surface or an algebraic curve.  Since
the world sheet lives in space-time, it seems reasonable to study
algebraic curves living
inside manifolds or smooth varieties. However, 
preserving the full data of the world sheet requires using {\em maps}
of curves into space-time instead. Considering the world sheet
as a curve already inside the ambient space
amounts to looking at the
image of the map. Much like representing the path of an object as a zero set
of equations
rather than with
a parametric curve, this viewpoint discards important information about the
the world sheet, such as how it may cross itself.
It turns out that
studying maps from algebraic curves
into a space is also necessary in order to 
obtain a sound mathematical theory of enumerative geometry of curves
in that space. 

Physicists want to determine various values 
associated with the 
space of all possible world sheets via the use of {\em correlation functions}.
They compute these correlation
functions by means
of Feynman integrals, which do not have a rigorous mathematical
definition. Developing a solid mathematical foundation for such computations
was the primary motivation behind the introduction of moduli spaces of
stable maps in \cite{KM}. Examples include {\em instanton numbers}, 
which intuitively count the number of 
holomorphic instantons on $X$ (nonconstant holomorphic maps from Riemann
surfaces to $X$). Instanton numbers are calculated using other values
called {\em Gromov-Witten invariants}.
Naively, Gromov-Witten invariants should count the number of curves of
a certain homology class and genus which pass through certain subvarieties of
the target space. More specifically,
let $X$ be a projective manifold, $\beta\in H_2(X)$, $g$ and $n$ nonnegative
integers. Let $\g_1,\ldots,\g_n\in H^*(X)$ be cohomology classes such that
there exist subvarieties $Z_1,\ldots,Z_n$
with $Z_i$ representing the Poincar\'{e} dual of $\g_i$.
Naively, the Gromov-Witten invariant $\bra\g_1,\ldots,\g_n\ket_{g,\b}$ should
count
the number of genus $g$ curves of class $\b$ that intersect all of the $Z_i$.
Also of interest are {\em gravitational correlators}, 
to be defined in Section \ref{sec:app},
which generalize Gromov-Witten invariants and include them as a special case.

Gravitational correlators are defined and computed mathematically as 
intersection numbers on the moduli space of stable maps. For example,
the Gromov-Witten invariant $\bra\g_1,\ldots,\g_n\ket_{g,\b}$ described
above is given by
\[\bra\g_1,\ldots,\g_n\ket_{g,\b}=\int_{[\mmm]^\text{vir}}\prod_{i=1}^n
\ev_i^*(\g_i)\text{.}\]
Of course, it is necessary to understand exactly what this integral
means and how to compute it. This brings us to the consideration of the
mathematical aspects of the moduli space of stable maps and concludes our
excursion into the physical motivation.

While physics provided the original impetus for their introduction,
moduli spaces of stable maps 
have also been used to give answers to many problems of enumerative
geometry that do not necessarily arise in 
physics and that were inaccessible by previous methods.
Enumerative geometry seeks to determine numbers of geometric 
objects of a given type that satisfy certain conditions. The most natural
method for answering such enumerative questions consists of the following
steps.  First, construct a moduli space parameterizing the type of objects
to be counted. Second, define an intersection theory on the moduli space,
including an appropriate fundamental class against which integrals are to
be evaluated.  Third, identify subspaces of the moduli space 
(and their associated classes in
the intersection theory) corresponding
to the various conditions imposed by the enumerative problem.
Finally, integrate the product of these classes against the fundamental class.
The result can only be nonzero if 
the sum of the codimensions of these subspaces is equal to the dimension
of the fundamental class. In this case the intersection of these subspaces
(with respect to
the fundamental class) is expected to
be a finite number of points. These points
correspond to solutions of the enumerative problem. The integral computes
the number of such points (with multiplicities), and hence the number
of solutions to the enumerative problem.

We will particularly direct our attention to  the enumerative geometry of 
curves and to the second step of the above process. 
This step involves defining an intersection ring for the moduli space.
Its algebraic incarnation, which we deal with for the most part, is also 
called the Chow ring. In the topological category it is usually
referred to as the cohomology ring.
The Chow ring $A^*(X)$ of a space X 
gives an algebraic way to compute intersections
in it. Intersections correspond to multiplications in the ring; unions 
correspond to addition. See Section \ref{sec:int} for the definition
of the  relevant Chow rings.

To ensure a satisfactory theory,
computations in enumerative geometry should occur
on a compact moduli space; this
prevents solutions from ``disappearing to infinity." 
(See \cite{Katz} for some simple examples.)
One is usually interested
in the enumeration of smooth curves, but the parameter space of smooth
curves in an ambient space is not compact in any nontrivial situation.
There are various ways to compactify this space, but Kontsevich's 
compactification by the moduli space of stable maps seems especially
well-suited for enumerative purposes.

In the ten years since Kontsevich introduced the concept in \cite{KM} and
\cite{Ko}, the moduli space of stable maps has been exploited to solve
a plethora of enumerative problems for curves via the 
above process. As a rule, these results were derived without a 
complete description of the Chow rings involved. Instead, the requisite
intersection numbers were calculated somewhat indirectly,
most often using the method of {\em localization}, which will be described
in Section \ref{sec:eq}.  Such a complete description would be the key step
in giving another, more direct, computation of these enumerative numbers
and possibly many others.
Since a presentation for a ring gives
an easy way to compute all products in the ring, giving
presentations for the Chow rings of moduli spaces of stable maps
is the clear path for attaining a full and direct knowledge of their
intersection theory.
As a consequence, this would also help give a new and more direct way of
determining values of instanton numbers, Gromov-Witten invariants, and 
gravitational correlators.

So far, presentations for Chow rings of moduli spaces of stable maps have
been given only in a few special cases. Most of these have projective space as
the target of the stable maps, and in this case the moduli space 
$\mg$ depends
on four nonnegative integer parameters: the genus $g$ of the curves, the
number $n$ of marked points on the curves, the dimension $r$ of the target
projective space, and the degree $d$ of the stable maps.  
Most impressive is Mustata's presentation in \cite{M} for 
$A^*(\mm_{0,1}(\PP^r,d))$ for arbitrary $d$ and $r$. (Recently A. Mustata and M. A. Mustata
have provided an extended description of this presentation in \cite{MM}.) Behrend and
O'Halloran give a presentation for $A^*(\mm_{0,0}(\PP^r,2))$ and
conjecture a presentation for $A^*(\mm_{0,0}(\PP^r,3))$ in \cite{BO}. 
Also of relevance,  Dragos Oprea has recently described a system of
{\em tautological subrings} of the cohomology (and hence Chow) 
rings in the genus zero
case and shown that, if the target $X$ is an $\SL$
flag variety, then all rational
cohomology classes on $\mznXb$ are tautological. This gives, at least in
principle, a set of generators for any such Chow ring, namely its tautological
classes. He furthermore describes an additive basis for the cohomology ring 
of any genus zero moduli space (with target a projective algebraic variety),
which is a substantial step toward giving
a presentation. Finally, he speculates that all relations between the
tautological generators are consequences of the topological recursion 
relations. These developments may provide direction for finding presentations
for the 
Chow rings of many more moduli spaces of stable maps in the near future.
See \cite{O} and \cite{O2} for more details.
More basic 
examples include $A^*(\mm_{0,n}(\PP^r,0))\iso A^*(\PP^r)\x A^*(\MM_{0,n})$,
where $\MM_{0,n}$ is the moduli space of stable curves.  This case reduces
to finding presentations for the rings $A^*(\MM_{0,n})$, and Keel does so
in \cite{K}. Also
$\mm_{0,0}(\PP^r,1)$ is isomorphic to $\GG(1,r)$, 
the Grassmannian of lines in projective
space, and $\mm_{0,1}(\PP^r,1)$ is isomorphic to $\FF(0,1;r)$, 
the flag variety of lines in
projective space together with a point on the line. The spaces 
$\mm_{0,n}(\PP^1,1)$ are Fulton-MacPherson compactifications of configuration
spaces of $\PP^1$. Presentations for their Chow rings 
were given by Fulton and MacPherson in \cite{FM}. 
Detailed descriptions of Chow rings of spaces $\mg$, with $g>0$, are
almost nonexistent. Additional complications arise in this case.

This dissertation gives the first known presentation for a Chow ring of a 
moduli space of stable maps of degree greater than one with more than one 
marked point. 
In particular, we give the following 
presentation for $A^*(\mm_{0,2}(\PP^1,2))$:
\[A^*(\mm_{0,2}(\PP^1,2))\iso\frac{\QQ[D_0,D_1,D_2,H_1,H_2,\psi_1,\psi_2]}
{\left({H_1^2, H_2^2,D_0\psi_1,D_0\psi_2,D_2-\psi_1-\psi_2,
\psi_1-\frac{1}{4}D_1-\frac{1}{4}D_2-D_0+H_1, \atop
\psi_2-\frac{1}{4}D_1-\frac{1}{4}D_2-D_0+H_2, (D_1+D_2)^3,
D_1\psi_1\psi_2}\right)}\text{.}\]

\vspace{.1in}

\noindent
Some steps involved in finding this presentation are extended to the case
of the Chow rings $A^*(\mm_{0,2}(\PP^r,2))$,
with target arbitrary dimensional projective space. 
Presentations for these Chow rings will be included in a future paper.

Chapter \ref{sec:mod} gives the background on moduli spaces of stable
maps and their intersection theory. Chapters \ref{sec:ser} through
\ref{sec:prez} provide a detailed construction of the presentation for
$A^*(\mm_{0,2}(\PP^1,2))$ and prove that this presentation is complete.
Knowing the Betti numbers of $\mm_{0,2}(\PP^1,2)$ is the
important first step in our computation, since it will give us a good idea of
how many generators and relations to expect in each degree.  We accomplish
this in Chapter \ref{sec:ser} by using the equivariant Serre polynomial method
of \cite{GP}. In fact, in this chapter we compute the Betti numbers of
$\mm_{0,2}(\PP^r,2)$ for arbitrary $r$.
 In Chapter \ref{sec:gen} we list some natural divisor classes that occur
in the Chow rings of all moduli spaces of stable maps.  
It will become clear later that in the case of $\mm_{0,2}(\PP^1,2)$
these classes generate the Chow ring.
Relations among these classes are found in Chapter \ref{sec:rel}.
In the course of proving these relations, computations on moduli spaces
of stable maps with fewer marked points or lower degree naturally arise.
Chapter \ref{sec:simpler} describes presentations for the Chow rings of these
simpler moduli spaces.  All of these pieces are compiled to give
the whole presentation in Chapter \ref{sec:prez}. This chapter also explains
why the presentation is complete, {\em ie.}, how we know it includes enough
generators and relations to capture the entire Chow ring.
Finally, Chapter \ref{sec:app} applies the presentation to give a new
computation of the genus zero, degree two, two-pointed gravitational
correlators of $\PP^1$. Algorithms for computing theses values have previously
been developed; see \cite{KM2} and \cite{CK}, for example.

\renewcommand{\baselinestretch}{1}
\chapter{Preliminaries on moduli spaces of stable maps}
\label{sec:mod}

\section{Moduli spaces of stable maps}
\label{sec:modintro}

We will work over the field $\CC$ of complex numbers.
All schemes will be algebraic schemes over $\CC$, and we 
let $(Sch/\CC)$ denote the category of such schemes.
In this section, the definitions 
could just as well be given for the category of schemes over
any field, and all the results about moduli stacks
hold over any field as well. 
Let $\und{n}=\NN\cap[1,n]$ be the 
initial segment consisting of the first $n$ natural numbers.
The basic notions of algebraic geometry used below can be found in
\cite{H}. For information about stacks, see \cite{LM}.

\begin{Def} \label{curve}
An {\em $n$-pointed prestable curve $(C,p_1, . . . ,p_n)$
of genus $g$}
 over $\CC$ is a connected, reduced,
projective, at worst nodal curve $C$ of arithmetic genus $g$ together with 
$n$ distinct,
nonsingular marked points $p_1, . . . ,p_n$ on $C$.
\end{Def}

We will often refer to $n$-pointed prestable curves of genus $g$ as $n$-pointed,
genus $g$ curves, or simply as {\em curves}.

\begin{Def}
The {\em special points} of a curve are the marked points and the nodes.
\end{Def}

\begin{Def}\label{def:sm}
Let X be a scheme. Let $\b\in H_2(X,\ZZ)$.
A {\em stable map $(C,x_1,\ldots,x_n,f)$ 
of class $\b$} from an $n$-pointed, genus $g$ curve  
$C$ to $X$ is a morphism $f:C\rightarrow X$ such that the push-forward
$f_*([C])$ of the fundamental class is 
$\b$ and, moreover, this data satisfies the stability condition:
If $E$ is an irreducible component of $C$ on which $f$ is constant, then

\begin{enumerate}
\item If $g(E)=0$, then $E$ contains at least three special points of $C$.

\item If $g(E)=1$, then $E$ contains at least one special point of $C$.
\end{enumerate}
\end{Def}

\begin{Def}\label{def:fsm}
A {\em family $(\pi:\C\ra S,s_1,\ldots,s_n,\mu)$ 
of stable maps from  $n$-pointed, genus $g$ curves to $X$ of 
class $\b$} over
a scheme $S$ consists of a flat, proper morphism $\pi:\C\rightarrow
S$, $n$ sections $s_1,...,s_n$ of $\pi$, and a morphism $\mu:\C
\rightarrow X$ such that for every geometric point $s\in S$, the fiber
$(\C_s,s_1(s),...,s_n(s),\mu|_{\C_s})$ is a stable map from an
 $n$-pointed, genus $g$ curve to $X$ of class $\b$.
\end{Def}

\begin{Def}\label{def:mor}
A {\em morphism} from a family of stable maps
$(\pi:\C\rightarrow S,s_1,...,s_n,\mu)$ to another family 
$(\pi\pr:\mathcal{C}\pr\rightarrow T,t_1,...,t_n,\nu)$ 
is a fiber diagram

\begin{table*}[h]
  \begin{center}
\begin{equation*}    
\leavevmode
 \xymatrix{{\C} \ar[d]^{\pi} \ar[r]^{\phi}
  & {\mathcal{C}\pr} \ar[d]^{\pi\pr}\\ 
{S} \ar[r]^{f} & {T}} \end{equation*}
  \end{center}
\end{table*}
such that $\nu\circ\phi=\mu$ and 
$t_i\circ f=\phi\circ s_i$ for every $i\in \und{n}$.
\end{Def}

The definitions of an isomorphism of stable maps and of an automorphism of
a stable map are clear from the definition of a morphism.
A stable map as in Definition \ref{def:sm} is a family of stable maps over
$\Spec(\CC)$.

The above definitions concerning stable maps produce a moduli problem. 
The associated moduli functor
\linebreak[3] $\eM:(Sch/\CC)^\circ \ra(Sets)$ to 
the category of sets is given on objects by 

\renewcommand{\baselinestretch}{1}
\small\normalsize

\vspace{.1in}

\[
\eM(S)=\left\{\begin{array}{c}
\text{isomorphism classes of families of} \\ 
\text{stable maps from $n$-pointed, genus $g$ } \\
\text{curves to $X$ of class $\b$ over $S$}
\end{array}\right\}\text{.} 
\]
\small\normalsize


\noindent
Given a morphism $g:T\ra S$, $\eM(g)$ takes a family
$(\pi:C\ra S,s_1,\ldots,s_n,\mu)$ to $(\pi\pr:C\x_S T\ra T,s_1\pr,
\ldots,s_n\pr, \mu\circ g\pr)$, where $\pi\pr$, $g\pr$, and the $s_i\pr$
are the morphisms naturally induced by the fiber product.

Unfortunately, the functor $\eM$ 
is not usually representable by a scheme; there is hardly ever a
fine moduli space
for this moduli problem. There are two natural routes to pursue in search
of some other sort of moduli space, and both prove fruitful. First,
$\eM$ does have a {\em coarse} moduli scheme $\MMM$,
at least if $X$ is projective. Second, by enlarging
our ambient category from complex schemes to complex stacks $(St/\CC)$, we
can consider the moduli stack $\mmm$ of stable maps. The moduli stack
captures all the data of the moduli problem, while the moduli scheme
loses some information, including that of
automorphisms of families. Since retaining
all of this data leads to a more beautiful, powerful, and complete theory,
we will work with the stack incarnations of the moduli spaces. 

Define a category $\mmm$ over $Sch/\CC$ whose objects
are families of stable maps from  $n$-pointed, genus $g$ curves to $X$ of 
class $\b$ over complex schemes and whose morphisms are as in Definition
\ref{def:mor}.

\begin{prop}
The category $\mmm$ is a complex stack.
\end{prop}

\noindent
First, $\mmm$ is a groupoid by definition: Given a family over a scheme $S$ 
and a morphism $T\ra S$,
the fiber product as in Definition
\ref{def:mor} always exists, and it is unique up to a unique isomorphism.
Further, isomorphisms are a sheaf and every descent datum is effective.
These properties follow from Grothendieck's descent theory.  (See 
\cite{Gr} and \cite[Chapter V]{Man}.)
The following two basic properties of moduli stacks of stable maps were first
proven by Kontsevich in \cite{Ko}. More detailed proofs appear in
\cite{FP}, although the language of stacks is avoided there.

\begin{prop}
Let $X$ be a projective scheme of finite type over $\CC$. Then $\mmm$ is 
a proper Deligne-Mumford stack of finite type.
\end{prop}

For the next property, we need two definitions.

\begin{Def}
Let $\mathcal{M}_{g,n}(X,\b)$ 
be the substack of $\mmm$ corresponding to stable
maps from smooth curves.
\end{Def}

\begin{Def}
A smooth, complete variety $X$ is {\em convex} if, for any morphism
$f:\PP^1\ra X$, $H^1(\PP^1,f^*(T_X))=0$, where $T_X$ is the tangent
bundle of $X$.
\end{Def}

\noindent
The most important examples of convex varieties are 
{\em homogeneous} varieties, that is, varieties which are the quotient of
an algebraic group by a parabolic subgroup. In particular,
projective spaces are homogeneous, thus convex.
Moduli spaces of stable maps to convex varieties have many nice properties,
especially when $g=0$.

\begin{prop}
\label{nicemod}
Let $X$ be a smooth, proper, convex scheme. Then the stack \linebreak[4]
$\mznXb$ is smooth,
and the complement of $\mathcal{M}_{0,n}(X,\b)$ 
is a divisor with normal crossings.
\end{prop}

We think of $\mznXb$ as a compactification of $\mathcal{M}_{0,n}(X,\b)$;
indeed, stable maps were defined with this in mind. Besides allowing 
degenerations of the maps themselves, the compactification allows marked
points to approach nodes and each other without coinciding. In the limit,
new components of the curve appear at the place where such coincidence
would otherwise occur. We illustrate the two main types of such ``sprouting"
of new components.

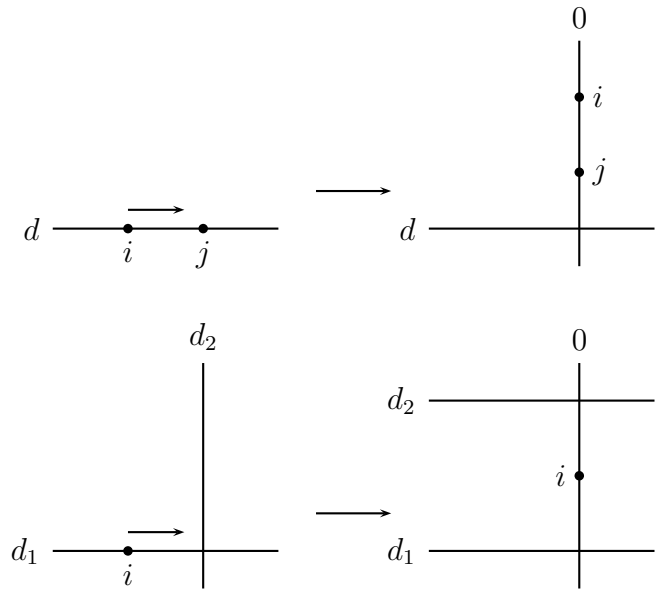
\begin{center}
\begin{pspicture}(0,0)(9,4)
\pnode(.5,1){a}
\pnode(3.5,1){b}
\dotnode(1.5,1){c}
\dotnode(2.5,1){d}
\ncline{a}{b}
\uput{5pt}[d](1.5,1){$i$}
\uput{5pt}[d](2.5,1){$j$}
\uput{5pt}[l](.5,1){$d$}
\pnode(1.5,1.25){e}
\pnode(2.25,1.25){f}
\ncline{->}{e}{f}
\pnode(4,1.5){g}
\pnode(5,1.5){h}
\ncline{->}{g}{h}
\pnode(5.5,1){i}
\pnode(8.5,1){j}
\pnode(7.5,.5){k}
\pnode(7.5,3.5){l}
\dotnode(7.5,1.75){m}
\dotnode(7.5,2.75){n}
\ncline{i}{j}
\ncline{k}{l}
\uput{5pt}[r](7.5,2.75){$i$}
\uput{5pt}[r](7.5,1.75){$j$}
\uput{5pt}[u](7.5,3.5){$0$}
\uput{5pt}[l](5.5,1){$d$}
\end{pspicture}

\vspace{0.1in}

\begin{pspicture}(0,0)(9,4)
\pnode(.5,1){a}
\pnode(3.5,1){b}
\pnode(2.5,.5){c}
\pnode(2.5,3.5){d}
\ncline{a}{b}
\ncline{c}{d}
\dotnode(1.5,1){z}
\uput{5pt}[d](1.5,1){$i$}
\uput{5pt}[l](.5,1){$d_1$}
\uput{5pt}[u](2.5,3.5){$d_2$}
\pnode(1.5,1.25){e}
\pnode(2.25,1.25){f}
\ncline{->}{e}{f}
\pnode(4,1.5){g}
\pnode(5,1.5){h}
\ncline{->}{g}{h}
\pnode(5.5,1){i}
\pnode(8.5,1){j}
\pnode(7.5,.5){k}
\pnode(7.5,3.5){l}
\pnode(8.5,3){m}
\pnode(5.5,3){n}
\ncline{i}{j}
\ncline{k}{l}
\ncline{m}{n}
\dotnode(7.5,2){y}
\uput{5pt}[l](7.5,2){$i$}
\uput{5pt}[l](5.5,1){$d_1$}
\uput{5pt}[l](5.5,3){$d_2$}
\uput{5pt}[u](7.5,3.5){$0$}
\end{pspicture}
\end{center}

\noindent
The new components arising in this way always have degree zero and contain
any marked points involved in the limit. Of course, larger numbers of marked
points can simultaneously approach each other as well. The
degenerations become more
complicated. (See \cite{FM}.) 

Behrend and Manin prove further basic properties of these moduli stacks
in \cite{BM}, including a description of the universal family 
$(\pi:\C\ra\mmm,\s_1,\ldots,\s_n,\mu:\C)$ over
the moduli stack $\mmm$. This description can be conveniently
expressed using contraction (or ``forgetful") morphisms, which we will now
introduce.
En route to doing so, we will briefly describe another class of 
moduli spaces that will play an auxiliary role in this dissertation, the
moduli spaces of stable curves.

\begin{Def}
An $n$-pointed, genus $g$ curve $C$ (Definition \ref{curve}) is {\em stable}
if, whenever
$E$ is an irreducible component of $C$, then

\begin{enumerate}
\item If $g(E)=0$, then $E$ contains at least three special points of $C$.

\item If $g(E)=1$, then $E$ contains at least one special point of $C$.
\end{enumerate}
\end{Def}

\noindent
Comparing with Definition \ref{def:sm}, we see that stable curves
correspond to stable maps to a point. The definitions of families 
of stable curves and morphisms between these families are analogous to
Definitions \ref{def:fsm} and \ref{def:mor}. Families of $n$-pointed,
genus $g$ stable curves, together with their morphisms, form a category
$\mm_{g,n}$. All of the $\mm_{g,n}$ are Deligne-Mumford stacks 
except for $\mm_{0,0}$, $\mm_{0,1}$, $\mm_{0,2}$, and $\mm_{1,0}$. Stable 
curves do not exist in these cases with the definition above, although
these moduli spaces do exist as
Artin stacks. Coarse moduli schemes $\MM_{g,n}$ also exist with the
same exceptions. In fact, for $n\geq 3$, Knudsen shows in \cite{Kn} that
$\MM_{0,n}$ is a fine moduli scheme
and a complete, nonsingular variety. 
The stack $\mm_{0,n}$ itself may be considered as
a variety in this case. We will use the notation $\MM_{0,n}$ for these moduli
spaces in recognition of these nice properties.

Knudsen defines and constructs contraction morphisms for families of
stable curves. 
Assume that $2g-2+n>0$, and let $\f$ be a morphism over $S$ 
from
a family $(\C\ra S,s_1,\ldots,s_{n+1})$ of $(n+1)$-pointed,
genus $g$ stable curves to a family $(\C\pr\ra S,s\pr_1,\ldots,s\pr_{n})$ of 
$n$-pointed,
genus $g$ stable curves.

\begin{Def}\label{def:contract}
With notation as above, let $s\in S$ be a geometric point, $E_s\sub\C_s$ the
component of the fiber over $s$ containing the $(n+1)$'st marked point 
$s_{n+1}(s)$. The morphism $\f$ is a {\em contraction} if $s_i\pr=\f s_i$ for
$i\in\und{n}$ and $\f_s:\C_s\ra\C_s\pr$ is an isomorphism except in the
following situation: If $g(E_s)=0$ and $E_s$ has only three special points,
then $E_s$ maps to a closed 
point ({\em i.e.}, is contracted), and $\f_s$ is an
isomorphism between the complements of $E_s$ and its image.
\end{Def}

Knudsen further shows that these contraction morphisms are unique.
These morphisms give rise to
contraction morphisms $\pi:\mm_{g,n+1}\ra\mm_{g,n}$ of the moduli stacks
via their universal families.
By permuting labels on the marked points and applying the above repeatedly,
we can extend the definition of contraction to morphisms that forget the 
marked points labeled by any
subset $B\sub\und{n+1}$ of the labeling set (as long as the target is not
one of the four exceptions listed earlier).

Similarly, given $A\sub\und{n}$,
there is a contraction morphism $\pi_A:\mm_{g,n}(X,\b)\ra\mm_{g,A}(X,\b)$.
(If $\b=g=0$, we require $|A|\geq 3$, and if $\b=0$ and $g=1$, we require 
$A\neq\emptyset$.) Here $\mm_{g,A}(X,\b)$ parameterizes
stable maps whose marked points are indexed by $A$, although contraction
morphisms often implicitly include a monotonic relabeling of the marked 
points of the target by $\und{|A|}$. This map is given 
pointwise by
forgetting all the marked points in the complement of $A$ and contracting
components that become unstable. 
Behrend and Manin show these are morphisms by modifying
the argument of  \cite{Kn}.
Contraction morphisms on the coarse moduli spaces are defined in the same
way.
In case only the $i$'th 
marked point is forgotten, we write
the contraction as $\pi_i$. 

The family of curves 
$\pi:\C\ra\mmm$ involved in the universal
family is given by $\pi_{n+1}:\mm_{g,n+1}(X,\b)\ra\mmm$. 
On the fiber over $(C,x_1,\ldots,x_n,f)$, 
the universal stable map $\mu:\C\ra X$ is induced by $f$ since there is
an inclusion $C\ra\mm_{g,n+1}(X,\b)$.
For $i\in\und{n}$,
the image of a stable map $(C,x_1,\ldots,x_n,f)$ under
the universal section $\s_i$ is the {\em stabilization} of the prestable
map $(C,x_1,\ldots,x_n,x_i,f)$ where the $i$'th and $n+1$'st marked points
agree. Stabilization is achieved by replacing $x_i$ with a rational component
containing the $i$'th and $(n+1)$'st marked points and mapping this new 
component to the image of $x_i$.  Stabilization is also a morphism (\cite{Kn},
\cite{BM}).

There is another forgetful morphism $\st:\mmm\ra\mm_{g,n}$ which forgets
the map data, remembering only the source curves and their marked points,
and contracts components which become unstable as a result. 
Manin calls this the {\em absolute stabilization map} and shows that it
is a morphism in \cite[Chapter V]{Man}. More detail is available in
\cite{BM}.

Composing each universal section 
with the universal map of the universal family of stable maps over the moduli
space, we get evaluation maps $\ev_1, . . . ,\ev_n$; 
$\ev_i=\mu_i\circ\s_i$. The universal map $\mu$ considered above can be
identified with $\ev_{n+1}$. Pointwise we have $\ev_i(C,x_1,\ldots,x_n,f)
=f(x_i)$. 

The last collection of morphisms we need is that of
the {\em gluing morphisms}.
Fiberwise, these involve gluing two disjoint curves together at a distinguished
marked point of each. For stable maps, the images of these marked points must
agree. (There is another kind of gluing morphism that identifies
two marked points on the same curve. However, since this type of gluing always
increases the arithmetic genus of the curve, it will not play a role in our
study.) For families, this amounts to gluing along two sections.  Proofs that
these are morphisms can be found in \cite{Kn} for stable curves and \cite{Man}
for 
stable maps. We end up with gluing morphisms
\[\mm_{g_1,n_1+1}\x\mm_{g_2,n_2+1}\ra\mm_{g_1+g_2,n_1+n_2}
\]
and
\[\mm_{g_1,n_1+1}(X,\b_1)\x_X\mm_{g_2,n_2+1}(X,\b_2)\ra
\mm_{g_1+g_2,n_1+n_2}(X,\b_1+\b_2)
\]
among moduli spaces. For the latter maps, the fiber product
is with respect to the evaluation morphisms corresponding to the two markings
being glued. The image of a gluing map is a boundary divisor whose generic
stable map has domain curve with two irreducible components, as described
in Section \ref{sec:gen} for some special cases.  The attributes of each
component are inherited from the corresponding factor in the domain of the
gluing map.

\section{Intersection theory}
\label{sec:int}

We need to define the Chow rings of the moduli stacks $\mg$ and describe
their basic properties. Fulton's book \cite{F} presents a comprehensive
introduction to intersection theory on algebraic schemes.
The requisite extensions to an analogous theory on Deligne-Mumford 
stacks were developed by Vistoli in \cite{V}. 
Manin gives a wonderfully clear exposition of Vistoli's theory in
\cite[Chapter V,\S\S 6--8]{Man}. Since Manin records much of the
information we need on Chow groups, the first subsection 
will essentially reproduce parts of his
exposition for the reader's convenience. The second subsection is for the
most part taken directly from \cite{V}.
Below we only review
the relevant definitions and properties, referring the reader to \cite{V}
and \cite{Man} for proofs and details. All stacks are assumed to be of finite
type over a fixed field $k$ unless otherwise specified.

\subsection{Chow groups}

\begin{Def}
Let $F$ be a stack. A {\em cycle of dimension $n$} on $F$ is an element
of the free abelian group $Z_n(F)$ generated by the symbols $[G]$ for all 
$n$-dimensional integral closed 
substacks $G$ of $F$ .
\end{Def}

\noindent
Here $[G]$ is the cycle associated with 
the closed substack $G$ via equivalence
groupoids.

A {\em rational function} on an integral
stack $G$ is an equivalence class of morphisms
$G\pr\ra\AA^1$, where $G\pr$ is an open dense substack of $G$, and morphisms
are equivalent if they agree on a dense open substack. The rational functions
on $G$ form a field $k(G)$. Let $W_n(F)=\+_Gk(G)$, 
the sum being taken over all integral substacks $G$
of $F$ of dimension $n+1$. If $X$ is a scheme, there is a homomorphism
\[\div_n^X:W_n(X)\rightarrow Z_n(X)\]
that takes a rational function to the cycle associated to its Weil divisor.
Thus there is a morphism $\div_n$ of presheaves on the 
\'{e}tale topology of $F$ given on an open set $X\ra F$ by $\div_n^X$.
The presheaves $W_n$ and $Z_n$ are sheaves by descent theory, and their
groups of global sections are $W_n(F)$ and $Z_n(F)$, respectively.

\begin{Def}
The {\em group $R_n(F)$ of $n$-dimensional cycles rationally equivalent
to zero} 
on a stack $F$ is the image of $\div_n$ on global sections.
\end{Def}

\begin{Def}
The {\em Chow group} of $F$ is $A_*(F)=\+_{n\geq 0}A_n(F)$, where $A_n(F)$
is the quotient $Z_n(F)/R_n(F)$.
The {\em rational Chow group} $A_*(F)_\QQ$ of $F$ is $A_*(F)\*\QQ$.
\end{Def}

The basic operations of Vistoli's intersection theory regularly introduce 
fractional coefficients to cycles and their classes in reaction to
nontrivial automorphism groups of the corresponding substacks. Thus
it only makes sense to work with the rational Chow groups $A_*(F)_\QQ$.
We will use rational Chow groups throughout this dissertation, and from
now on we write $A_*(F)_\QQ$ as simply $A_*(F)$. All cycle groups will
be tensored with $\QQ$ as well. (Kresch has developed
an integer-valued intersection theory on Deligne-Mumford stacks in \cite{Kr},
but one must still tensor with $\QQ$ in order to do enumerative geometry.)

For an integral stack $G$ of a stack $F$, let $\d(G)$ be the degree 
of the automorphism
group of a generic point of $G$. Define the {\em normalized fundamental
cycle} of $G$ to be $[G]_{\nor}=\d(G)[G]$, where $[G]$ is the fundamental
cycle associated to $G$ as above. We will use the same symbols for
the corresponding classes in the Chow group of $F$.
Thus there are two notions of the fundamental
class of a substack, and it turns out both are important. This may have first
been pointed out by Mumford in \cite{Mu}, who described the classes 
$[G]_{\nor}$ as the appropriate ones for determining rational equivalences
and the classes $[G]$ as the right ones for computing intersections. Although
these classes only differ by a rational number, we must be careful to properly
distinguish between them later. 
Otherwise the results will be confused and useless.

Let $f:F\rightarrow G$ be a separated dominant morphism of finite type of 
integral stacks. It gives an imbedding $k(G)\ra k(F)$. We define the
{\em degree} of $f$ by
\[\deg(f)=\deg(F/G)=\frac{\d(G)}{\d(F)}[k(F):k(G)]\text{.}\]

Some basic operations of intersection theory are proper pushforward, flat
pullback, and Gysin maps. We now give the definitions of the first two
in Vistoli's theory. We will sketch the construction of Gysin maps and list
some of their properties in the next subsection.

\begin{Def}
Let $f:F\rightarrow G$ be a morphism of stacks.
\begin{enumerate}
\item If $f$ is flat (and of constant relative dimension), 
we define the flat pullback 
\[f^*:Z_*(G)\rightarrow Z_*(F)\]
by $f^*[G\pr]=[G\pr\x_G F]$ for any closed integral substack $G\pr$ of $G$.

\item If $F$ is proper, the proper pushforward
\[f_*:Z_*(F)\rightarrow Z_*(G)\]
is defined by $f_*[F\pr]=\deg(F\pr/G\pr)[G\pr]$, where $F\pr$ is an integral
substack of $F$ and $G\pr$ is its image in $G$.

\end{enumerate}
\end{Def}

\begin{prop}
The flat pullback and proper pushforward defined above pass to rational
equivalence.
\end{prop}

Thus we get flat pullback $f^*:A_*(G)\rightarrow A_*(F)$ and proper 
pushforward $f_*:A_*(F)\rightarrow A_*(G)$.

\subsection{Cones, local regular embeddings, and Gysin maps}

Constructing Gysin maps for Chow groups of stacks requires additional theory,
which we will only sketch. The most important players are cones and
regular local embeddings. We quote the relevant material from \cite{V}.

Cones can be constructed by descent: Let

\vspace{0.1in}

\renewcommand{\baselinestretch}{1}\small\normalsize
\psset{arrows=->}
\begin{center}
\begin{psmatrix}
$U\x_F U$ & $U$
\ncline[offset=-3pt]{1,1}{1,2}\nbput{$p_1$}
\ncline[offset=3pt]{1,1}{1,2}\naput{$p_2$}
\end{psmatrix}
\end{center}
\psset{arrows=-}

\vspace{0.1in}

\noindent
be a presentation of $F$. If $C_U$ is a cone on $U$, and there is
an isomorphism of cones $p_1^*C_U\iso p_2^*C_U$ satisfying the cocycle
condition, then $C_U$ is the pullback of a canonically defined cone
on $F$.

\begin{lem}
\label{locemb}
Let $f:F\ra G$ be a representable morphism of finite type of stacks.
Then $f$ is unramified if and only if there are atlases $U\ra F$ and 
$V\ra G$ together with an embedding $U\ra V$ compatible with $f$.
\end{lem}

For this reason, 
a representable, unramified morphism of finite type 
of stacks is called a {\em local embedding}.

\begin{Def}
A local embedding of stacks $f:F\rightarrow G$ is {\em regular of 
codimension $d$} if we can choose  $U$,
$V$ and the local embedding $f\pr:U\rightarrow V$ as in the previous lemma
such that $f\pr$ is a regular embedding of codimension $d$.
\end{Def}

\begin{eg}
If $F$ is a smooth stack over a scheme $S$ of constant relative dimension $d$,
then the diagonal $F\ra F\x_S F$ is a regular local embedding of 
codimension $d$.
\end{eg}

In the situation of Lemma \ref{locemb}, 
consider the normal cone $C_{U/V}$ to $U$ in $V$. One can
show that there is an isomorphism from $p_1^*C_{U/V}$ to $p_2^*C_{U/V}$
that satisfies the cocycle condition.

\begin{Def}
The cone $C_{F/G}$ obtained by descent from $C_{U/V}$ is called the 
{\em normal cone} to $F$ in $G$. If $f$ is a regular local embedding,
$C_{F/G}$ is a vector bundle called the normal bundle and written 
$N_{F/G}$.
\end{Def}

Gysin maps can be defined in a manner quite similar to that used
for schemes in \cite{F}.
Suppose we have a fiber square of stacks

\vspace{0.2in}

\renewcommand{\baselinestretch}{1}\small\normalsize
\psset{arrows=->}
\begin{equation}\label{square}
\text{
\begin{psmatrix}
$F\pr$ & $G\pr$ \\
$F$ & $G$
\ncline{1,1}{2,1}\naput{$p$}
\ncline{1,1}{1,2}\naput{$g$}
\ncline{2,1}{2,2}\naput{$f$}
\ncline{1,2}{2,2}\naput{$q$}
\end{psmatrix}
}\end{equation}
\psset{arrows=-}

\vspace{0.2in}

\noindent
with $f$ a local regular imbedding of codimension $d$. Vistoli assumes
that $G\pr$ is a scheme for his construction.  We can reduce to
the case where $G\pr$ is a scheme by taking atlases.
Thus we can use Vistoli's approach.

First assume $G\pr=T$ is a purely $k$-dimensional scheme, so that $F\pr=U$
is as well. Then $g$ is a local imbedding of schemes. There is a natural
closed imbedding $C_{U/T}\ra N$, where $N=p^*N_{F/G}$ is the pullback of the
normal bundle of $F$ in $G$. Let $s:U\ra N$ be the zero section, and let
$s^*:A_*(N)\ra A_*(U)$ be the Gysin homomorphism defined in \cite[\S 3.3]{F}
as the inverse to the isomorphism induced by the vector bundle projection.
Define $(F.T)\in A_{k-d}(U)$ to be $(F.T)=s^*[C_{U/T}]$.

Now suppose $G\pr$ is an arbitrary scheme. Let 
$y\pr=\sum_im_i[T_i]\in A_*(G\pr)$, where
the $T_i$ are subvarieties of $G\pr$. Let $h_i:F\x_G T_i\ra F\pr$ be the
natural imbeddings.

\begin{Def}
The {\em Gysin map} $f^!:Z_*(G\pr)\ra A_*(F\pr)$ is a homomorphism defined by
\[f^!y\pr=\sum_im_i{h_i}_*(F.T_i)\text{.}\]
\end{Def}

\noindent
Gysin maps pass to rational equivalence, and we call the resulting
homomorphisms $f^!:A_*(G\pr)\ra A_*(F\pr)$ Gysin maps as well.
These statements are also true if $G\pr$, and thus $F\pr$, are allowed to 
be stacks, as we assume for the rest of the subsection. 

\begin{prop}\label{Gys}
Gysin maps satisfy the following properties:

\begin{enumerate}
\item {\bf (Compatibility with proper pushforwards)} 
If $p$ is proper and $f$ is a local regular imbedding 
in the fiber diagram (\ref{Gys1}) below, 
then $f^!p_*=q_*f^!$.
\item {\bf (Compatibility with flat pullbacks)} 
If $p$ is flat and $f$ is a local regular imbedding
in the fiber diagram (\ref{Gys1}) below, 
then $f^!p^*=q^*f^!$.
\item {\bf (Commutativity)} If $f$ and $j$ are local regular imbeddings 
in the fiber diagram (\ref{Gys2}) below, then $f^!j^!=j^!f^!$.
\end{enumerate}

\vspace{0.2in}

\renewcommand{\baselinestretch}{1}\small\normalsize
\psset{arrows=->}
\begin{equation}\label{Gys1}
\text{
\begin{psmatrix}
$F\dpr$ & $G\dpr$ \\
$F\pr$ & $G\pr$ \\
$F$ & $G$
\ncline{1,1}{2,1}\naput{$q$}
\ncline{1,1}{1,2}\naput{$f\dpr$}
\ncline{2,1}{2,2}\naput{$f\pr$}
\ncline{1,2}{2,2}\naput{$p$}
\ncline{2,1}{3,1}
\ncline{2,2}{3,2}\naput{$g$}
\ncline{3,1}{3,2}\naput{$f$}
\end{psmatrix}
}\end{equation}
\psset{arrows=-}

\vspace{0.2in}

\psset{arrows=->}
\begin{equation}\label{Gys2}
\text{
\begin{psmatrix}
$F\dpr$ & $G\dpr$ & $X$ \\
$F\pr$ & $G\pr$ & $Y$ \\
$F$ & $G$
\ncline{1,1}{2,1}
\ncline{1,1}{1,2}
\ncline{2,1}{2,2}
\ncline{1,2}{2,2}
\ncline{2,1}{3,1}
\ncline{2,2}{3,2}
\ncline{3,1}{3,2}\naput{$f$}
\ncline{1,2}{1,3}
\ncline{2,2}{2,3}\naput{$g$}
\ncline{1,3}{2,3}\naput{$j$}
\end{psmatrix}
}\end{equation}
\psset{arrows=-}
\renewcommand{\baselinestretch}{2}\small\normalsize

\begin{flushright}
$\Box$
\end{flushright}
\end{prop}
Gysin maps also satisfy 
other properties, but they will not be relevant for us.
In case $G\pr=G$ and $q=\id_G$ in Diagram (\ref{square}), 
we use the notation $f^*:A_*(G)\ra A_*(F)$
for the Gysin map.

\subsection{Intersection product}

If $F$ is a smooth stack of dimension $n$, 
then the diagonal imbedding $\d:F\ra F\x F$ is
a local regular imbedding of codimension $n$. Define a product morphism 
$A_k(F)\* A_l(F)\ra A_{k+l-n}(F)$ by
\[x\cdot y=\d^*(x\x y)\text{.}\]
Set $A^k(F)=A_{n-k}(F)$ and $A^*(F)=\sum_k A^k(F)$. Then this product
makes $A^*(F)$ into a commutative and graded ring with identity element $[F]$,
called the {\em Chow ring of $F$}.

\subsection{Additional intersection theory facts}

In this subsection, we mention some other intersection theoretic objects
and ideas that we will need, including the excision sequence, Chern classes,
expected dimension and virtual fundamental classes of moduli spaces, and
homology isomorphism.

\begin{prop}[Excision sequence]
Let $i:G\ra F$ be a closed substack of a stack $F$ with complement $j:U\ra F$.
Then the sequence

\psset{arrows=->}
\begin{equation}\label{eksiz}
\text{
\begin{psmatrix}
$A_k(G)$ & $A_k(F)$  & $A_k(U)$ & 0
\ncline{1,1}{1,2}\naput{$i_*$}
\ncline{1,2}{1,3}\naput{$j^*$}
\ncline{1,3}{1,4}
\end{psmatrix}
}
\end{equation}
\psset{arrows=-}

\noindent
of Chow groups is exact for all $k$.
\end{prop}

The excision sequence is a useful tool for learning about
the structure of the Chow group of a space by decomposing it into simpler
spaces.

Chern classes of vector bundles on stacks can be defined almost exactly as
done in \cite{F} for schemes. We will take the definition 
from \cite{Man}. Let $L$ be an invertible sheaf on a stack $F$ of
dimension $n$. Then 
$c_1(L)\cap [F]\in A_{n-1}$ 
is defined to be the Cartier divisor associated to $L$. More generally, for
a closed substack $f:G\ra F$, set 
\[c_1(L)\cap [G]=f_*(c_1(f^*(L))\cap G)\text{.}\]
Thus $c_1(L)\cap$ is an operator on $A_*(F)$. If $E\ra F$ 
is a vector bundle
of rank $e+1$, let $p:P\ra F$ be its projectivization with tautological
bundle $\O_P(1)$. First, we define Segre classes 
$s_i(E)\cap$ by
\[s_i(E)\cap y=p_*(c_1(\O_P(1))^{e+i}\cap p^*(y))\]
for $y\in A_*(F)$. The Segre polynomial of $E$ is $s_t(E)=\sum_i s_i(E)t^i$.
The Chern polynomial $c_t(E)$ is defined similarly; its coefficients are
the Chern classes of $E$. We define the Chern classes by the formula
$c_t(E)=(s_t(E))^{-1}$. These classes satisfy the usual properties of Chern
classes. We will often slightly abuse terminology and notation by 
calling the cycle class $c_i(E)\cap[F]$ the $i$'th Chern class of $E$ and
writing it as just $c_i(E)$.

The moduli spaces $\mmm$ have an {\em expected dimension} from deformation
theory, given by the formula
\[(\dim(X)-3)(1-g)-\int_{\b}K_X + n\text{,}\]
where $K_X$ is the canonical class of $X$. If $X$ is convex, then the spaces
$\mznXb$ always have the expected dimension. (See \cite{FP}.)
For $\m0$ this dimension is $d+r+dr+n-3$.
All of our work will occur
in this pleasant situation. However, in general the dimension of $\mmm$
can be larger than its expected dimension.  This is analogous to the 
situation where two subvarieties of a variety do not intersect properly.
In this case the ordinary fundamental class $[\mmm]$ is not the correct
class to integrate against (nor is $[\mmm]_{\nor}$) 
in order to compute gravitational 
correlators or do enumerative geometry in general. We must introduce a third
type of fundamental class $[\mmm]^{\text{vir}}$ 
called the {\em virtual fundamental
class}, which lives in the Chow group of the expected dimension
and gives enumerative
geometry and  
gravitational correlators the desired properties. Construction 
of virtual fundamental classes
is very complicated. See \cite[\S 7.1.4]{CK} for an overview and further
references. We mention the virtual fundamental class only because it appears
in some of our general definitions and formulas.

Let $H^*(F)$ denote the rational de Rham cohomology ring of a Deligne-Mumford
stack $F$.

\begin{prop}[Homology isomorphism]
\label{hi}
Let $X$ be a flag variety. Then
there is a ring isomorphism

\begin{equation}\label{a2h}
A^*(\mznXb)\ra H^*(\mznXb)
\end{equation}
\end{prop}

\noindent
See \cite{O} for a proof.
In particular, this holds for $X=\PP^r$. It allows us to switch freely
between cohomology and Chow rings. We should note that the isomorphism
doubles degrees. The degree of a $k$-cycle in $A^k(\mznXb)$, called the
algebraic degree, is half the degree of its image in $H^{2k}(\mznXb)$. 

\renewcommand{\baselinestretch}{1}
\section{Equivariant Cohomology and Localization on $\m0$}
\label{sec:eq}

Localization often greatly simplifies gravitational correlator and 
other enumerative geometry calculations. The calculations we are
interested in take place in intersection rings of moduli
stacks $\mmm$.  We will concentrate on the relatively simple case where
computations occur in the intersection rings of spaces $\m0$.

\subsection{Equivariant cohomology and the localization 
theorem of Atiyah and Bott}

Equivariant cohomology was originally constructed in the topological setting
by Atiyah and Bott (\cite{AB}). More recently, the corresponding algebraic
theory of equivariant Chow rings was developed by Edidin and Graham 
(\cite{EG}). It follows from Proposition \ref{hi} that the cycle map of 
\cite{EG}
gives an isomorphism between the equivariant Chow ring and the equivariant
cohomology ring of a moduli space of stable maps, 
so we can switch freely between these settings. Our description 
below follows the exposition of \cite{CK}, which takes a topological
perspective. 

Let $X$ be a topological space 
and $G$ a connected Lie group with classifying bundle $EG\ra BG$.
Setting $X_G=X\x_{G}EG$, we define the {\em equivariant cohomology} of $X$ to
be $H_G^*(X)=H^*(X_G)$. 
We now state some basic facts about equivariant cohomology.
First, $H_G^*(\point)=H^*(BG)$.
By pullback via $X\ra\point$, $H_G^*(X)$ is an $H^*(BG)$-module. We can 
regard $H^*(BG)$ as the coefficient ring for equivariant cohomology.
Note that inclusion of a fiber $i_X:X\ra X_G$ induces a ``forgetful
map'' $i_X^*:H_G^*(X)\ra H^*(X)$.

We will consider the case where $G$ is the torus $T=(\CC^*)^n$. Let 
$M(T)$ be the character group of $T$. For each $\r\in M(T)$, we get a 
1-dimensional vector space $\CC_\r$ with a $T$-action given by $\r$. If 
$L_\r=(\CC_\r)_T$ is the corresponding line bundle over $BT$, then the
assignment $\r\mapsto -c_1(L_\r)$ defines an isomorphism
$\psi:M(T)\ra H^2(BG)$, which induces a ring isomorphism 
$\Sym(M(T))\iso H^*(BG)$.  We call $\psi(\r)$ the {\em weight} of $\r$.

Let $\r_i\in M(T)$ be the character given by the $i$'th projection, and let
$\la_i$ be the weight of $\r_i$. Then
\[H^*(BG)\iso \CC[\la_1,\ldots,\la_n]\text{.}\]
The map $i_X^*$ can be thought of as a
nonequivariant limit that maps all the $\la_i$ to 0.
We denote the line bundle $L_{\r_i}$ by $\O(-\la_i)$, so that 
$\la_i=c_1(\O(\la_i))$. 

Let $T$ act on a smooth manifold $X$. The fixed point locus $X^T$ is a
union of smooth connected components $Z_j$. We have inclusions
$i_j:Z_j\ra X$ and normal bundles $N_j=N_{Z_j/X}$ which are equivariant.
Inclusion induces $i_j^*:H_T^*(X)\ra H_T^*(Z_j)$. Since
$Z_j$ is a submanifold of 
$X$, we also have a Gysin map ${i_j}_!:H_T^*(Z_j)\ra H_T^*(X)$.

If $E$ is a $G$-equivariant vector bundle of rank $r$ on $X$, the top
Chern class $\text{Euler}_G(E)=c_r^G(E)$ is called the equivariant Euler
class of $E$.
Let $\R_T\iso\CC(\la_1,\ldots,\la_n)$ be the field of fractions of 
$H^*(BT)$.  Atiyah and Bott \cite{AB} have shown that $\eqeul(N_j)$ is 
invertible in the localization $H_T^*(Z_j)\*\R_T$ for all $j$.

\begin{thm}[Localization Theorem of Atiyah-Bott]
There is an isomorphism
\[H_T^*(X)\*\R_T\iso\+_j H_T^*(Z_j)\*\R_T\]
induced by the map $\a\mapsto(i_j^*(\a)/\eqeul(N_j))_j$. The inverse is
induced by $(\a_j)_j\mapsto\sum_j{i_j}_!(\a_j)$. Thus for any 
$\a\in H_T^*(X)\*\R_T$ we have
\begin{equation}\label{locexp}
\a=\sum_j {i_j}_!\left(\frac{i_j^*(\a)}{\eqeul(N_j)}\right)\text{.}
\end{equation}
\end{thm}

\noindent
{\bf Idea of Proof:} The self-intersection formula says
$i_j^*\circ{i_j}_!(\g)=\g\cup\eqeul(N_j)$ for any 
$\g\in H_T^*(Z_j)$. Let $\g=i_j^*(\a)$.$\Box$

For any variety $X$ with $T$-action, $X\ra\point$ induces an equivariant
projection $\pi_X:X_T\ra BT$. The pushforward map ${\pi_X}_!$ will be called
the equivariant integral and written
\[\int_{X_T}:H_T^*(X)\ra H^*(BT)\text{.}\]

\begin{cor}
For any $\a\in H_T^*(X)\*\R_T$,
\[\int_{X_T}\a=\sum_j \int_{(Z_j)_T}\frac{i_j^*(\a)}{\eqeul(N_j)}\text{.}\]
\end{cor}

\noindent
{\bf Proof.} Apply ${\pi_X}_!$ to both sides of (\ref{locexp}).$\Box$

A more topological or analytical approach to moduli problems leads to the 
notion of an {\em orbifold}.
Orbifolds correspond to smooth Deligne-Mumford stacks. 
A variety $X$ is an orbifold if it
admits local (analytic) charts $U/H$ with $U$ smooth and $H$ a small
subgroup of $\GL(n,\CC)$ acting on $U$.
Deligne-Mumford stacks always admit stratification into quotient substacks.
Viewing these quotients as varieties rather than stacks, 
it is not hard to imagine that one can get an orbifold
associated to the stack if the stack is smooth. 
A $T$-action on a stack gives rise to local
$T$-actions on the charts $U$.
By working with the charts, the technical issues 
that would arise in considering localization
on smooth stacks directly can be avoided.
We will always take this approach, accounting for the local quotients at the
end by dividing answers by the order of the quotienting group.
The resulting formula is otherwise the same as that for varieties.

\begin{cor}
\label{stackloc}
Let $X$ be an orbifold which is the variety underlying a smooth
stack with a $T$-action. 
If $\a\in H_T^*(X)\*\R_T$, then
\[\int_{X_T}\a=\sum_j \int_{(Z_j)_T}\frac{i_j^*(\a)}{a_j\eqeul(N_j)}\text{,}\]
where $a_j$ is the order of the group $H$ occurring in a local chart at the
generic point of $Z_j$.
\end{cor}

\subsection{Localization in $\m0$}

The natural action of 
$T=(\CC^*)^{r+1}$ on $\PP^r$ induces a $T$-action
on $\m0$ by composition of the action with stable maps. 
This moduli space is a smooth Deligne-Mumford stack, so
we can apply Corollary \ref{stackloc}. 

The fixed point loci and their equivariant normal bundles were determined
for $\m0$ by Kontsevich \cite{Ko}. A $T$-fixed point of $\m0$ corresponds
to a stable map $(C, p_1,\ldots, p_n, f)$ where each component of $C_i$ of $C$
is either mapped to a $T$-fixed point of $\PP^r$ or multiply covers a 
coordinate line.  Also, each marked point $p_i$, each node of $C$, and each
ramification point of $f$ is mapped to a $T$-fixed point of $\PP^r$. 
This implies that coordinates on $C_i$ and its image can be chosen so the
cover is given by $(x_0,x_1)\mapsto(x_0^{d_i},x_1^{d_i})$. As a result,
the fixed point components $Z_j$ of the $T$-action can be described by 
combinatorial data. Let $q_0,\ldots, q_r$ be the fixed points of $\PP^r$ 
under this $T$-action,
so that $q_i=(0:\ldots:0:1:0\ldots:0)$, with the 1 in the $i$'th position.
The connected components of $\m0^T$ are in 1--1 
correspondence with connected trees $\G$ of the following type:
The vertices $v$ of $\G$ are in 1--1 correspondence with the connected 
components $C_v$ of $f^{-1}(\{q_0,\ldots,q_r\})$, so each $C_v$ is either a point
or a connected union of irreducible components of $C$. The edges
$e$ of $\G$ correspond to irreducible components $C_e$ of $C$ which are
mapped onto some coordinate line $\ell_e$ in $\PP^r$.

The graph $\G$ has the following labels: Associate to each vertex $v$
the number $i_v$ defined by $f(C_v)=q_{i_v}$, as well as the set $S_v$
consisting of those $i$ for which the marked point $p_i$ is in $C_v$.
Associate to each edge $e$ the degree $d_e$ of the map $f|_{C_e}$.
Finally, we impose the following three conditions:

\begin{enumerate}
\item If an edge $e$ connects $v$ and $v\pr$, then $i_v\neq i_{v\pr}$, and
$\ell_e$ is the coordinate line joining $q_{i_v}$ and $q_{i_{v\pr}}$.

\item $\sum_e d_e=d$.

\item $\coprod_v S_v=\und{n}$.
\end{enumerate}

The locus of stable maps with graph $\G$ is a fixed point substack $\mm_\G$ of
$\m0$. Fix $(C, p_1,\ldots, p_n, f)\in\mm_\G$. For each $v$ such that $C_v$
is a curve, $C_v$ has $n(v)=|S_v|+\val(v)$ special points. The data consisting
of $C_v$ plus
these $n(v)$ points forms a stable curve, giving an element of $\MM_{0,n(v)}$.
Using the data of $\G$, we can construct a morphism
\[\psi_\G:\prod_{v:\dim C_v=1} \MM_{0,n(v)}\ra \mm_\G\text{.}\]
Define $\MM_\G$ to be the above product (which is a point if all components
are contracted). The morphism $\psi_\G$ is finite.  Indeed, there
is a finite group of automorphisms $A_\G$ acting on $\MM_\G$ such that
$\mm_\G$ is the quotient stack $[\MM_\G/A_\G]$. We have an exact sequence
\[0\lra\prod_e \ZZ/d_e\ZZ\lra A_\G\lra\Aut(\G)\lra 0\text{,}\]
where $\Aut(\G)$ is the group of automorphisms of $\G$ which preserve
the labels. The left-hand term comes from sheet interchanges of multiply
covering components.
When Corollary \ref{stackloc} is applied to $\m0$, the factor
$a_{\G}$ appearing in the denominator of the term corresponding to $\mm_\G$
is the order of $A_\G$.

The last ingredients needed in order to use localization on $\m0$ are the
Euler classes of the fixed components. Denote the normal bundle of
$\mm_\G$ by $N_\G$. Define a {\em flag} $F$ of a graph to be a pair
$(v,e)$ such that $v$ is a vertex of $e$. Put $i(F)=v$ and let $j(F)$ be
the other vertex of $e$. Set
\[\w_F=\frac{\la_{i_{i(F)}}-\la_{i_{j(F)}}}{d_e}\text{.}\]
This corresponds to the weight of the $T$-action on the tangent space
of the component $C_e$ of $C$ at the point $p_F$ lying over $i_v$. 
%
%
Let
$e_F$ be the first Chern class of the bundle on $\mm_\G$ whose fiber
is the cotangent space to the component associated to $v$ at $p_F$.
(More information about this type of class is given in Section \ref{sec:psi}.)
If $\val(v)=1$, let $F(v)$ denote
the unique flag containing $v$.  If $\val(v)=2$, let $F_1(v)$ and $F_2(v)$
denote the two flags containing $v$. Similarly, let $v_1(e)$ and $v_2(e)$
be the two vertices of an edge $e$.

\begin{thm}\label{norm}
The equivariant Euler class of the normal bundle $N_\G$ is a product
of contributions from the flags, vertices and edges:
\[\eqeul(N_\G)=e_{\G}^{\text{F}}e_{\G}^{\text{v}}e_{\G}^{\text{e}}\text{,}\]
where 
\[e_{\G}^{\text{F}}=\frac{\prod_{F:n(i(F))\geq 3}(\w_F-e_F)}
{\prod_{F}\prod_{j\neq i_{i(F)}}(\la_{i_{i(F)}}-\la_j)}\]
\[e_{\G}^{\text{v}}=\left(\prod_v \prod_{j\neq i_v}(\la_{i_v}-\la_j)
\right)
\left(\prod_{{\val(v)=2\atop S_v=\emptyset}}(\w_{F_1(v)}+\w_{F_2(v)})\right)
/\prod_{{\val(v)=1\atop S_v=\emptyset}} \w_{F(v)}\]
\[e_{\G}^{\text{e}}=\prod_e\left(\frac{(-1)^{d_e}(d_e!)^2
(\la_{i_{v_1(e)}}-\la_{i_{v_2(e)}})^{2d_e}}{d_e^{2d_e}}
\prod_{{a+b=d_e\atop k\neq i_{v_j(e)}}}
\left(\frac{a\la_{i_{v_1(e)}}+b\la_{i_{v_2(e)}}}{d_e}-\la_k\right)\right)
\text{.}
\]
\end{thm}

This formula won't be very efficient in most particular examples 
because there will be many cancellations. However, it is sufficient for
our purposes. We are now equipped with the theory necessary to compute
integrals over moduli spaces of stable maps using localization. This will
be undertaken in Chapters \ref{sec:simpler} and \ref{sec:rel}.

\chapter{The Betti numbers of $\mm_{0,2}(\PP^r,2)$}
\label{sec:ser}

\renewcommand{\baselinestretch}{1}
\section{Serre polynomials and the Poincar\'{e} polynomial 
of $\mm_{0,2}(\PP^r,2)$}
\label{sec:poincare}

This section owes much to Getzler and Pandharipande, who provide the 
framework for computing the Betti numbers of all the spaces $\m0$
in their unpublished work 
\cite{GP}. (They have recently completed these computations in \cite{GP2}.)
However, we will take the definitions and basic results from
other sources, and prove their theorem in the special case that we need. 
We will compute a formula for the Poincar\'{e} polynomials of the moduli
spaces $\mm_{0,2}(\PP^r,2)$ using what are called Serre polynomials in
\cite{GP} and \cite{Ge}. (These polynomials are also known as virtual
Poincar\'{e} polynomials or E-polynomials.)
Serre polynomials are defined for 
quasi-projective varieties over $\CC$ via the mixed Hodge theory of Deligne
(\cite{D}). Serre conjectured the existence of polynomials satisfying
the first two key properties given below. A formula was later given by
Danilov and Khovanski\u{\i} in \cite{DK}. If $(V,F,W)$ is a mixed Hodge
structure over $\CC$, set
\[V^{p,q}=F^p\gr_{p+q}^WV\cap \bar{F}^q\gr_{p+q}^WV
\]
and let $\X(V)$ be the Euler characteristic of $V$ as a graded vector space.
Then 
\[\ser(X)=\sum_{p,q=0}^\infty u^pv^q\X(H_c^\bullet(X,\CC)^{p,q})\text{.}
\]
If $X$ is a smooth projective variety, then the 
Serre polynomial of $X$ is just its Hodge polynomial:
\[\ser(X)=\sum_{p,q=0}^\infty(-u)^p(-v)^q\dim H^{p,q}(X,\CC)\text{.}\]
If $X$ further satisfies $H^{p,q}(X,\CC)=0$ for $p\neq q$, then we can 
substitute a new variable $q=uv$ for $u$ and $v$.
In this case, the coefficients of the Serre polynomial of $X$ give its
Betti numbers, so that $\ser(X)$ is the Poincar\'{e} polynomial of $X$.  

We will use two additional key properties of Serre 
polynomials from \cite{Ge}. The first gives a compatibility
with decomposition: If $Z$ is a closed subvariety of $X$, then 
$\ser(X)=\ser(X\backslash Z)+\ser(X)$.  Second, it respects products:
$\ser(X\x Y)=\ser(X)\ser(Y)$.  (This is actually a consequence of the previous
properties.) It follows from these two properties 
that the Serre polynomial
of a fiber space is the product of the Serre polynomials of the base and the
fiber. The definition and properties above come from \cite{Ge}.
We also use the following consequence of the Eilenberg-Moore spectral
sequence, which is essentially Corollary 4.4 in \cite{Sm}.

\begin{prop}
Let $Y\ra B$ be a fiber space with $B$ simply connected, and let $X\ra B$
be continuous. If $H^*(Y)$ is a free $H^*(B)$-module, then
\[H^*(X\x_B Y)\iso H^*(X)\*_{H^*(B)}H^*(Y)\]
as an algebra.
\end{prop}

Since we deal exclusively with cases where the isomorphism (\ref{a2h})
holds, there is never any torsion in the cohomology. Thus we have the
following.

\renewcommand{\baselinestretch}{1}
\begin{cor}
\label{serfiber}
Let $X$ and $Y$ be varieties over a simply connected base $B$, and suppose
either $X$ or $Y$ is locally trivial over $B$. Then
\[\ser(X\x_B Y) = \frac{\ser(X)\ser(Y)}{\ser(B)}\text{.}\]
\end{cor}

\noindent
We will sometimes use the notation $Y/B$ for the fiber of a fiber space
$Y\ra B$.

To extend this setup to Deligne-Mumford stacks, where automorphism
groups can be nontrivial (but still finite), 
{\em equivariant} Serre polynomials 
are needed.  Let $G$ be a finite group acting on a (quasiprojective) 
variety $X$.
The idea is this: The action  of $G$ on $X$
induces an action on its cohomology (preserving the mixed Hodge 
structure), which in turn gives a representation of $G$ on each (bi)graded
piece of the cohomology. The cohomology of the quotient variety $X/G$, 
and hence of
the quotient stack $[X/G]$, is the part of the cohomology of $X$ which
is fixed by the $G$-action, {\em i.e.}, in each degree the subspace on which
the representation is trivial. 

Our definition comes from 
\cite{Ge}. 
The equivariant Serre polynomial $\ser(X,G)$ of $X$ is given by the formula
\[\ser_g(X)=\sum_{p,q=0}^\infty u^pv^q
\sum_i(-1)^i\Trace(g|(H_c^i(X,\CC))^{p,q})\text{.}
\]
for each element $g\in G$.
We can also describe the equivariant Serre polynomial more compactly
with the formula
\[\ser(X,G)=\sum_{p,q=0}^\infty u^pv^q
\sum_i(-1)^i[H_c^i(X,\CC)^{p,q}]\text{,}\]
taken from \cite{GP}.
In the
case $G=S_n$, we write $\ser_n(X)$ for $\ser(X,S_n)$.
A $G$-equivariant Serre polynomial takes values
in $R(G)[u,v]$, where $R(G)$ is the virtual representation ring of $G$.
The augmentation morphism $\epsilon:R(G)\rightarrow\ZZ$, which extracts the 
coefficient of the trivial representation $\one$ from an element
of $R(G)$, extends to an
augmentation morphism $R(G)[u,v]\rightarrow\ZZ[u,v]$. If $G$ acts on 
a quasi-projective variety $X$, the Serre polynomial of the
quotient stack $[X/G]$ is the 
augmentation of the equivariant Serre polynomial of $X$.

Every virtual representation ring $R(G)$ 
has the extra structure
of a {\em $\la$-ring}. Our definition of a $\la$-ring comes from \cite{Knutson}.
First, let $\xi_1,\ldots,\xi_q,\n_1,\ldots,\n_r$ be variables, and let
$s_i$ and $\s_i$ be the $i$'th elementary symmetric functions in the $\xi_j$'s
and the $\n_j$'s, respectively. Define $P_n(s_1,\ldots,s_n,\s_1,\ldots\s_n)$
be the coefficient of $t^n$ in
\[\prod_{i,j}(1+\xi_i\n_jt)
\]
and $P_{n,d}(s_1,\ldots,s_{nd})$ to be the coefficient of $t^n$ in
\[\prod_{1\leq i_1<\ldots<i_d\leq q}(1+\xi_{i_1}\cdots\xi_{i_d}t)
\]
\begin{Def}\label{lring}
A {\em $\la$-ring} is a commutative ring $R$ with a series of operations
$\la_k:R\ra R$ for $k\in\{0\}\cup\NN$ satisfying the following properties.
\begin{enumerate}
\item For all $x\in R$, $\la_0(x)=1$. 

\item For all $x\in R$, $\la_1(x)=x$.

\item \label{lamsum}
For all $x,y\in R$, $\la_n(x+y)=\sum_{k=0}^n \la_k(x)\la_{n-k}(y)$.

\item $\la_t(1)=1+t$.

\item For all $x,y\in R$ and $n\in\{0\}\cup\NN$,
$\la_n(xy)=P_n(\la_1x,\la_2x,\ldots\la_nx,\la_1y,\ldots,\la_ny)$.

\item  For all $x\in R$ and $n,m\in\{0\}\cup\NN$, 
$\la_m(\la_n(x))=P_{m,n}(\la_1x,\ldots,\la_{mn}x)$.

\end{enumerate}
\end{Def}

\noindent
Here $\la_t(x)$ is the formal power series $\sum_{i=0}^\infty \la_i(x)t^i$.
If the $\la_i$ satisfy the first three properties, $R$ is called a
{\em pre-$\la$-ring}. 

Let $V$ be a $G$-module. Then $\la_i(V)$ is the {\em $i$'th exterior power}
$\La^iV$
of $V$, where we define $g\in G$ to act by 
$g(v_1\wedge\ldots\wedge v_i)=gv_1\wedge\ldots\wedge gv_i$. Define $\la_0(V)$
to be the trivial one-dimensional representation. 
(One can similarly define a $G$-module structure
on the $i$'th symmetric power $S^iV$.) Knutson proves in 
\cite[Chapter II]{Knutson} that
these exterior power operations give $R(G)$ the structure of a $\la$-ring
for any finite group $G$. Addition is given by $[V]+[W]=[V\+W]$, and the
product is $[V]\cdot[W]=[V\* W]$, both with the naturally induced actions.

Knutson also shows that $\ZZ$ is a $\la$-ring with $\la$-operations
given via $\la_t(m)=(1+t)^m$. For $m,n\geq 0$, this definition gives 
$\la_n(m)=\lchoose{m}{n}$. Finally, he shows that
if $R$ is a $\la$-ring, then there is
a unique structure of $\la$-ring on $R[x]$ under which 
$\la_k(rX^n)=\la_k(r)X^{nk}$ for $n,k\in\NN\cup\{0\}$ and $r\in R$.
This gives a $\la$-ring structure on $\ZZ[q]$. The augmentation morphism 
$\e:R(G)\rightarrow\ZZ$ is a {\em map of $\la$-rings}; it commutes with the
$\la$-operations.

We will use the following facts about Serre polynomials and equivariant
Serre polynomials.
For $n\in\NN$, let $[n]=\frac{q^n-1}{q-1}$. Then $[n+1]$ is the Serre
polynomial of $\PP^n$, as is clear from the presentation for its Chow ring.  
Getzler and Pandharipande prove that the 
Serre polynomial of the Grassmannian $G(k,n)$ of 
$k$-planes in $\CC^n$ is the $q$-binomial coefficient
\[\chews{n}{k}=\frac{[n]!}{[k]![n-k]!}\text{,}\]
where $[n]!=[n][n-1]\cdot\cdot\cdot[2][1]$. We will prove this formula in
the special case $k=2$.

\begin{lem}\label{grasser}
The Serre polynomial of $G(2,n)$ is $\chews{n}{2}$.
\end{lem}

\noindent
{\bf Proof.} 
By projectivizing the ambient affine space, $G(2,n)\iso\GG(1,n-1)$, the
Grassmannian of lines in $\PP^{n-1}$. We will use the projective viewpoint
in this proof.  The universal $\PP^1$ bundle over $\GG(1,n-1)$ is isomorphic
to $\FF(0,1;n-1)$ of pairs $(p,\ell)$ of a point $p$ and a line $\ell$ in
$\PP^{n-1}$ with $p\in\ell$. On the other hand, there is a projection
$\FF(0,1;n-1)\ra\PP^{n-1}$ taking $(p,\ell)$ to $p$. Its fiber over a point
$p$ is $\{\ell\, |\, p\in\ell\}$, 
which is isomorphic to $\PP^{n-2}$. (To see this
isomorphism, fix a hyperplane $H\sub\PP^{n-1}$ not containing $p$ and map
each line to its intersection with $H$.) It follows that $\ser(\FF(0,1;n-1))
=[n][n-1]$. Since $\ser(\FF(0,1;n-1))=\ser(\GG(1,n-1))[2]$ also, 
we are able to 
conclude that $\ser(\GG(1,n-1))=[n][n-1]/[2]$.
$\Box$

Next, since $\PGL(2)$
is the complement of a quadric surface in $\PP^3$, 
$\ser(\PGL(2))=[4]-[2]^2=q^3-q$.

In addition to the $\la$-operations, every $\la$-ring $R$
has {\em $\s$-operations}
as well. These can be defined
in terms of the $\la$-operations by $\s_k(x)=(-1)^k\la_k(-x)$. 
Routine checking shows that the $\s$-operations give $R$ the structure of
a pre-$\la$-ring. 
Here
we simply note the following formulas for the $\la$-ring $\ZZ[q]$.
\[\s_k([n])=\chews{n+k-1}{k} \text{ and } 
\la_k([n])=q^{k \choose 2}\chews{n}{k}\text{.}\]
Proofs of these formulas can be found in \cite[Section I.2]{Mac}.
Next we explain why these formulas are relevant. Let $\e$ be the sign
representation of $S_n$. Note the identity $\e^2=\one$.
We will prove the following claim from \cite{GP}.

\begin{lem}
If $X$ is smooth and $S_2$ acts on $X^2$ by switching the factors, then
\[\ser_2(X^2)=\s_2(\ser(X))\text{\em \one}+\la_2(\ser(X))\e\text{.}\]
\end{lem}

\noindent
{\bf Proof.}
Let $V$ be a vector space. Now $V\* V=S^2V\+\La^2V$ as $S_2$-modules,
with $S_2$ acting by switching the factors of $V\* V$,
trivially on $S^2V$, and by sign on $\La^2V$.
If 0 is the zero representation, certainly $\la_i(0)=0$ for $i>0$.
We use this and the properties of $\la$-rings in Definition \ref{lring}
to obtain
\begin{eqnarray*}
0 & = & \la_2([V]-[V]) \\
  & = & \one\cdot\la_2(-[V])+[V]\cdot(-[V])+\la_2[V]\cdot\one\\
  & = & \la_2(-[V])-[S^2V]-[\La^2V]+\la_2[V]
\end{eqnarray*}

Since $\s_2[V]=\la_2(-[V])$, this implies $\s_2[V]=[S^2V]$. Since $X$ is 
smooth $H^*(X^2)=H^*(X)\* H^*(X)$, with the action of $S_2$ switching the
factors. Applying the above with $V=H^*(X)$ gives 
$[H^*(X^2)]=\s_2[H^*(X)]+\la_2[H^*(X)]$. Breaking this down by 
(cohomological) degree,
we have $[H^i(X^2)]q^i=[\s_2[H^*(X)]]_iq^i+[\la_2[H^*(X)]]_iq^i$. 
We need to show
that $[H^i(X^2)]q^i=[\s_2(\ser(X))]_i\one +[\la_2(\ser(X))]_i\e$. We will show
the equality of
the first summands of each expression; showing equality of
the terms involving $\la_2$ is easier. First, by induction the identity
$\la_2(-m)=\la_2(m+1)$ holds. Second, note that any pre-$\la$-operation $\la_2$
acts on sums by $\la_2(\sum_i x_i)=\sum_i \la_2(x_i)+\sum_{i<j}x_ix_j$.
Third, note that vector spaces in the following computation live in the
graded algebra $H^*(X)\* H^*(X)$, and we will apply the usual rules for
grading in a tensor product. Finally, all of the representations below
are trivial.
We find 
\begin{eqnarray*}
&   & [\s_2[H^*(X)]]_iq^i \\
& = & \left[\s_2\left[\sum_j H^j(X)\right]\right]_iq^i\\
& = & \left[\sum_j[S^2H^j(X)]+\sum_{j<k}[H^j(X)\* H^k(X)]\right]_iq^i\\
& = & \begin{cases}
\left([S^2H^{i/2}(X)]+\sum_{\stackrel{\scriptstyle j+k=i}{j<k}}
[H^j(X)\* H^k(X)]\right)q^i
& \text{if $i$ is even,}\\
\left(\sum_{\stackrel{\scriptstyle j+k=i}{j<k}}[H^j(X)\* H^k(X)]\right)q^i
& \text{if $i$ is odd,}
\end{cases}\\
& = & \begin{cases}
\left(\lchoose{h^{i/2}(X)+1}{2}\one
+\sum_{\stackrel{\scriptstyle j+k=i}{j<k}}h^j(X)h^k(X)\one\right)q^i
& \text{if $i$ is even,}\\
\left(\sum_{\stackrel{\scriptstyle j+k=i}{j<k}}h^j(X)h^k(X)\one\right)q^i
& \text{if $i$ is odd.}
\end{cases}
\end{eqnarray*}
On the other hand,
\begin{eqnarray*}
&   & [\s_2(\ser(X))]_i\one \\
& = & [\la_2(-\sum h^j(X)q^j)]_i\one \\
& = & \left[\sum \la_2(-h^j(X))q^{2j}+\sum_{j<k}h^j(X)h^k(X)q^{j+k}
	\right]_i\one \\
& = & \begin{cases}
\left(\lchoose{h^{i/2}(X)+1}{2}q^i
+\sum_{\stackrel{\scriptstyle j+k=i}{j<k}}h^j(X)h^k(X)q^i\right)\one
& \text{if $i$ is even,}\\
\left(\sum_{\stackrel{\scriptstyle j+k=i}{j<k}}h^j(X)h^k(X)q^i\right)\one
& \text{if $i$ is odd.}
\end{cases}
\end{eqnarray*}

\begin{flushright}
$\Box$
\end{flushright}

\noindent
As a corollary, the ordinary Serre polynomial of $[X^2/S_2]$ is 
$\s_2(\ser(X))$.

A key fact used in our computations is the following.
\begin{prop}
\label{Serre00}
If $d>0$, $\ser(\nobarmm_{0,0}(\PP^r,d))=q^{(d-1)(r+1)}\chews{r+1}{2}$.
\end{prop}

\noindent
This follows from Pandharipande's proof in \cite{P4} that the
Chow ring of the nonlinear Grassmannian $M_{\PP^k}(\PP^r,d)$ is
isomorphic to the Chow ring of the ordinary Grassmannian $\GG(k,r)$.
If $k=1$, the nonlinear Grassmannian is $\nobarmm_{0,0}(\PP^r,d)$. (The Serre
polynomial grades by dimension rather than codimension. This is why the
shifting factor $q^{(d-1)(r+1)}$ appears.)
Notice that the proposition refers to the locus of stable maps with 
smooth domain curve, which is a proper (dense) subset of the compactified
moduli space $\mm_{0,0}(\PP^r,d)$.  

Recall that 
$\nobarmm_{0,n}(\PP^r,0)\iso M_{0,n}\x\PP^r$, so that the Serre polynomials
of these spaces are easy to compute.

Finally,
let $F(X,n)$ be the configuration space of $n$ distinct labeled points in 
a nonsingular variety $X$. Fulton and MacPherson show in \cite{FM} that 
\[\ser(F(X,n))=\prod_{i=0}^{n-1}(\ser(X)-i)
\text{.}\]

In order to compute
the Serre polynomial of a moduli space of stable maps,  we can stratify
it according to the degeneration types of the maps and compute the 
Serre polynomial of each stratum separately.  The degeneration types 
of maps in $\m0$ 
are in 1--1 correspondence with stable $(n,d)$-trees via taking the
dual graph of a stable map.  These concepts were defined by Behrend
and Manin. A summary of their definitions follows; see 
\cite{BM} for a full development.

\begin{Def}
A {\em graph} $\tau$ is a quadruple $(V_{\tau},F_{\tau},j_{\tau},
\partial_{\tau})$, where $V_{\tau}$ and $F_{\tau}$ are sets, 
$j_{\tau}:F_{\tau}\ra F_{\tau}$ is an involution, and 
$\partial_{\tau}:F_\t\ra V_\t$. Elements of $V_{\tau}$ are the {\em vertices}
of $\t$, and elements of $F_{\tau}$ are its {\em flags}. 
\end{Def}

Two additional sets associated to a graph $\t$ are the set of {\em tails}
$S_\t=\{f\in F_\t|j_\t(f)=f\}$ and the set of {\em edges}
$E_\t=\{(f_1,f_2)|f_i\in F_\t,j_\t(f_1)=f_2\}$.

Geometrically, a flag can be thought of as half of an edge. 
Every flag $f$ has an associated vertex $\dd_{\t}(f)$ 
attached to one of its ends. 
The edges in the graph are given by gluing two flags $f_1$ and $f_2$ 
together at their open ends whenever $j_\t(f_1)=f_2$.  The fixed points
of $j_\t$ remain half-edges, or tails, with a vertex at only one end.
This interpretation leads to the {\em geometric realization} $|\t|$ of $\t$,
which is a topological graph.

\begin{Def}
A {\em tree} is a graph $\t$ whose geometric realization is simply connected.
Equivalently, $|\t|$ is connected and $|E_\t|=|V_\t|-1$.
\end{Def}

\begin{Def}
Let $n$ and $d$ be non-negative integers. An {\em $(n,d)$-tree} is a tree
$\t$ together with a bijection $\nu:S_\t\ra\und{n}$ and a map
$d:V_\t\ra\NN\cup\{0\}$ such that
\[\sum_{v\in V_\t}d(v)=d\text{.}\]
We call $d(v)$ the {\em degree} of $v$.
\end{Def}

The {\em valence} of a vertex $v$ is $n(v)=|\{f\in F_\t|\dd(f)=v\}|$.
An $(n,d)$-tree $\t$ is {\em stable} if whenever $v$ 
is a vertex of $\t$ with $d(v)=0$, then the valence of $v$ is at
least three.
An automorphism of a stable $(n,d)$-tree is a graph automorphism that fixes
the tails and preserves the degree labels.

The {\em dual graph} of a stable map $(C,x_1,\ldots,x_n,f)$ 
has one vertex for each irreducible
component of $C$. An edge connects two vertices 
whenever the corresponding components intersect. This includes the
possibility of a loop edge, with both ends incident to the same vertex, 
if the corresponding component intersects itself. However,
this won't happen for stable maps from genus zero curves. The flags incident
to a vertex are in 1--1 correspondence with the marked points lying on that
component. Each vertex is labeled with the homology class of the pushforward
of the fundamental class of the corresponding component. 
If $X=\PP^r$, this becomes the degree of a vertex, which is just the degree 
of the restriction of $f$
to the corresponding component. It follows that
 the dual graph of a genus zero, 
$n$-pointed stable map to $\PP^r$ of degree $d$ is a stable $(n,d)$-tree.

A family of stable maps is said to {\em degenerate} over a point if the curve
in the fiber  over that point
has more nodes than the curve over the general fiber. Thus in the genus
zero case, the curve of the degenerate fiber will have extra components.
Two stable maps are 
said to have the same degeneration type if and only if their
dual graphs are identical. 

We are now ready to compute the Poincar\'{e} polynomials of some moduli spaces 
of stable maps.

\begin{prop}
The Poincar\'{e} polynomial of $\mm_{0,2}(\PP^r,2)$ is
\[\ser(\mm_{0,2}(\PP^r,2))= 
\left(\sum_{i=0}^r q^i\right)\left(\sum_{i=0}^{r-1} q^i\right)
\left(\sum_{i=0}^{r+2} q^i+2\sum_{i=1}^{r+1} q^i+2\sum_{i=2}^{r} q^i\right)
\text{,}\]
and the Euler characteristic of $\mm_{0,2}(\PP^r,2)$ is $r(r+1)(5r+3)$.
\end{prop}

\noindent
{\bf Proof.}
We begin by stratifying $\mm_{0,2}(\PP^r,2)$ according to the degeneration
type of the stable maps. Since the strata are locally closed,
the compatibility of Serre polynomials
with decomposition allows us to compute the Serre polynomial of each
stratum separately and add up the results to obtain 
$\ser(\mm_{0,2}(\PP^r,2))$.

Each stratum is isomorphic to a finite group quotient of a
fiber product of moduli spaces $\nobarm0$ via an inverse procedure
to the gluing defined in Section \ref{sec:modintro}.
Given a stable map $(C,x_1,x_2,f)$, consider 
the normalization of $C$.
It consists of a disjoint union of smooth curves $C_i$ corresponding to the
components of $C$, and there are maps $f_i$
from each curve to $\PP^r$ naturally
induced by $f$.  Furthermore, auxiliary marked points are added to retain
data about the node locations. The result is a collection of stable maps with
smooth domain curves,
one for each component. The evaluations of auxiliary marked points 
corresponding to the same node must agree. This gives rise to
a fiber product of
moduli spaces of stable maps from smooth domain curves, together with a
morphism onto the stratum coming from the normalization map.
There can be automorphisms of the stable maps in the stratum that are not
accounted for by the fiber product. These occur when there is a collection
of connected unions $U_i$ of connected components that all intersect a common
component or point and satisfy the following conditions:

\begin{enumerate}
\item None of the $U_i$ contain any marked points.
 
\item The restrictions of $f$ to $U_i$ and $U_j$ give isomorphic stable maps
for any $i$ and $j$.
\end{enumerate}

\noindent
These automorphisms correspond exactly to the automorphisms of the dual graph
$\G$ of the stratum.
Proving that $\Aut(\G)$ is the right group to quotient by appears quite
complicated in general. We prove it for the strata of $\mm_{0,2}(\PP^r,2)$
by brute force.  All the strata are listed below, and the assertion clearly
holds in each case.
So we can compute the Serre polynomials of the strata
using Corollary \ref{serfiber} and Proposition \ref{Serre00}.

When stratified according to the dual graphs of stable maps, 
$\mm_{0,2}(\PP^r,2)$ has 9 types of strata. The corresponding graphs are

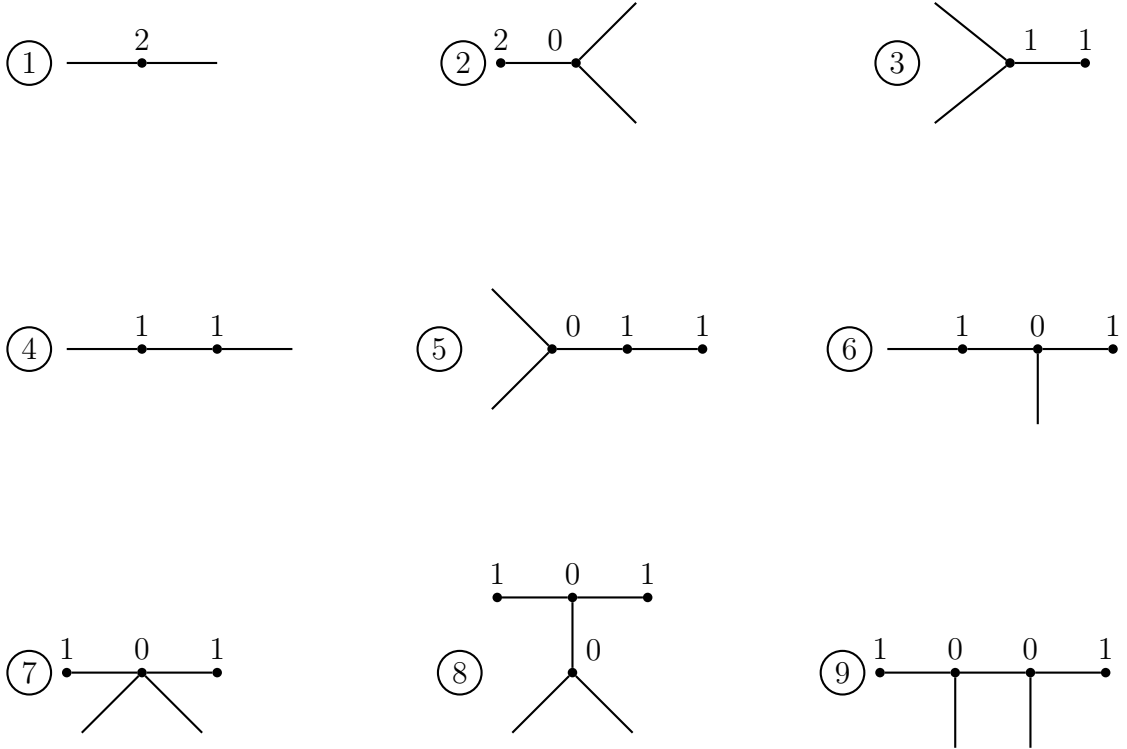
\vspace{0.5in}

\begin{pspicture}(0,0)(3,2)
\rput(0,1.5){\pscirclebox{1}}
\pnode(0.5,1.5){a}
\dotnode(1.5,1.5){c}
\pnode(2.5,1.5){b}
\ncline{a}{c}
\ncline{c}{b}
\uput{5pt}[u](1.5,1.5){2}
\end{pspicture}
\hspace{1in}
\begin{pspicture}(0,0)(3,2.5)
\rput(0,1.5){\pscirclebox{2}}
\dotnode(0.5,1.5){a}
\dotnode(1.5,1.5){b}
\pnode(2.3,2.3){c}
\pnode(2.3,0.7){d}
\ncline{a}{b}
\ncline{b}{c}
\ncline{b}{d}
\uput{5pt}[u](0.5,1.5){2}
\uput{7pt}[ul](1.5,1.5){0}
\end{pspicture}
\hspace{1in}
\begin{pspicture}(0,0)(3,2.5)
\rput(0,1.5){\pscirclebox{3}}
\pnode(0.5,0.7){c}
\dotnode(1.5,1.5){a}
\pnode(.5,2.3){d}
\dotnode(2.5,1.5){b}
\ncline{c}{a}
\ncline{a}{d}
\ncline{a}{b}
\uput{7pt}[ur](1.5,1.5){1}
\uput{5pt}[u](2.5,1.5){1}
\end{pspicture}

\vspace{0.5in}

\begin{pspicture}(0,0)(4,1)
\rput(0,1.5){\pscirclebox{4}}
\pnode(0.5,1.5){c}
\dotnode(1.5,1.5){a}
\dotnode(2.5,1.5){b}
\pnode(3.5,1.5){d}
\ncline{c}{a}
\ncline{a}{b}
\ncline{b}{d}
\uput{5pt}[u](1.5,1.5){1}
\uput{5pt}[u](2.5,1.5){1}
\end{pspicture}
\hspace{.5in}
\begin{pspicture}(0,0)(4,2.5)
\rput(0,1.5){\pscirclebox{5}}
\pnode(.7,.7){c}
\dotnode(1.5,1.5){a}
\pnode(.7,2.3){d}
\dotnode(2.5,1.5){b}
\dotnode(3.5,1.5){e}
\ncline{c}{a}
\ncline{a}{d}
\ncline{a}{b}
\ncline{b}{e}
\uput{7pt}[ur](1.5,1.5){0}
\uput{5pt}[u](2.5,1.5){1}
\uput{5pt}[u](3.5,1.5){1}
\end{pspicture}
\hspace{.5in}
\begin{pspicture}(0,0)(4,2)
\rput(0,1.5){\pscirclebox{6}}
\pnode(.5,1.5){c}
\dotnode(1.5,1.5){a}
\dotnode(2.5,1.5){b}
\dotnode(3.5,1.5){d}
\pnode(2.5,.5){e}
\ncline{c}{a}
\ncline{a}{b}
\ncline{b}{d}
\ncline{b}{e}
\uput{5pt}[u](1.5,1.5){1}
\uput{5pt}[u](2.5,1.5){0}
\uput{5pt}[u](3.5,1.5){1}
\end{pspicture}

\vspace{.5in}

\begin{pspicture}(0,0)(3,2)
\rput(0,1.5){\pscirclebox{7}}
\pnode(.7,.7){c}
\dotnode(.5,1.5){a}
\dotnode(1.5,1.5){b}
\dotnode(2.5,1.5){e}
\pnode(2.3,.7){d}
\ncline{a}{b}
\ncline{b}{e}
\ncline{b}{c}
\ncline{b}{d}
\uput{5pt}[u](0.5,1.5){1}
\uput{5pt}[u](1.5,1.5){0}
\uput{5pt}[u](2.5,1.5){1}
\end{pspicture}
\hspace{1in}
\begin{pspicture}(0,0)(3,3)
\rput(0,1.5){\pscirclebox{8}}
\pnode(.7,.7){c}
\dotnode(.5,2.5){a}
\dotnode(1.5,2.5){b}
\dotnode(2.5,2.5){e}
\dotnode(1.5,1.5){f}
\pnode(2.3,.7){d}
\ncline{a}{b}
\ncline{b}{e}
\ncline{b}{f}
\ncline{f}{d}
\ncline{f}{c}
\uput{5pt}[u](0.5,2.5){1}
\uput{5pt}[u](1.5,2.5){0}
\uput{5pt}[u](2.5,2.5){1}
\uput{7pt}[ur](1.5,1.5){0}
\end{pspicture}
\hspace{.75in}
\begin{pspicture}(0,0)(4,2)
\rput(0,1.5){\pscirclebox{9}}
\pnode(1.5,.5){c}
\dotnode(.5,1.5){a}
\dotnode(1.5,1.5){b}
\dotnode(2.5,1.5){e}
\dotnode(3.5,1.5){f}
\pnode(2.5,.5){d}
\ncline{a}{b}
\ncline{b}{e}
\ncline{e}{f}
\ncline{b}{c}
\ncline{e}{d}
\uput{5pt}[u](.5,1.5){1}
\uput{5pt}[u](1.5,1.5){0}
\uput{5pt}[u](2.5,1.5){0}
\uput{5pt}[u](3.5,1.5){1}
\end{pspicture}

There are actually 10 strata, because there are two strata of type 6 depending
on which marked point is identified with which tail.  
We use the same numbers to label the strata as 
those labeling the corresponding graphs above.
Eight of the strata have no automorphisms, so we can directly compute ordinary
Serre polynomials in these cases. The strata corresponding to
Graphs 7 and 8 have automorphism group $S_2$.  
Calculating the $S_2$-equivariant Serre
polynomials of these strata is necessary as
an intermediate step.
We now compute the Serre polynomials of the strata.

Stratum 1 is $\nobarmm_{0,2}(\PP^r,2)$. It is an $F(\PP^1,2)$-bundle over
$\nobarmm_{0,0}(\PP^r,2)$. Thus Stratum 1 has 
Serre polynomial 
\begin{eqnarray*}
\ser(F(\PP^1,2))\ser(\nobarmm_{0,0}(\PPr,2)) 
& = & (q^2+q)q^{r+1}\frac{[r+1][r]}{[2]}\\
& = & q^{r+2}[r+1][r] 
\text{.}
\end{eqnarray*}

Stratum 2 is isomorphic to the fiber product
\[\nobarmm_{0,1}(\PPr,2)\x_{\PPr}\nobarmm_{0,3}(\PPr,0)\text{.}\]
Now $\nobarmm_{0,3}(\PPr,0)\iso\PPr$, so the Serre polynomial of this stratum
is just
\[\ser(\nobarmm_{0,1}(\PPr,2))=\ser(\PP^1)\ser(\nobarmm_{0,0}(\PPr,2))
=q^{r+1}[r+1][r]\]
since $\nobarmm_{0,1}(\PPr,2)$ is a $\PP^1$-bundle over $\nobarmm_{0,0}(\PPr,2)$.

Stratum 3 is isomorphic to the fiber product
\[\nobarmm_{0,3}(\PP^r,1)\x_{\PP^r}\nobarmm_{0,1}(\PP^r,1)
\text{.}\]
The $F(\PP^1,3)$-bundle $\nobarmm_{0,3}(\PP^r,1)$ over $\nobarmm_{0,0}(\PP^r,1)$ has
Serre polynomial $(q^3-q)\chews{r+1}{2}$. Similarly, 
$\ser(\nobarmm_{0,1}(\PP^r,1))=(q+1)\chews{r+1}{2}=[r+1][r]$. Thus
Stratum 3
has Serre polynomial 
\[\frac{(q^3-q)\chews{r+1}{2}[r+1][r]}{[r+1]}
=(q^2-q)[r+1][r]^2\text{.}\]

Stratum 4 is isomorphic to the fiber product
\[\nobarmm_{0,2}(\PP^r,1)\x_{\PP^r}\nobarmm_{0,2}(\PP^r,1)
\text{.}\]
The $F(\PP^1,2)$-bundle $\nobarmm_{0,2}(\PP^r,1)$  over $\nobarmm_{0,0}(\PP^r,1)$
has Serre polynomial $(q^2+q)\chews{r+1}{2}$. Thus Stratum 4
has Serre polynomial 
\[\frac{(q^2+q)^2[r+1]^2[r]^2}{[r+1][2]^2}=q^2[r+1][r]^2\text{.}\]

Stratum 5 is isomorphic to the fiber product
\[\nobarmm_{0,3}(\PP^r,0)\x_{\PP^r}\nobarmm_{0,2}(\PP^r,1)\x_{\PP^r}\nobarmm_{0,1}(\PP^r,1)
\text{,}\]
and this in turn is isomorphic to 
$\nobarmm_{0,2}(\PP^r,1)\x_{\PP^r}\nobarmm_{0,1}(\PP^r,1)$. So Stratum 5
has Serre polynomial 
\[\frac{(q^2+q)\chews{r+1}{2}(q+1)\chews{r+1}{2}}{[r+1]}
=q[r+1][r]^2\text{.}\]

A stratum of type 6 is isomorphic to the fiber product
\[\nobarmm_{0,2}(\PP^r,1)\x_{\PP^r}\nobarmm_{0,3}(\PP^r,0)\x_{\PP^r}\nobarmm_{0,1}(\PP^r,1)
\text{.}\]
This is isomorphic to  $\nobarmm_{0,2}(\PP^r,1)\x_{\PP^r}\nobarmm_{0,1}(\PP^r,1)$, so 
each stratum of type 6 has Serre polynomial 
\[\frac{(q^2+q)\chews{r+1}{2}(q+1)\chews{r+1}{2}}{[r+1]}
=q[r+1][r]^2\text{.}\]
Thus the total 
contribution from strata of type 6 is 
\[2q[r+1][r]^2\text{.}\]

Stratum 9 is isomorphic to the fiber product
\[\nobarmm_{0,1}(\PP^r,1)\x_{\PP^r}\nobarmm_{0,3}(\PP^r,0)\x_{\PP^r}
\nobarmm_{0,3}(\PP^r,0)\x_{\PP^r}\nobarmm_{0,1}(\PP^r,1)
\text{.}\]
It
has Serre polynomial 
\[\frac{(q+1)^2\chews{r+1}{2}^2}{[r+1]}
=[r+1][r]^2\text{.}\]

We now turn our attention to the two strata with automorphisms.
Stratum 8 is isomorphic to the quotient of
\[X=\nobarmm_{0,3}(\PP^r,0)\times_{\PP^r}\nobarmm_{0,3}(\PP^r,0)\times_{(\PP^r)^2}
\nobarmm_{0,1}(\PP^r,1)^2\]
by the action of $S_2$.  The first copy of $\nobarmm_{0,3}(\PP^r,0)$
is superfluous. The action of 
$S_2$ on the cohomology of the second copy of 
$\nobarmm_{0,3}(\PP^r,0)$ is trivial. The action
switches the two factors of $\nobarmm_{0,1}(\PP^r,1)$ as well as the two factors
in $\PP^r\times\PP^r$. Since $\nobarmm_{0,1}(\PPr,1)$ is a fiber space
over $\PP^r$, we can use Corollary \ref{serfiber} in
computing the equivariant Serre polynomial of $X$ to be
\begin{eqnarray*}
& & \ser_2(\nobarmm_{0,3}(\PPr,0))\ser_2((\nobarmm_{0,1}(\PPr,1)/\PPr)^2)\\
& = & [r+1]\left(\s_2\left(\frac{\ser(\nobarmm_{0,1}(\PPr,1))}{\ser(\PPr)}\right)\one
+\la_2\left(\frac{\ser(\nobarmm_{0,1}(\PPr,1))}{\ser(\PPr)}\right)\e\right) \\
& = & [r+1](\s_2([r])\one+\la_2([r])\e) \\
& = & [r+1](\chews{r+1}{2}\one+q\chews{r}{2}\e) 
\end{eqnarray*}
As in the proof of Lemma \ref{grasser}, the fiber 
$\nobarmm_{0,1}(\PPr,1)/\PPr$ is isomorphic to $\PP^{r-1}$.
Now augmentation gives
\[\frac{[r+1]^2[r]}{[2]}\]
as the Serre polynomial of Stratum 8.
Stratum 7 is isomorphic to the quotient of 
\[Y=\nobarmm_{0,1}(\PP^r,1)^2\x_{(\PP^r)^2}\nobarmm_{0,4}(\PP^r,0)\]
by the action of $S_2$, which again switches the squared factors. In addition,
it switches two of the four marked points in $\nobarmm_{0,4}(\PP^r,0)$. 
Now $\nobarmm_{0,4}(\PP^r,0)\iso M_{0,4}\x\PP^r$, and $S_2$ acts trivially
on the $\PP^r$ factor. Furthermore, 
$M_{0,4}\iso\PP^1\setminus\{0,1,\infty\}$
has Serre polynomial $q-2$. But we need to know $\ser_2(M_{0,4})$
under an $S_2$-action switching 
two of the deleted points. 
It is not hard to imagine that
$\ser_2(M_{0,4})=(q-1)\one-\e$, but this takes some work to prove. 
Considering $M_{0,4}$ as the parameter space of four distinct points
in $\PP^1$ modulo automorphisms of $\PP^1$, we obtain
$M_{0,4}\iso F(\PP^1,4)/\PGL(2)$. Now $\PGL(2)$ acts freely on $F(\PP^1,4)$.
As a result,
\begin{equation}\label{mm04}
\ser_2(M_{0,4})=\frac{\ser_2(F(\PP^1,4))}{\ser_2(\PGL(2))}\text{.}
\end{equation}
Since the cohomology of $\PGL(2)$ is not affected by the action,
\begin{equation}\label{pgl2}
\ser_2(\PGL(2))=\ser(\PGL(2))=q^3-q\text{.}
\end{equation}

We can stratify $(\PP^1)^4$ into fifteen cells whose closures are 
respectively $(\PP^1)^4$,
the six large diagonals, the seven ``medium diagonals" where two coordinate
identifications are made, and the small diagonal, so that $F(\PP^1,4)$
is the complement of the union of all the cells corresponding to diagonals. 
We examine how the action
affects cells of each type, subtracting the polynomials for cells that are
removed.
For concreteness, suppose the first two marked points are switched. Any
pairs of Chow ring generators of $(\PP^1)^4$ which differ exactly by a 
multiple of $H_2-H_1$ are switched, so it is not hard to get
\[\ser_2((\PP^1)^4)=(q^4+3q^3+4q^2+3q+1)\one+(q^3+2q^2+q)\e\text{.}\]
How does the action affect the diagonals removed from $(\PP^1)^4$?
Exactly two pairs, $(\D_{13},\D_{23})$ and $(\D_{14},\D_{24})$, of 
the six large diagonals are switched, so the
corresponding cells contribute
\[(-4q^3+4q)\one+(-2q^3+2q)\e\]
to the equivariant Serre polynomial, since
these diagonals have been removed.
Exactly two pairs, $(\D_{134},\D_{234})$ and $(\D_{(13)(24)},\D_{(14)(23)})$,
among the seven diagonals with two identifications are 
switched as well. The corresponding cells contribute
\[(-5q^2-5q)\one+(-2q^2-2q)\e\]
to the equivariant Serre polynomial. 
The small diagonal is not affected by the action, so it contributes
\[(-q-1)\one\text{.}\]
Putting these together gives 
\[\ser_2(F(\PP^1,4))=(q^4-q^3-q^2+q)\one+(-q^3+q)\e\text{.}\]
Then by (\ref{mm04}) and (\ref{pgl2}), we have the desired result
$\ser_2(M_{0,4})=(q-1)\one-\e$. Using Corollary \ref{serfiber} again,
we thus calculate the equivariant Serre polynomial of $Y$ to be
\begin{eqnarray*}
&  & \ser_2(\nobarmm_{0,4}(\PPr,0))\ser_2((\nobarmm_{0,1}(\PPr,1)/\PPr)^2) \\
& = & [r+1]((q-1)\one-\e)\left(\chews{r+1}{2}\one+q\chews{r}{2}\e\right) \\
& = & [r+1]\left(\left((q-1)\chews{r+1}{2}-q\chews{r}{2}\right)\one
+\left((q^2-q)\chews{r}{2}-\chews{r+1}{2}\right)\right)\e\text{.}
\end{eqnarray*}
Augmentation gives
\begin{eqnarray*}
[r+1]\left((q-1)\chews{r+1}{2}-q\chews{r}{2}\right)
& = & \frac{[r+1][r]}{[2]}((q-1)[r+1]-q[r-1])\\
& = & \frac{[r+1][r]}{[2]}(q^{r+1}+q^r)-\frac{[r+1]^2[r]}{[2]}\\
& = & [r+1][r]q^r-\frac{[r+1]^2[r]}{[2]}\\
\end{eqnarray*}
as the Serre polynomial of Stratum 7. 

To get the Serre polynomial for the whole moduli space, we add together
the contributions from all the strata.  

\begin{eqnarray}
& & \ser(\mm_{0,2}(\PP^r,2))\nonumber \\
& = & q^{r+2}[r+1][r]+q^{r+1}[r+1][r]+(q^2-q)[r+1][r]^2
+q^2[r+1][r]^2+q[r+1][r]^2\nonumber \\
& & +2q[r+1][r]^2+[r+1][r]^2+\frac{[r+1]^2[r]}{[2]}
+[r+1][r]q^r-\frac{[r+1]^2[r]}{[2]}\nonumber \\
& = & [r+1][r](q^{r+2}+q^{r+1}+(q^2-q)[r]+q^2[r]+3q[r]+[r]+q^r)\nonumber \\
& = & [r+1][r](q^{r+2}+q^{r+1}+q^r+[r](2q^2+2q+1))\nonumber \\
& = & [r+1][r](q^{r+2}+q^{r+1}+q^r+2\sum_{i=2}^{r+1}q^i+2\sum_{i=1}^{r}q^i
+\sum_{i=0}^{r-1}q^i)\nonumber \\
& = & \left(\sum_{i=0}^r q^i\right)\left(\sum_{i=0}^{r-1} q^i\right)
\left(\sum_{i=0}^{r+2} q^i+2\sum_{i=1}^{r+1} q^i+2\sum_{i=2}^{r} q^i\right)
\text{.}
\label{ser2r2}
\end{eqnarray}

\noindent
Evaluating this sum at $q=1$ gives the Euler characteristic 
$(r+1)r(5r+3)$.$\Box$

\section{Betti numbers of flag varieties of
pointed lines}


Let $\a_i$ denote the $i$'th Betti number of the flag variety $\FF(0,1;r)$
of point-line pairs in $\PP^r$ such that the point lies on the line.
Recall from the proof of Lemma \ref{grasser} that
$\ser(\FF(0,1;r))=[r+1][r]$. 
The product $[r+1][r]$ also
appears as a factor in the Serre polynomial (\ref{ser2r2}) 
of $\mm_{0,2}(\PP^r,2)$, making
its coefficients especially relevant to our computations.
It is easy to see that
the Betti numbers  of $\FF(0,1;r)$ initially follow the pattern

$$(1, 2, 3, ...)\text{,}$$
so that for the first half of the Betti numbers we have $\a_i=i+1$.
Since $\dim\FF(0,1;r)=2r-1$ is always odd, it always has an even number
of Betti numbers.  By Poincar\'{e} duality, it
follows that the middle two Betti numbers are both $r$, and the Betti numbers
then decrease back to 1.
It can be checked that all the Betti numbers are given by the following 
formula.
\[\a_i=r+\frac{1}{2}-|r-\frac{1}{2}-i|\]
for $i\in{0}\cup\und{2r-1}$, and $\a_i=0$ otherwise.

\section{Formulas for the Betti numbers of $\mm_{0,2}(\PP^r,2)$}

Let $\b_j$ be the $j$'th Betti number of $\mm_{0,2}(\PP^r,2)$.  We
can write $\b_j=\b_j^1+\b_j^2+\b_j^3$, where $\b_j^i$ is the contribution 
from the $i$'th sum in the last factor of Equation \ref{ser2r2}.
We compute each contribution separately.

We can factor $2q$ out of the second sum and $2q^2$ out of the third sum,
and then temporarily ignore these factors.
It will be easy to recover their effect later by shifting and doubling
some contributions as described below. 

We have reduced our computation to finding coefficients
of expressions of the form $[r+1][r][m]$, where $m\in\{r-1,r+1,r+3\}$.
The degree $j$ term of this polynomial is given by
\[\sum_{i=\max\{0,j-m+1\}}^{\min\{j,2r-1\}}(\a_iq^i)q^{j-i}\text{.}\]
So each coefficient is the sum of some of the Betti numbers $\a_i$ of 
$\FF(0,1;r)$. The sum can have at most $m$ nonzero terms, since this is 
the number of terms in $[m]$.  One could follow the convention that all
sums have $m$ terms, allowing some of the terms to be $\a_i$ that are
zero for dimension reasons.  However, we are restricting
our expressions to include only nonzero $\a_i$, so that the indices
in the sums will be in the range $\{0\}\cup\und{2r-1}$. 
We also adopt the usual convention that
empty sums---those whose upper index is smaller than their lower index---are
zero. Let $\g_i^j$ be the coefficients resulting from these computations.
Reinserting the factors of $2q$ and $2q^2$ 
factored out of the second and third 
pieces, we find $\b_j^1=\g_i^1$, $\b_j^2=2\g_{j-1}^2$, and 
$\b_j^3=2\g_{j-2}^3$.
Hence we have 

\renewcommand{\baselinestretch}{1}
\small\normalsize
\[
\b_j^1=
\begin{cases}
\sum_{i=0}^{j}\a_i   & \text{ if $j\leq r+2$} \\
\sum_{i=j-r-2}^{j}\a_i   & \text{ if $r+2\leq j\leq 2r-1$} \\
\sum_{i=j-r-2}^{2r-1}\a_i   & \text{ if $2r-1\leq j\leq 3r+1$} 
\end{cases}\text{,}
\]
\[
\b_j^2=
\begin{cases}
2\sum_{i=0}^{j-1}\a_i   & \text{ if $j\leq r+1$} \\
2\sum_{i=j-r-1}^{j-1}\a_i   & \text{ if $r+1\leq j\leq 2r$} \\
2\sum_{i=j-r-1}^{2r-1}\a_i   & \text{ if $2r\leq j\leq 3r$} 
\end{cases}\text{,}
\]
\[
\b_j^3=
\begin{cases}
2\sum_{i=0}^{j-2}\a_i   & \text{ if $j\leq r$} \\
2\sum_{i=j-r}^{j-2}\a_i   & \text{ if $r\leq j\leq 2r+1$} \\
2\sum_{i=j-r}^{2r-1}\a_i   & \text{ if $2r+1\leq j\leq 3r-1$.}
\end{cases}
\]


\vspace{0.05in}

Thus we get the following formulas for the Betti numbers of 
$\mm_{0,2}(\PP^r,2)$:

\renewcommand{\baselinestretch}{1}
\small\normalsize
\[
\b_{j}=
\begin{cases}
 \sum_{i=0}^{j}\a_i+2\sum_{i=0}^{j-1}\a_i+2\sum_{i=0}^{j-2}\a_i   
&\text{ if } j\leq r \\
\\
\sum_{i=0}^{r+1}\a_i+2\sum_{i=0}^{r}\a_i+2\sum_{i=1}^{r-1}\a_i  
&\text{ if } j=r+1\\
\\
 \sum_{i=j-r-2}^{j}\a_i+2\sum_{i=j-r-1}^{j-1}\a_i+2\sum_{i=j-r}^{j-2}
\a_i  
&\text{ if } r+2\leq j\leq 2r-1 \\
\\
 \sum_{i=r-2}^{2r-1}\a_i+2\sum_{i=r-1}^{2r-1}\a_i+2\sum_{i=r}^{2r-2}
\a_i   &\text{ if } j=2r\\
\\
\sum_{i=j-r-2}^{2r-1}\a_i+2\sum_{i=j-r-1}^{2r-1}\a_i
+2\sum_{i=j-r}^{2r-1}\a_i   &\text{ if } 2r+1\leq j\leq 3r+1\text{.}
\end{cases}
\]


\vspace{0.05in}

We can come up with an especially explicit description of $\b_j$ for 
$j<r$ since we know $\a_i=i+1$ for $i<r$. Also, $\a_r=r$, which gives the
second part below.

\begin{cor}
\begin{enumerate}
\item For $j<r$, the $j$'th Betti number of $\mm_{0,2}(\PP^r,2)$ is 
\[\b_j=\frac{5}{2}j^2+\frac{3}{2}j+1\text{.}\]

\item Furthermore
\[\b_r=\frac{5}{2}r^2+\frac{3}{2}r\text{.}\]

\end{enumerate}
\end{cor}

\noindent
{\bf Proof.} Simplify
\[\b_j=\frac{(j+1)(j+2)}{2}+j(j+1)+j(j-1)\]
and
\[\b_r=\frac{(r+1)(r+2)}{2}+r(r+1)+r(r-1)-1\text{.}\]
\begin{flushright}$\Box$\end{flushright}

\noindent
As a consequence of this, a particular Betti number of $\mm_{0,2}(\PP^r,2)$
stabilizes as $r$ becomes large.

\begin{cor}
For all $r>j$, the $j$'th Betti number of $\mm_{0,2}(\PP^r,2)$ is
$\b_j=\frac{(j+1)(j+2)}{2}+j(j+1)+j(j-1)$.
\end{cor}

\noindent
Let $\bar{\b}_j$ be this limiting value. We have
\[\bar{\b}_0=1,
\bar{\b}_1=5,
\bar{\b}_2=14,
\bar{\b}_3=28,
\bar{\b}_4=47,
\bar{\b}_5=71,\ldots\]

\renewcommand{\baselinestretch}{1}
\section{Poincar\'{e} polynomials of $\mm_{0,1}(\PP^r,2)$ and 
$\mm_{0,2}(\PP^r,2)$
for small $r$}

Using the same procedure as above,
one can easily compute the Poincar\'{e} polynomial of 
$\mm_{0,1}(\PP^r,2)$.

\begin{prop}
If $r$ is even, the Poincar\'{e} polynomial of $\mm_{0,1}(\PP^r,2)$ is
\[\ser(\mm_{0,1}(\PP^r,2))= 
\left(\sum_{i=0}^r q^i\right)\left(\sum_{i=0}^{(r-2)/2} q^{2i}\right)
\left(\sum_{i=0}^{r+2} q^i+\sum_{i=1}^{r+1} q^i+\sum_{i=2}^{r} q^i\right)
\text{,}\]
and if $r$ is odd, the Poincar\'{e} polynomial of $\mm_{0,1}(\PP^r,2)$ is
\[\ser(\mm_{0,1}(\PP^r,2))= 
\left(\sum_{i=0}^{r-1} q^i\right)\left(\sum_{i=0}^{(r-1)/2} q^{2i}\right)
\left(\sum_{i=0}^{r+2} q^i+\sum_{i=1}^{r+1} q^i+\sum_{i=2}^{r} q^i\right)
\text{.}\]
\end{prop}

\noindent
Thus, for small values of $r$, we get the explicit Poincar\'{e} 
polynomials listed in Table \ref{ep012r2}.

\renewcommand{\baselinestretch}{1}
\small\normalsize

\begin{sidewaystable}
\begin{tabular}{|rrl|}
\hline
$r$ & $\chi(X)$ & $\ser(X)$ \\
1 &  6 &  $1+2q+2q^2+q^3$ \\
2 & 27 & $1+3q+6q^2+7q^3+6q^4+3q^5+q^6$ \\
3 & 72 & $1+3q+7q^2+11q^3+14q^4+14q^5+11q^6+7q^7+3q^8+q^9$ \\
4 & 150 & $1+3q+7q^2+12q^3+18q^4+22q^5+24q^6+22q^7+18q^8+12q^9+7q^{10}+3q^{11}+q^{12}$ \\
5 & 270 & $1+3q+7q^2+12q^3+19q^4+26q^5+32q^6+35q^7$ \\
& & \hspace{0.1in} $+35q^8+32q^9+26q^{10}+19q^{11}+12q^{12}+7q^{13}+3q^{14}+q^{15}$ \\
6 & 441 & $1+3q+7q^2+12q^3+19q^4+27q^5+36q^6+43q^7+48q^8+49q^9$ \\
& & \hspace{0.1in} $+48q^{10}+43q^{11}+36q^{12}+27q^{13}+19q^{14}+12q^{15}+7q^{16}+3q^{17}+q^{18}$ \\
7 & 672 & $1+3q+7q^2+12q^3+19q^4+27q^5+37q^6+47q^7+56q^8+62q^9+65q^{10}$ \\
& & \hspace{0.1in} $+65q^{11}+62q^{12}+56q^{13}+47q^{14}+37q^{15}+27q^{16}+19q^{17}+12q^{18}+7q^{19}+3q^{20}+q^{21}$ \\
8 & 972 & $1+3q+7q^2+12q^3+19q^4+27q^5+37q^6+48q^7+60q^8+70q^9+78q^{10}+82q^{11}+84q^{12}$ \\
& & \hspace{0.1in} $+82q^{13}+78q^{14}+70q^{15}+60q^{16}+48q^{17}+37q^{18}+27q^{19}+19q^{20}+12q^{21}+7q^{22}+3q^{23}+q^{24}$ \\
\hline
\end{tabular}

\begin{tabular}{|rrl|}
\hline
$r$ & $\chi(Y)$ & $\ser(Y)$ \\
1 & 16 & $1+4q+6q^2+4q^3+q^4$  \\
2 & 78 & $1+5q+13q^2+20q^3+20q^4+13q^5+5q^6+q^7$   \\
3 & 216 & $1+5q+14q^2+27q^3+39q^4+44q^5+39q^6+27q^7+14q^8+5q^9+q^{10}$  \\
4 & 460 & $1+5q+14q^2+28q^3+46q^4+63q^5+73q^6+73q^7+63q^8+46q^9+28q^{10}+14q^{11}+5q^{12}+q^{13}$   \\
5 & 840 & $1+5q+14q^2+28q^3+47q^4+70q^5+92q^6+107q^7+112q^8$ \\
& & \hspace{0.1in} $+107q^9+92q^{10}+70q^{11}+47q^{12}+28q^{13}+14q^{14}+5q^{15}+q^{16}$   \\
6 & 1386 & $1+5q+14q^2+28q^3+47q^4+71q^5+99q^6+126q^7+146q^8+156q^9$ \\
& & \hspace{0.1in} $+156q^{10}+146q^{11}+126q^{12}+99q^{13}+71q^{14}+47q^{15}+28q^{16}+14q^{17}+5q^{18}+q^{19}$   \\
7 & 2128 & $1+5q+14q^2+28q^3+47q^4+71q^5+100q^6+133q^7+165q^8+190q^9+205q^{10}+210q^{11}$ \\
& & \hspace{0.1in} $+205q^{12}+190q^{13}+165q^{14}+133q^{15}+100q^{16}+71q^{17}+47q^{18}+28q^{19}+14q^{20}+5q^{21}+q^{22}$   \\
8 & 3096 & $1+5q+14q^2+28q^3+47q^4+71q^5+100q^6+134q^7+172q^8+209q^9+239q^{10}+259q^{11}+269q^{12}$ \\
& & \hspace{0.1in} $+269q^{13}+259q^{14}+239q^{15}+209q^{16}+172q^{17}+134q^{18}+100q^{19}+71q^{20}$ \\
& & \hspace{0.1in} $+47q^{21}+28q^{22}+14q^{23}+5q^{24}+q^{25}$ \\
\hline
\end{tabular}
\caption{Euler characteristics and Poincar\'{e} polynomials for 
$X=\mm_{0,1}(\PP^r,2)$ and $Y=\mm_{0,2}(\PP^r,2)$ \label{ep012r2}}
\end{sidewaystable}



\chapter{Generators}
\label{sec:gen}

This section describes three types of divisor classes found in all
Chow rings \linebreak[4] $A^*(\m0)$. We will
see in Section \ref{sec:comp} that divisors 
of these types generate the Chow
ring of the moduli space $\mm_{0,2}(\PP^r,2)$. All of these divisors
have equivariant lifts, produced by taking the equivariant first Chern classes
of their corresponding equivariant line bundles. We will use the same
symbols for these equivariant versions; the meaning should be clear from
the context.

\section{Hyperplane pullbacks}

Let $\ev_1,\ldots,\ev_n$ be the evaluation maps on $\mmm$ defined in Section
\ref{sec:modintro}. 
Any cohomology class on $X$ pulls back under any evaluation
map to a class on $\mmm$, which
may be considered as an element in the Chow ring under the isomorphism 
(\ref{hi}).  In the case $X=\PP^r$, we 
can pull the hyperplane class $H$ back under each evaluation, 
getting the $n$ divisor classes $H_i=\ev_i^*(H)$.  

\section{Boundary divisors}
\label{sec:bound}

The boundary of $\mm_{0,n}(X,\b)$ by definition consists
of the locus of stable maps with reducible domain curves. 
By Proposition \ref{nicemod}, it is a divisor with normal crossings.
The irreducible components of 
the boundary of $\m0$
are in 1--1 correspondence with quadruples $(A,d_A,B,d_B)$,
where $A,B\sub\und{n}$ partition $\und{n}$, $d_1+d_2=d$, and if 
$d_A=0$ (resp. $d_B=0$), then $A$ (resp. $B$) has at least two elements. 
Such a boundary divisor and its
class in the Chow ring are both denoted $D_{A,d_A,B,d_B}$. Geometrically,
the divisor $D_{A,d_A,B,d_B}$ corresponds to the closure of the
locus of stable maps where the domain curve has two 
components, one having marked points labeled by $A$ and mapping to $\PP^r$
with degree $d_A$, and the other having marked points labeled by $B$ and 
mapping to $\PP^r$ with degree $d_B$. We represent this divisor by the
following picture.

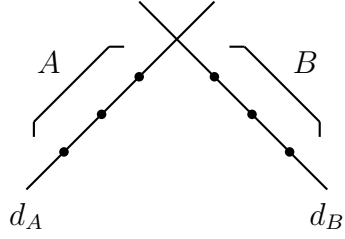
\begin{center}
\begin{pspicture}(0,0)(5,3.5)
\pnode(0.5,.5){a}
\dotnode(1,1){b}
\dotnode(1.5,1.5){c}
\dotnode(2,2){d}
\pnode(3,3){e}
\pnode(2,3){f}
\dotnode(3,2){g}
\dotnode(3.5,1.5){h}
\dotnode(4,1){i}
\pnode(4.5,0.5){j}
\pnode(.6,1.2){k}
\pnode(.6,1.4){l}
\pnode(1.6,2.4){m}
\pnode(1.8,2.4){n}
\pnode(4.4,1.2){o}
\pnode(4.4,1.4){p}
\pnode(3.4,2.4){q}
\pnode(3.2,2.4){r}
\ncline{a}{e}
\ncline{f}{j}
\ncline{k}{l}
\ncline{l}{m}\naput{$A$}
\ncline{m}{n}
\ncline{o}{p}
\ncline{q}{p}\naput{$B$}
\ncline{q}{r}
\uput{5pt}[d](.5,.5){$d_A$}
\uput{5pt}[d](4.5,.5){$d_B$}
\end{pspicture}
\end{center}
Note that the domain curves of stable maps lying in a boundary
divisor may have more than two components. In the limit, some combinations
of marked points and the node may coincide, causing new components to sprout.
The marked points involved in the ``collision" will appear on this new
component. Observe that, although marked points can thus
migrate to newly added
components in the limit, they cannot move onto components that already existed
as long as the node separating the components is maintained.
Additionally, the map itself may degenerate in such a way that the number of
components of the domain curve increases. 
The diagram representation given
above for divisors can easily be extended to describe the closures of
other degeneration loci. This description is an alternative to the
dual graphs used in Section \ref{sec:poincare}. We will use dual graphs when
referring to the degeneration loci themselves, and we will use this diagram
representation when referring to their closures. Such a diagram therefore
directly describes only a generic element of the locus it represents. For
example, the diagram

\begin{center}
\begin{pspicture}(0,0)(5,2)
\pnode(1,1){a}
\dotnode(2,1){c}
\dotnode(3,1){d}
\pnode(4,1){b}
\ncline{a}{b}
\uput{5pt}[r](4,1){2}
\end{pspicture}
\end{center}
represents (a generic element of) $\mm_{0,2}(\PP^r,2)$ itself.

In $A^*(\mm_{0,2}(\PP^r,2))$, there are exactly three boundary divisors.
We use the notation  
$D_0=D_{\underline{2},0,\emptyset,2}$, $D_1=D_{\underline{2},1,\emptyset,1}$, 
and
$D_2=D_{\underline{1},1,\{2\},1}$ for these divisors.  
Thus the domain of a generic stable map in $D_0$ has one collapsed component 
containing both marked points and one component of degree two. 
Generic elements of $D_1$ are maps which have degree one on each of the two
components of the domain curve, with both 
marked points lying on the same component.
Finally, $D_2$ is the boundary
divisor whose generic maps have degree one on each component with one marked 
point on each component. Note that $D_i$ corresponds to curves with $i$
degree one components containing marked points; this aids in remembering the 
notation.
The diagrams for these three divisors 
follow.

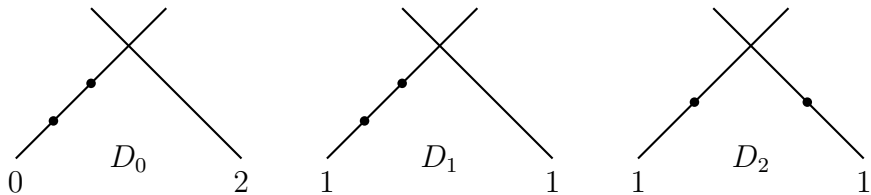
\begin{center}
\begin{pspicture}(0,0)(4,3)
\rput(2,0.5){$D_0$}
\pnode(0.5,.5){a}
\dotnode(1.5,1.5){c}
\dotnode(1,1){d}
\pnode(2.5,2.5){e}
\pnode(1.5,2.5){f}
\pnode(3.5,0.5){j}
\ncline{a}{e}
\ncline{f}{j}
\uput{5pt}[d](.5,.5){0}
\uput{5pt}[d](3.5,.5){2}
\end{pspicture}
\begin{pspicture}(0,0)(4,3)
\rput(2,0.5){$D_1$}
\pnode(0.5,.5){a}
\dotnode(1.5,1.5){c}
\dotnode(1,1){d}
\pnode(2.5,2.5){e}
\pnode(1.5,2.5){f}
\pnode(3.5,0.5){j}
\ncline{a}{e}
\ncline{f}{j}
\uput{5pt}[d](.5,.5){1}
\uput{5pt}[d](3.5,.5){1}
\end{pspicture}
\begin{pspicture}(0,0)(4,3)
\rput(2,0.5){$D_2$}
\pnode(0.5,.5){a}
\dotnode(1.25,1.25){c}
\dotnode(2.75,1.25){d}
\pnode(2.5,2.5){e}
\pnode(1.5,2.5){f}
\pnode(3.5,0.5){j}
\ncline{a}{e}
\ncline{f}{j}
\uput{5pt}[d](.5,.5){1}
\uput{5pt}[d](3.5,.5){1}
\end{pspicture}
\end{center}

We use $D_1$ as an example to 
illustrate the further degeneration that can occur within a boundary
divisor. Contained within $D_1$ are loci with the following diagrams.

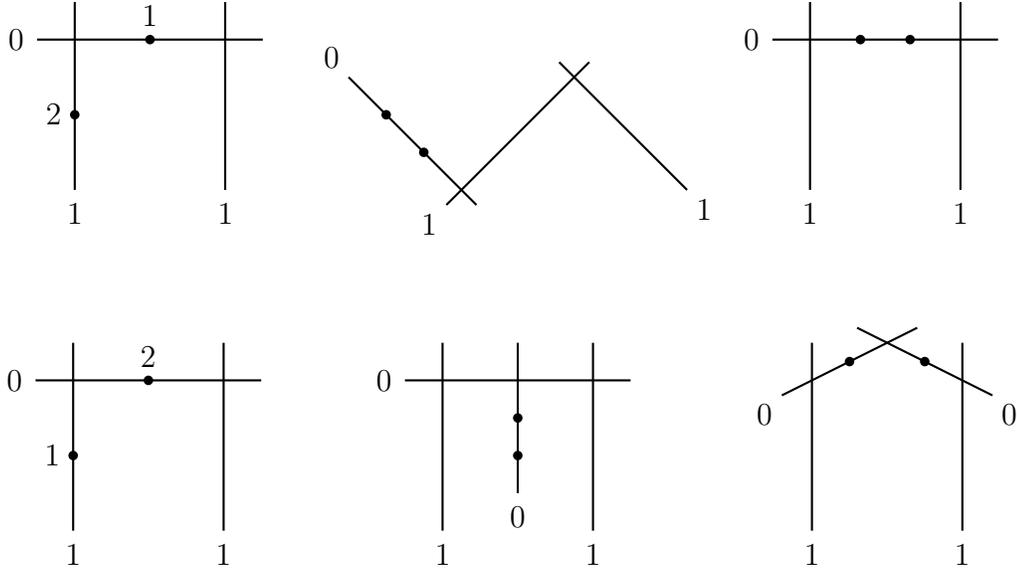
\begin{center}
\begin{pspicture}(0,0)(4,3.5)
\pnode(1,.5){a}
\pnode(1,3){b}
\pnode(.5,2.5){c}
\pnode(3.5,2.5){d}
\pnode(3,3){e}
\pnode(3,.5){f}
\dotnode(1,1.5){i}
\dotnode(2,2.5){j}
\ncline{a}{b}
\ncline{c}{d}
\ncline{e}{f}
\uput{5pt}[d](1,.5){1}
\uput{5pt}[l](.5,2.5){0}
\uput{5pt}[d](3,.5){1}
\uput{5pt}[l](1,1.5){2}
\uput{5pt}[u](2,2.5){1}
\end{pspicture}
\begin{pspicture}(0,0)(5.5,2.5)
\pnode(.5,2){a}
\pnode(2.2,.3){b}
\pnode(1.8,.3){c}
\pnode(3.7,2.2){d}
\pnode(3.3,2.2){e}
\pnode(5,.5){f}
\dotnode(1,1.5){i}
\dotnode(1.5,1){j}
\ncline{a}{b}
\ncline{c}{d}
\ncline{e}{f}
\uput{5pt}[ul](.5,2){0}
\uput{5pt}[dl](1.8,.3){1}
\uput{5pt}[dr](5,.5){1}
\end{pspicture}
\begin{pspicture}(0,0)(4,3.5)
\pnode(1,.5){a}
\pnode(1,3){b}
\pnode(.5,2.5){c}
\pnode(3.5,2.5){d}
\pnode(3,3){e}
\pnode(3,.5){f}
\dotnode(1.67,2.5){i}
\dotnode(2.33,2.5){j}
\ncline{a}{b}
\ncline{c}{d}
\ncline{e}{f}
\uput{5pt}[d](1,.5){1}
\uput{5pt}[l](.5,2.5){0}
\uput{5pt}[d](3,.5){1}
\end{pspicture}

\begin{pspicture}(0,0)(4,3.5)
\pnode(1,.5){a}
\pnode(1,3){b}
\pnode(.5,2.5){c}
\pnode(3.5,2.5){d}
\pnode(3,3){e}
\pnode(3,.5){f}
\dotnode(1,1.5){i}
\dotnode(2,2.5){j}
\ncline{a}{b}
\ncline{c}{d}
\ncline{e}{f}
\uput{5pt}[d](1,.5){1}
\uput{5pt}[l](.5,2.5){0}
\uput{5pt}[d](3,.5){1}
\uput{5pt}[l](1,1.5){1}
\uput{5pt}[u](2,2.5){2}
\end{pspicture}
\hspace{.25in}
\begin{pspicture}(0,0)(4,3.5)
\pnode(1,.5){a}
\pnode(1,3){b}
\pnode(.5,2.5){c}
\pnode(3.5,2.5){d}
\pnode(3,3){e}
\pnode(3,.5){f}
\pnode(2,1){g}
\pnode(2,3){h}
\dotnode(2,1.5){i}
\dotnode(2,2){j}
\ncline{a}{b}
\ncline{c}{d}
\ncline{e}{f}
\ncline{g}{h}
\uput{5pt}[d](1,.5){1}
\uput{5pt}[l](.5,2.5){0}
\uput{5pt}[d](3,.5){1}
\uput{5pt}[d](2,1){0}
\end{pspicture}
\hspace{.25in}
\begin{pspicture}(0,0)(4,4.5)
\pnode(1,.5){a}
\pnode(1,3){b}
\pnode(.6,2.3){c}
\pnode(2.4,3.2){d}
\pnode(1.6,3.2){e}
\pnode(3.4,2.3){f}
\pnode(3,.5){g}
\pnode(3,3){h}
\dotnode(1.5,2.75){i}
\dotnode(2.5,2.75){j}
\ncline{a}{b}
\ncline{c}{d}
\ncline{e}{f}
\ncline{g}{h}
\uput{5pt}[d](1,.5){1}
\uput{5pt}[dl](.6,2.3){0}
\uput{5pt}[dr](3.4,2.3){0}
\uput{5pt}[d](3,.5){1}
\end{pspicture}
\end{center}

\noindent
The marked points are not labeled in diagrams where the distinction does
not affect the boundary class, either because both marked points lie on
the same component or because of symmetry.

We will show in Section \ref{sec:lla} that the three boundary divisor classes 
together with the hyperplane pullbacks $H_1$ and $H_2$ generate the linear
part of the ring $A^*(\mm_{0,2}(\PP^1,2))$. Even better, we will see 
in Section \ref{sec:comp} that these classes
generate this entire Chow ring.

\section{The $\psi$-classes}
\label{sec:psi}

We use definitions and properties from \cite{HM}. Every flat family 
$\phi:\mathcal{C}\rightarrow B$ of nodal curves over a scheme has a 
relative dualizing sheaf $\w_{\mathcal{C}/B}$, or $\w_{\phi}$, 
defined to be the sheaf of rational relative 
differentials. This sheaf is invertible, so we may also consider it as a line
bundle.  If the total space is smooth, we can write
\[\w_\phi=K_\mathcal{C}\*\phi^*K_B^\vee\text{,}\]
where $K_\mathcal{C}$ and $K_B$ are the canonical bundles.

If the family is equipped with sections $s_1,...,s_n$, we may 
consider the bundles $L_i=s_i^*(\w_{\phi})$, called the cotangent line bundles.
At a point $b\in B$, the fiber of $L_i$ is the cotangent space to the curve
$\mathcal{C}_b$ at the point $s_i(b)$. We also define the {\em $\psi$-class}
$\psi_i$ to be the first Chern class $c_1(L_i)$ for each $i$.

We extend this setup to moduli stacks of stable maps and the universal 
curves of their universal families.  We saw in Section \ref{sec:mod} that
the universal curve of $\m0$ is 
$\mm_{0,n+1}(\PP^r,d)$, which is a smooth stack. 
Everything said above carries
over to this
stack setting.  The universal curve is equipped with $n$ universal
sections $\s_1,\ldots,\s_n$, so we have $n$ naturally defined $\psi$-classes. 
Furthermore, it is straightforward to check that these $\psi$-classes are
universal as well: given any morphism $g:S\ra\m0$, the pullbacks 
$g^*(\psi_1),\ldots,g^*(\psi_n)$ are the $\psi$-classes on the induced family.

Although the $\psi$-classes are not strictly necessary as generators for
the rings $A^*(\mm_{0,2}(\PP^1,2))$, we include them because their geometric
nature makes some of the relations much easier to understand and state.
Compare the presentations of
Theorem \ref{thm:prez} and Proposition \ref{altprez} for immediate and
visible testimony to the value of such incorporation.
We will give here one example of the usefulness of the $\psi$-classes
in describing geometric conditions.

To make this example easier to state, we first introduce a slight
modification of the concept of $\psi$-classes.
Restricting to the closure of a particular degeneration locus, let $p$
be a node at the intersection of components $E_1$ and $E_2$, like in the
diagrams of Section \ref{sec:bound}. Then we can define classes
$\tilde{\psi}_{p,E_i}$ essentially just like we defined $\psi$-classes.  The
only difference here is that we
additionally specify which branch to consider $p$ as lying on, 
so that the cotangent
space is one-dimensional. In fact, this type of class is defined at the end of
Section \ref{sec:eq}, where it is denoted $e_F$. 
If, reversing the gluing process described in
Section \ref{sec:modintro}, we remove, say, $E_2$ and then the 
associated connected component of the curve, replacing the node with an
auxiliary marked point $s_{\bullet}$, then $\tilde{\psi}_{p,E_1}$ becomes
a legitimate $\psi$-class $\psi_{\bullet}$ on the resulting moduli space.

An important fact which will be used many times throughout the dissertation
is that a collapsed rational component with exactly three special points is a
{\em rigid object}. In other words, the marked points and
nodes on such a component
cannot be moved around internally on the component. 
Among other things, this says that, once such a component appears in
a degeneration, it will remain in any further degeneration. We have already
used this fact implicitly, for example, in describing the possible further 
degenerations of $D_1$ in Section \ref{sec:bound}.

This rigidity is intuitively clear from at least two different perspectives.
First, there exists an automorphism of $\PP^1$ taking any three distinct
points to any other three distinct points. So, up to isomorphism, any two data
consisting of three marked points in $\PP^1$ are equivalent. We might as well
always take the points to be 0, 1, and $\infty$ (or any three arbitrarily
chosen points). Second, if such points and nodes were allowed to
move on the domain of a stable map, they could come together and force the
sprouting of new components. However, the result would no longer be a
stable map; there aren't enough special points left for the collapsed 
components.

More rigorously, the rigidity of such a component is equivalent to the 
vanishing of the corresponding $\psi$-classes and $\tilde{\psi}$-classes
on that component.
For $\psi$-classes, this is because the cotangent 
line bundles are trivial if and only if the marked points are fixed.
A similar statement holds for $\tilde{\psi}$-classes. We will prove in
Section \ref{sec:geomrel} that  $\psi_1$ and $\psi_2$
vanish on $D_0$. This is the most important special case of the rigidity
just described. From the argument given there, it is easy to see that 
this characterization of the rigidity of a component holds in general.

\renewcommand{\baselinestretch}{1}
\chapter{Presentations for the Chow rings of some simpler spaces}
\label{sec:simpler}

\renewcommand{\baselinestretch}{1}
\section{Presentations for $A^*(\mm_{0,1}(\PP^1,1))$, 
$A^*(\mm_{0,2}(\PP^1,1))$, and $A^*(\mm_{0,3}(\PP^1,1))$}
\label{sec:0n11}

The study of degree one
stable maps to $\PP^1$ is relatively simple because every degree one
morphism from $\PP^1$ to $\PP^1$ is an isomorphism. More generally, genus 
zero, degree one stable maps never have nontrivial automorphisms. Thus the
corresponding moduli spaces may just as well be considered as fine moduli
schemes, as no loss of information results. 
Let us record the standard fact that
\begin{equation}\label{projprod}
A^*(\prod_{i=1}^n\PP^{r_i})=
\frac{\QQ[H_1,\ldots,H_n]}{\left(H_1^{r_1+1},\ldots,H_n^{r_n+1}\right)}
\text{,}
\end{equation}
where $H_i$ is the pullback of the hyperplane class under the $i$'th 
projection. See \cite[Chapter 8]{F}.

First, we mention that
the moduli space $\mm_{0,0}(\PP^1,1)$ is a point because the domain
of such a stable map is always $\PP^1$, and, thanks to the opening comment
above, all degree one stable 
maps from $\PP^1$ to itself are isomorphic to the 
identity.

Second, this seems like a good place to note that 
$\mm_{0,3}(\PP^1,0)\iso\MM_{0,3}\x\PP^1\iso\PP^1$ since
$\MM_{0,3}$ is also a point.
Since the automorphism group of $\PP^1$ is three-dimensional, there is
an automorphism mapping
any three distinct points to any other three distinct points. Thus we
can fix three points, say 0, 1, and $\infty$ of $\PP^1$, and any data
$(\PP^1,x_1,x_2,x_3)$ of three distinct points in $\PP^1$ is isomorphic
to our fixed data. In effect, the marked points of $\MM_{0,3}$, and hence
by pullback those of $\mm_{0,3}(\PP^1,0)$, are
not allowed to vary.
(See the discussion at the end of Section \ref{sec:psi}.
Incidentally, $\MM_{0,4}\iso\PP^1$ roughly because, after
fixing the first three marked points, the fourth is allowed to vary over
$\PP^1$. See \cite{K} for more detail.)

Thus
\[A^*(\mm_{0,0}(\PP^1,1))\iso\QQ \text{\hspace{.5in}and\hspace{.5in}}
A^*(\mm_{0,3}(\PP^1,0))\iso\QQ[H]/(H^2)\text{.}\]
The classes in the former are just multiples of its fundamental class. In
the latter, $H$ corresponds to the hyperplane pullback under any of the three
evaluation morphisms, which all simply 
record the image of the trivial stable map.
We now proceed to the three main subjects of this section.

\begin{lem} 
$\mm_{0,1}(\PP^1,1)\iso \PP^1$.
\end{lem}

\noindent
{\bf Proof.} The family of one-pointed, degree one stable maps

\vspace{0.3in}

\begin{center}
\psset{arrows=->}
\begin{psmatrix}
$\PP^1\times\mathbb{P}^1$ & $\mathbb{P}^1$ \\ 
$\mathbb{P}^1$
\ncline[offset=-3pt]{1,1}{2,1}\naput{$\pj_1$}
\ncarc[arcangle=45]{2,1}{1,1}\naput{$\D$}
\ncline{1,1}{1,2}\naput{$\pj_2$}
\end{psmatrix}
\end{center}

\vspace{0.1in}

\noindent
gives a morphism $\a:\PP^1\ra\mm_{0,1}(\PP^1,1)$. The evaluation morphism
$\ev:\mm_{0,1}(\PP^1,1)\ra\PP^1$ is easily checked to be inverse to $\a$,
and hence is the desired isomorphism.$\Box$

\noindent
Under this isomorphism, the hyperplane $H$ in $\PP^1$ naturally corresponds
to its pullback $H_1=\ev^{-1}H$.
Therefore
\[A^*(\mm_{0,1}(\PP^1,1))\iso\frac{\QQ[H_1]}{(H_1^2)}\text{.}\]

\begin{lem} \label{m0211}
$\mm_{0,2}(\PP^1,1)\iso \PP^1\x\PP^1$
\end{lem}

\noindent
{\bf Proof.} We start as before with a family

\vspace{0.2in}

\begin{center}
\psset{arrows=->}
\begin{psmatrix}
$(\PP^1)^2\times\mathbb{P}^1$ & $\mathbb{P}^1$ \\ 
$(\mathbb{P}^1)^2$
\ncline[offset=-3pt]{1,1}{2,1}\naput{$\pj_1$}
\ncarc[arcangle=45]{2,1}{1,1}\naput{$s_i=(\id,\pj_i)$}
\ncline{1,1}{1,2}\naput{$\pj_2$}
\end{psmatrix}\text{.}
\psset{arrows=-}
\end{center}

\vspace{0.2in}

\noindent
This is not a family of stable maps because the images of the $s_i$
coincide on the diagonal.  To fix this, let $f:\Bl_{\D}(\PP^1)^3\ra(\PP^1)^3$
be the
blowup of $(\PP^1)^3$ along its small diagonal.
The family $((f_1,f_2):\Bl_{\D}(\PP^1)^3\ra(\PP^1)^2,\tilde{s}_1,\tilde{s}_2,
f_3)$
induced by $f$ from the original family {\em is} a family of 
stable maps. (Here the sections $\tilde{s}_i$ are proper transforms.)
To see this, fix $a\in\PP^1$ and consider the hypersurface $Y=Z(z_3-a)$ in 
$(\PP^1)^3$. The blowup $f$ restricts to the blowup of $Y$
at $(a,a,a)$. The restrictions of the sections $s_i$ to $Y$ have distinct
tangent directions at $(a,a,a)$, so their proper transforms will be disjoint
on the exceptional $\PP^1$ over $(a,a,a)$.

Let $\a:(\PP^1)^2\ra\mm_{0,2}(\PP^1,1)$ be the morphism induced by this 
family.
We will check that $\ev=(\ev_1,\ev_2)$ is inverse to $\a$, and hence
is the desired isomorphism. An alternate way to compactify
$\nobarmm_{0,2}(\PP^1,1)$ is by allowing the marked points to coincide.
The resulting moduli space is equivalent to $\mm_{0,2}(\PP^1,1)$. Once again,
this is due to the rigidity of a collapsed component with three special
points, as discussed in Section \ref{sec:psi}, since such a component
results as a limit of stable maps when the marked points would otherwise 
coincide. Using this equivalent description of $\mm_{0,2}(\PP^1,1)$,
we may assume that the map associated to any
stable map in $\mm_{0,2}(\PP^1,1)$ is the identity. Thus $\a$ takes
an ordered pair $(z_1,z_2)$ to the stable map $(\PP^1,z_1,z_2,\id)$. Simple
observation now confirms that the compositions $(\ev_1,\ev_2)\circ\a$ and
$\a\circ(\ev_1,\ev_2)$ are both identities.
$\Box$

\noindent
Under this isomorphism, the hyperplanes $H_i$ in $(\PP^1)^2$ naturally
correspond to the hyperplane pullbacks $H_i$ in $\mm_{0,2}(\PP^1,1)$.
Hence
\[A^*(\mm_{0,2}(\PP^1,1))\iso\frac{\QQ[H_1,H_2]}{(H_1^2,H_2^2)}\text{.}\]

\begin{prop}
\label{m0311}
We have an isomorphism of $\QQ$-algebras 

\begin{eqnarray*}
& & A^*(\mm_{0,3}(\PP^1,1)) \\
& \iso & \frac{\QQ[H_1,H_2,H_3,D]}{(H_1^2,H_2^2,H_3^2, 
(H_1+H_2-D)(H_2+H_3-D), D(H_1-H_2), D(H_2-H_3))}
\end{eqnarray*}

\end{prop}

\noindent
{\bf Proof.}  Since $\mm_{0,3}(\PP^1,1)$ is the universal curve over
$\mm_{0,2}(\PP^1,1)$, the proof of Lemma \ref{m0211} shows that
$\mm_{0,3}(\PP^1,1) \iso\Bl_\D{(\PP^1)^3}$.
Now $\D\iso\PP^1$, and the restriction map 
corresponding to $i:\D\hookrightarrow(\PP^1)^3$ sends each $H_i$
to the hyperplane class in $\D$. So $i^*:A^* ((\PP^1)^3)\rightarrow
A^*(\D)$ is surjective, and we can take $H_1-H_2$ and $H_2-H_3$ as
generators for $\ker{i^*}$. The small diagonal is the complete
intersection of any two of the large diagonals. So we may apply Keel's
Lemma 1 from
\cite{K}, which says that whenever $X$ is a complete intersection of two
divisors $D_1$,$D_2$ in a scheme $Y$ and the restriction map 
$i^*:A^*(Y)\ra A^*(X)$ is
surjective, then
\[A^*(\tilde{Y})=A^*(Y)[T]/((D_1-T)(D_2-T), \ker{i^*}\cdot T)\text{,}\]
where $\tilde{Y}$ is the blowup of $Y$ along $X$. 
Here $T$ corresponds to the exceptional divisor. We know from \cite{F} that 
$A^*((\PP^1)^3)=\QQ[H_1,H_2,H_3]/(H_1^2,H_2^2,H_3^2)$ and that 
we can express two of the 
large 
diagonal classes as $D_i=H_i+H_{i+1}$ for $i\in\{1,2\}$.  
The expression in the proposition results.$\Box$

\noindent
As before, the $H_i$ in the presentation are naturally identified with 
the corresponding hyperplane pullbacks. Furthermore $D$ corresponds to
the boundary divisor $D=D_{\und{3},0,\emptyset,1}$ whose generic stable map
has all three marked points
lying on the same collapsed component.

\renewcommand{\baselinestretch}{1}
\section{The Chow ring $A^*(\mm_{0,1}(\PP^1,2))$ via equivariant 
cohomology}
\label{sec:0112}

The goal of this subsection
is to give a presentation for $A^*(\mm_{0,1}(\PP^1,2))$.
Our main tool is the method of localization in equivariant cohomology, 
developed in Section \ref{sec:eq}. We will employ additional methods from
\cite[Section 9.2]{CK} to aid in applying localization to this particular
moduli space.
From Table \ref{ep012r2}, 
\[\ser(\mm_{0,1}(\PP^1,2))=q^3+2q^2+2q+1\text{.}\]
By Section \ref{sec:gen}, 
we have the divisor classes $H_1=\ev_1^*(H)$, the pullback of the 
hyperplane class, and $D=D_{\underline{1},1,\emptyset,1}$, the lone boundary
divisor. We will show that these two classes generate the Chow ring, and we
will also find two relations involving them.

Let $\mm=\mm_{0,1}(\PP^1,2)$. Let $T=(\CC^*)^2$. 
Consider the usual $T$-action on $\mm$. We want to compute the equivariant
integrals of degree three monomials in the classes above.
 First,  we have to find the 
equivariant Euler classes of the normal bundles of the fixed point components.
 We do this using Theorem \ref{norm}.
There are six fixed
components, all of them isolated points.  We will label the graphs 
corresponding to the fixed components as follows.

\[\G_1=\text{
\psset{labelsep=2pt, tnpos=b,radius=2pt}
\pstree[treemode=R]{\TC*~{0}~[tnpos=a]{\{1\}}}
{
\TC*~{1}\taput{2} 
}
}\]

\[\G_2=\text{
\psset{labelsep=2pt, tnpos=b,radius=2pt}
\pstree[treemode=R]{\TC*~{1}~[tnpos=a]{\{1\}}}
{
\TC*~{0}\taput{2} 
}
}\]

\[\G_3=\text{
\psset{labelsep=2pt, tnpos=b,radius=2pt}
\pstree[treemode=R]{\TC*~{1}~[tnpos=a]{\{1\}}}
{
\pstree{\TC*~{0}\taput{1}}
{
\TC*~{1}\taput{1} 
}
}
}\]

\[\G_4=\text{
\psset{labelsep=2pt, tnpos=b,radius=2pt}
\pstree[treemode=R]{\TC*~{0}~[tnpos=a]{\{1\}}}
{
\pstree{\TC*~{1}\taput{1}}
{
\TC*~{0}\taput{1} 
}
}
}\]

\[\G_5=\text{
\psset{labelsep=2pt, tnpos=b,radius=2pt}
\pstree[treemode=R]{\TC*~{1}}
{
\pstree{\TC*~{0}~[tnpos=a]{\{1\}}\taput{1}}
{
\TC*~{1}\taput{1} 
}
}
}\]

\[\G_6=\text{
\psset{labelsep=2pt, tnpos=b,radius=2pt}
\pstree[treemode=R]{\TC*~{0}}
{
\pstree{\TC*~{1}~[tnpos=a]{\{1\}}\taput{1}}
{
\TC*~{0}\taput{1} 
}
}
}\]

\vspace{0.1in}

For all $i$, let $Z_i$ denote the fixed component corresponding to $\G_i$.
As we will see, all degree three classes restrict to zero on $Z_1$
and $Z_2$, 
so their equivariant Euler classes are not needed for our computations. Before
beginning computations for the other components, we note the following fact
about the term $e_F$ that appears in the formula for $e_{\G}^{\text{F}}$ in
Theorem \ref{norm}.  
Since it is a first Chern class with zero weight, we must have
$e_F=0$ on any fixed component which is a point. We 
label the vertices and edges of the remaining graphs as follows.

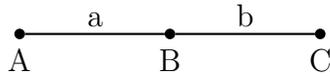
\begin{center}
\psset{labelsep=2pt, tnpos=b,radius=2pt}
\pstree[treemode=R]{\TC*~{A}}
{
\pstree{\TC*~{B}\taput{a}}
{
\TC*~{C}\taput{b} 
}
}
\end{center}

Then for $Z_3$ we have
\[e_{\G_{3}}^{\text{F}}=\frac{1}{(\la_0-\la_1)^2(\la_1-\la_0)^2}
\text{,}\]

\[e_{\G_{3}}^{\text{v}}=\frac{(\la_1-\la_0)^2(\la_0-\la_1)(\la_0-\la_1+\la_0-\la_1)}
{\la_1-\la_0}=2(\la_1-\la_0)^3
\text{,}\]
and
\[e_{\G_{3}}^{\text{e}}=(\la_0-\la_1)^4
\text{.}\]
Thus
\[\eqeul(N_{\G_{3}})=\frac{2(\la_1-\la_0)^3(\la_0-\la_1)^4}
{(\la_1-\la_0)^2(\la_0-\la_1)^2}=2(\la_1-\la_0)^3
\text{.}\]
A similar calculation shows that $\eqeul(N_{\G_{4}})=2(\la_0-\la_1)^3$.
Now for $Z_5$,
\[e_{\G_{5}}^{\text{F}}=\frac{\w_{B,a}\w_{B,b}}{(\la_0-\la_1)^2(\la_1-\la_0)^2}
=\frac{(\la_0-\la_1)^2}{(\la_0-\la_1)^2(\la_1-\la_0)^2}=\frac{1}{(\la_1-\la_0)^2}
\text{,}\]

\[e_{\G_{5}}^{\text{v}}=\frac{(\la_1-\la_0)^2(\la_0-\la_1)}{\w_{A,a}\w_{C,b}}
=\frac{(\la_1-\la_0)^2(\la_0-\la_1)}{(\la_1-\la_0)^2}=\la_0-\la_1
\text{,}\]
and
\[e_{\G_{5}}^{\text{e}}=(\la_0-\la_1)^4
\text{.}\]
Thus
\[\eqeul(N_{\G_{5}})=\frac{(\la_0-\la_1)(\la_0-\la_1)^4}{(\la_1-\la_0)^2}
=(\la_0-\la_1)^3
\text{.}\]
A similar calculation shows that $\eqeul(N_{\G_{6}})=(\la_1-\la_0)^3$.

Next
we need the restrictions of the divisor classes to the fixed components.
``Restriction" applies in a loose sense here, since the weights given
are actually those on the pullback to an atlas in cases where the fixed
components have automorphisms. This is reconciled later by including 
the factors
$a_j$ in the residue formula of Corollary \ref{stackloc}.
We will focus our comments on the particular case of $\mm$ for simplicity,
although the same ideas apply much more generally with slight modifications.
Each $Z_i$ maps to a fixed point in $\PP^1$ under evaluation. 
Thus the restriction of 
$H_1$ to $Z_i$ amounts to the pullback via the evaluation morphism of
the first Chern class of the restriction of $\O(1)$ to that 
point in $\PP_T^1$.  
The restriction of $\O(1)$ to a point is certainly a trivial 
bundle, but the Chern class is still
non-zero since this bundle carries a $T$-action.
In particular, restricting to the fixed point $q_j$ gives the weight
$\la_j$ for the first Chern class. (See Section \ref{sec:eq} for
notation.) Therefore, we need only look at the image
of the marked point to determine which $\la_j$ to put in the Table 
\ref{rest0112}.

Since $Z_1$ and $Z_2$ correspond to stable maps with smooth domain curves, 
they do not lie on $D$.  Thus 
restricting the class 
$D$ to them gives zero.  The other four fixed points all lie on
$D$. Restricting $D$ to one of these points is equivalent to taking the first
Chern class of the restriction of $\O(D)$ to that point. 
In what follows, we will factor this
restriction
into two steps, first restricting to $D$ and then restricting
to the fixed point.
Let $(C,x_1,\ldots,x_n,f)\in\mm$ be a stable map.
In case $D$ is smooth at $f$, $\O(D)$ restricts to the normal bundle
$N_{D/\mm}$ on $D$. Restricting this further to $f$ gives a one-dimensional
vector space which is a quotient of the tangent space $T_{\mm,f}$ of
$\mm$ at $f$.

 The tangent space  $T_{\mm,f}$
is the same as the space of first order deformations of the stable
map $f$. Each such deformation is a combination of three basic types: 
Deformations that
deform the map while preserving the curve and its marked points,
deformations that move the marked points and nodes while preserving the 
nodes, and deformations that smooth the nodes. 
We are interested in isolating the deformations that smooth the nodes.
This is because smoothing the node of a generic stable map in $D$
takes one outside 
the divisor $D$, and the space of deformations that do this corresponds to 
the normal space of $D$ at this point.
Toward this end, we mention that
there is a surjection
\[T_{\mm,f}\lra H^0(C,\und{\Ext}^1(\Omega_C,\O_C))\text{,}\]
where $H^0(C,\und{\Ext}^1(\Omega_C,\O_C))$ can be identified with the
node-smoothing deformations. 

We can describe the 
deformations that smooth the nodes as follows. Let $p_1,\ldots,p_m$ be
 the nodes of $C$. For each $i$, let $C_1^i$ and $C_2^i$ be the
components of $C$ intersecting at the node $p_i$. 
We make use of the following fact
stated by Kontsevich in \cite{Ko}. (See \cite[\S 3B]{HM} 
for more detail.)

\begin{lem}
There is a natural isomorphism

\begin{equation}\label{eq:sumtan}
H^0(C,\und{\Ext}^1(\Omega_C,\O_C))\iso\+_{i=1}^m(T_{p_i}C_1^i\*T_{p_i}C_2^i)
\end{equation}
respecting the natural $T$-actions.
\end{lem}

Now each node is mapped to
one of the fixed points $q_0$, $q_1\in\PP^1$. The weight of tangent space
to a component on which $f$ has degree one at such a special point 
is inherited directly by pullback from the 
weight of the tangent space to the image of the node.  We know from 
\cite[Chapter 9]{CK} that the equivariant Euler class of the
tangent space to $\PP^r$ at $q_i$ is 
\[\prod_{j\neq i} (\la_i-\la_j)\text{.}\]
Thus in our case smoothing a node 
will contribute a weight (or sum of two weights) of
the form $\la_0-\la_1$ or $\la_1-\la_0$. 
More generally, the node smoothing from any singleton fixed component will 
contribute weight (or sum of two weights)
of the form $\w_F$, in the notation of Section \ref{sec:eq}.

Since $T$ acts equivariantly on the vector spaces involved, the weight of
the quotient space $i_{Z_j}^*(\O(D))$ will be one of the weights of the
tangent space to $Z_j$. The products of these weights were computed above
as the equivariant Euler classes. 
We simply pick the weight that corresponds to smoothing the node as 
described above, since this is the 
weight associated to the restriction of the normal
bundle.
For example,
at $Z_3$ we
obtain the weight $2(\la_0-\la_1)$ from smoothing the node. 
 
Now at $Z_5$ and $Z_6$, $D$ is not smooth. 
This is actually no problem at all. 
Each of these fixed components is an $S_2$-quotient of a point.
Recall from Section 
\ref{sec:eq} that we simply carry out computations on this smooth variety
upstairs, and then account for the group quotient when integrating by
placing an extra factor of two in the denominator.
Since we noted in Section \ref{sec:eq} that
every $T$-fixed component in a moduli space $\m0$ is the quotient 
of a nonsingular
variety by a finite group, 
we can always utilize Expression (\ref{eq:sumtan}) 
in computing restrictions of boundary classes to fixed components that
do not intersect the corresponding boundary locus transversally.

Let $j_{Z_5}:\Spec \CC\ra D$  be the composition of the natural atlas for 
$Z_5$ and the inclusion of $Z_5$ into $D$.
It is a 2--1 morphism. It factors in a natural way 
through $\tilde{D}=\mm_{0,2}(\PP^1,1)\x_{\PP^1}\mm_{0,1}(\PP^1,1)$ 
via a 2--1 morphism 
$j:\Spec \CC\ra\tilde{D}$.
Let $\N_{D/\mm}$ be the normal sheaf of $D$ in $\mm$.
The two sheaves $\O(D)|_D$ and $\N_{D/\mm}$ are not isomorphic
on $D$, but their
pullbacks to $\tilde{D}$ {\em are} isomorphic bundles
since $\tilde{D}$ is smooth.
In other words, it now becomes clear which node of $Z_5$ to smooth: the one
which is glued under the gluing map $\tilde{D}\ra D$.
The collapsed component at this node 
contributes nothing to the smoothing weight
because it is fixed
by the $T$-action and thus has zero weight.
Hence the weight of $\O(D)$ pulled back to $\tilde{D}$
is again $\la_0-\la_1$, the 
weight that corresponds to smoothing the node of $D$.
Finally, pulling this back to $\Spec \CC$ gives $2(\la_0-\la_1)$ since
$j$ has degree two.
The computation for $Z_6$ is similar.
The results of all the restrictions are listed in Table 
\ref{rest0112}.


\renewcommand{\baselinestretch}{1}
\begin{table}
\begin{center}
\begin{tabular}{|c|c|c|c|c|c|c|}
\hline
Fixed component & $Z_1$ & $Z_2$ & $Z_3$ & $Z_4$ & $Z_5$ & $Z_6$
 \\ \hline 

$H_1$ & $\la_0$ & $\la_1$ &  $\la_1$ &  $\la_0$ &  $\la_0$ &  $\la_1$
 \\ \hline

$D$ & 0 & 0 & $2(\la_0-\la_1)$ & $2(\la_1-\la_0)$ 
& $2(\la_0-\la_1)$ & $2(\la_1-\la_0)$ \\ \hline

\end{tabular}
\caption{Restrictions of divisor classes in $A_T^*(\mm_{0,1}(\PP^1,2))$
to fixed components \label{rest0112}}
\end{center}
\end{table}
\renewcommand{\baselinestretch}{2}

We already know the relation $H_1^2=0$ since $H^2=0$, so the divisor classes
above give two degree three
classes, $D^3$ and $D^2H_1$. The restrictions of 
these degree three classes to the fixed components are given in
Table \ref{restdeg3}. These follow directly
from the divisor restrictions in Table \ref{rest0112}.

\renewcommand{\baselinestretch}{1}
\begin{table}
\begin{center}
\noindent
\begin{tabular}{|c|c|c|c|c|c|c|}
\hline
 & $Z_1$ & $Z_2$ & $Z_3$ & $Z_4$ & $Z_5$ & $Z_6$
 \\ \hline 

$D^3$ & 0 & 0 &  $8(\la_0-\la_1)^3$ &  $8(\la_1-\la_0)^3$ &  $8(\la_0-\la_1)^3$ &  
$8(\la_1-\la_0)^3$
 \\ \hline

$D^2H_1$ & 0 & 0 & $4\la_1(\la_0-\la_1)^2$ & $4\la_0(\la_1-\la_0)^2$ 
& $4\la_0(\la_0-\la_1)^2$ & $4\la_1(\la_1-\la_0)^2$ \\ \hline

\end{tabular}
\caption{Restrictions of degree 3 classes in $A_T^*(\mm_{0,1}(\PP^1,2))$
to fixed components
\label{restdeg3}}
\end{center}
\end{table}
\renewcommand{\baselinestretch}{2}

At last 
we are ready to compute the integrals of the degree three monomials. Note
that since components $Z_5$ and $Z_6$ have automorphism group $S_2$, we need to
put an extra factor of 2 in the denominators of integrands over these  
components when we apply localization. We get

\begin{eqnarray*}
\int_{\mm_T}D^3 
& = & \int_{(Z_{3})_T}\frac{8(\la_0-\la_1)^3}{2(\la_1-\la_0)^3}+
\int_{(Z_{4})_T}\frac{8(\la_1-\la_0)^3}{2(\la_0-\la_1)^3} \\
& = & +\int_{(Z_{5})_T}\frac{8(\la_0-\la_1)^3}{2(\la_0-\la_1)^3}
+\int_{(Z_{6})_T}\frac{8(\la_1-\la_0)^3}{2(\la_1-\la_0)^3} \\
& = & -4-4+4+4 = 0
\end{eqnarray*}
and

\begin{eqnarray*}
\int_{\mm_T}D^2H_1
& = & \int_{(Z_{3})_T}\frac{4\la_1(\la_0-\la_1)^2}{2(\la_1-\la_0)^3}+
\int_{(Z_{4})_T}\frac{4\la_0(\la_1-\la_0)^2}{2(\la_0-\la_1)^3} \\
& + & \int_{(Z_{5})_T}\frac{4\la_0(\la_0-\la_1)^2}{2(\la_0-\la_1)^3}
+\int_{(Z_{6})_T}\frac{4\la_1(\la_1-\la_0)^2}{2(\la_1-\la_0)^3} \\
& = & \frac{2\la_1}{(\la_1-\la_0)}+\frac{2\la_0}{(\la_0-\la_1)}
+\frac{2\la_0}{(\la_0-\la_1)}+\frac{2\la_1}{(\la_1-\la_0)} = 
2+2 = 4\text{.}
\end{eqnarray*}

Now we can construct a presentation for $A^*(\mm_{0,1}(\PP^1,2))$.
First we show that $D$ and $H_1$ are independent. Suppose we have a relation
\[aD+bH_1=0\text{,}\]
where $a,b\in\QQ$.  Multiplying both sides by $D^2$ gives $aD^3+bH_1D^2=0$.
Then integrating gives $4b=0$. Similarly, multiplying by $H_1D$
and integrating the result gives $4a=0$. Clearly $a=b=0$. Thus $D$ and 
$H_1$ generate $A^1(\mm)$.

Next we show that $DH_1$ and $D^2$ are independent. Suppose we have a relation
\[aDH_1+bD^2=0\text{,}\]
where $a,b\in\QQ$.  Multiplying both sides by $D$ and integrating gives 
$4a=0$.  Multiplying both sides by $H_1$ and integrating gives
$4b=0$.  We see once again that $b=a=0$.  Thus $DH_1$ and 
$D^2$ generate $A^2(\mm)$.

Finally, there must be a relation in $A^3(\mm)$.
In fact, since $\int D^3=0$ was computed above, 
we are already able to conclude that $D^3=0$ is the desired relation.
Nevertheless, we will compute this by linear algebra to further illustrate
the method.
Suppose we have a relation
\[aD^2H_1+bD^3=0\text{,}\]
where $a,b\in\QQ$. Integrating this gives $4a=0$, so that 
$a=0$. We can take $b=1$ to get the relation
$D^3=0$.
It is easy to see that $D^2H_1$ is not zero, so we have all 
generators and relations in degree three. 
Everything in higher degrees is zero,
so we can give a complete presentation
\[A^*(\mm_{0,1}(\PP^1,2))=\frac{\QQ[D,H_1]}{(H_1^2,D^3)}\text{.}\]
%

\renewcommand{\baselinestretch}{1}
\section{Expressions for $\psi$-classes in the above moduli spaces}
\label{sec:exppsi}

Below we express the $\psi$-classes in terms of the boundary
and hyperplane divisors. This is possible since we have seen that 
these divisor
classes generate the Chow rings in the cases under consideration. 
We will use the setup of 
Section \ref{sec:psi}.

To avoid confusion, we will use notation $\pi_n$ for universal projections to
a space of stable maps (or more generally for contraction morphisms that 
forget a marked point) and $\r_i:(\PP^1)^n\rightarrow\PP^1$ for projections
of $(\PP^1)^n$ to a factor.

We will make frequent use of several identities expressing pullbacks of
the standard divisor classes on $\m0$ under contraction morphisms in terms of
the standard divisor classes on $\mm_{0,n+1}(\PP^r,d)$.
The formula governing pullback of $\psi$-classes is 
$\psi_i=\pi_{n+1}^*(\psi_i)+D_{i,n+1}$, which is a well-known
extension of an identity in \cite{Wit}. On the
left side, $\psi_i\in A^*(\mm_{0,n+1}(\PP^r,d))$, while on the right side,
$\psi_i\in A^*(\m0)$ and $D_{i,n+1}$ is the divisor

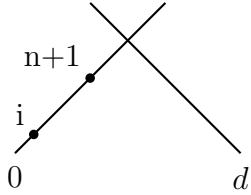
\begin{center}
\begin{pspicture}(0,0)(4,3)
\pnode(0.5,.5){a}
\dotnode(1.5,1.5){c}
\dotnode(.75,.75){d}
\pnode(2.5,2.5){e}
\pnode(1.5,2.5){f}
\pnode(3.5,0.5){j}
\ncline{a}{e}
\ncline{f}{j}
\uput{5pt}[d](.5,.5){0}
\uput{5pt}[d](3.5,.5){$d$}
\uput{5pt}[ul](.75,.75){i}
\uput{5pt}[ul](1.5,1.5){n+1}
\end{pspicture}
\end{center}
with the remaining marked points on the degree $d$ component. For $i<j$,
it's easy to see that $\ev_{i,n}\circ\pi_j=\ev_{i,n+1}$. Here the second
subscript of $\ev_i$ indicates the number of marked points associated to its
domain. It follows that $\pi_j^*(H_i)=H_i$ for $i<j$. A similar statement
holds when $i>j$ (and even when $i=j$ if the marked points are 
monotonically relabeled as in Section \ref{sec:modintro}), 
but attention must be given to how the indexing changes in the
wake of deleting one index.

Finally, we recall the basic fact about pullbacks of boundary divisors: 
If $\pi$ is any contraction morphism, then

\begin{equation}
\label{pulldiv}
\pi^*(D_{A,d_A,B,d_B})=\sum_{A\sub A\pr,B\sub B\pr} D_{A\pr,d_A,B\pr,d_B}
\text{.}
\end{equation}
Clearly the righthand side is the support of the pullback.
See \cite{FP}, for example, for relevant statements about why the pullback is
multiplicity-free.

\subsection{Expressions for $\psi$-classes in moduli spaces with $d=1$}

\subsubsection{The case $d=1$, $n=1$}
\label{sec:11}

Recall that $\mm_{0,1}(\PP^1,1)\iso\PP^1$ and 
its universal curve is $\mm_{0,2}(\PP^1,1)\iso\PP^1\x\PP^1$.
Furthermore,
the section is the diagonal map $\D$. The universal 
projection is the projection $\r_1$.
Then 

\begin{eqnarray*}
\w_{\r_1} = K_{\PP^1\x\PP^1}\*\r_1^*K_{\PP^1}^\vee
=\r_1^*(\O(-2))\*\r_2^*(\O(-2))\*\r_1^*(\O(-2))^\vee=\O(0,-2)
\text{.}\end{eqnarray*}
Now $c_1(\O(0,-2))=-2H_2$, and the pullback of each $H_i$ under $\D$ is $H_1$.
So in this case $\psi=\psi_1=\D^*(c_1(\O(0,-2)))=-2H_1$.

\subsubsection{The case $d=1$, $n=2$}

Using the results stated above, in $A^*(\mm_{0,2}(\PP^1,1))$ we have
\[\psi_1=\r_1^*(-2H_1)+D_{1,2}=-2H_1+H_1+H_2=H_2-H_1\text{.}\]
Here $D_{1,2}$ is the divisor corresponding to the locus where the marked
points have the same image. This is the diagonal of $\PP^1\x\PP^1$, and we
have used the fact that its class in $A^*(\PP^1\x\PP^1)$ is $H_1+H_2$.

By symmetry, that is, by forgetting the first marked point via the
projection $\r_2$ instead of the second using $\r_1$,
we find
\[\psi_2=H_1-H_2\text{.}\]

\subsubsection{The case $d=1$, $n=3$}
\label{sec:13}

Pulling back from $\mm_{0,2}(\PP^1,1)$, we have
\[\psi_1=\pi_3^*(\psi_1)+D_{1,3}=H_2-H_1+H_1+H_3-D=H_2+H_3-D
\text{,}\]
\[\psi_2=\pi_3^*(\psi_2)+D_{2,3}=H_1-H_2+H_2+H_3-D=H_1+H_3-D
\text{,}\]
and, by symmetry,
\[\psi_3=H_1+H_2-D
\text{.}\]
Here $D_{i,j}$ is the divisor corresponding to the closure of locus where the 
$i$'th and $j$'th marked
points have the same image, but the image of the remaining marked point
is different. This is the 
proper transform in $\Bl_\D{(\PP^1)^3}$ of the large
diagonal $\D_{ij}$ in $(\PP^1)^3$.
Its class in $A^*(\Bl_\D{(\PP^1)^3})$ is $H_i+H_j-D$.

\subsection{Expression for the $\psi$-class in the case $d=2$, $n=1$}

In this subsection let $\mm=\mm_{0,1}(\PP^1,2)$. We know that
$A^1(\mm)=\Div(\mm)$ is
generated by the boundary divisor class $D$ and the hyperplane
pullback class $H_1$. So $\psi=aD+bH_1$ for some $a,b\in\QQ$.  We use
a method for finding relations developed by
Pandharipande in \cite{P} to determine $a$ and
$b$.  Let $C$ be the curve in $\mm$ corresponding to the family

\vspace{0.3in}

\begin{center}
\psset{arrows=->}
\begin{psmatrix}
$\PP^1\times\mathbb{P}^1$ & $\mathbb{P}^1$ & $\mathbb{P}^1$ \\ 
$\mathbb{P}^1$
\ncline[offset=-3pt]{1,1}{2,1}\naput{$\r_1$}
\ncarc[arcangle=45]{2,1}{1,1}\naput{$\D$}
\ncline{1,1}{1,2}\naput{$\r_2$}
\ncline{1,2}{1,3}\naput{$f$}
\end{psmatrix}
\psset{arrows=-}
\end{center}

\vspace{0.1in}

\noindent
of stable maps, where $f$ is a fixed double cover of $\PP^1$ by
itself. Since $C$ is contained in the substack $\nobarmm$ 
of $\mm$ corresponding to
stable maps with smooth domain curves,
$D\cdot C=0$.  Let $i:C\ra\mm$ be the inclusion. Next, 
by definition the restriction $i_C^*\psi$ of
$\psi$ to $C$ is the first Chern class of the pullback of the
sheaf of relative differentials of the family under the section
$\D$. (See Section \ref{sec:psi}.)
But we computed exactly this in Section \ref{sec:11}, finding it to be
$-2H_1$. Finally, the degree of $C\cdot H_1$ is two, 
because the map has degree two: two points lie
above a general point in $\PP^1$.  So
\[2b=\int C\cdot bH_1=\int C\cdot\psi=\int i_C^*\psi=-2\text{,}\]
and we conclude that $b=-1$.

Now we construct another curve $C\pr$ which intersects $D$ but not $H_1$. 
(Here we consider $H_1$ as the pullback of a fixed hyperplane.)
As above, we give a family of curves together with a map from the total
space to $\PP^1$. We
can avoid $H_1$ by arranging for the marked point to have the same 
image under this map on each fiber.  In order to cause $C\pr$ to intersect 
$D$, we will initially specify a rational
map with base points. Blowing up to eliminate
the base points gives rise to some reducible fibers.
Let $\pi=\r_1:\PP^1\x\PP^1\rightarrow\PP^1$ be our 
initial family of curves. 
Let 
$t$ be a coordinate on the base, considering this $\PP^1$ as the one-point
compactification of $\CC$, and let $(x:y)$ be homogeneous coordinates on the
other factor. Associate to this family the constant section $s$ given 
by $s(t)=(t,(0:1))$.
Define $f:\PP^1\x\PP^1\rightarrow\PP^1$ by
$f(t,(x:y))=((t+1)x^2+ty^2:tx^2+(t+1)xy)$.  

Base points will occur whenever there is a common factor in the two
components. There are two ways this can happen, corresponding to the two 
factors $x$ and $tx+(t+1)y$ of the second component.
In the fiber over $t=0$, the map
is given by $(x^2,xy)$, which has a base point at $(0:1)$ since
the first component also has a factor of $x$.  This is the only base point
of this type. The existence of three base 
points of the second type is more subtle. These arise when, for a certain $t$,
$tx+(t+1)y$ divides $(t+1)x^2+ty^2$. For simplicity, we may assume $y=1$.
A base point will occur where $tx+(t+1)$ and $(t+1)x^2+t$ simultaneously 
vanish.  Performing some algebra, we see that this happens when
$2t^3+3t^2+3t+1=0$.  This has solutions $-\frac{1}{2}$, $e^{i\pi/3}$, and 
$e^{i2\pi/3}$. The corresponding base points are $(-\frac{1}{2},(1:1))$,
$(e^{i\pi/3},(-e^{-i\pi/6}:1))$, and $(e^{i2\pi/3},(e^{i\pi/6}:1))$. 

Notice 
that the section $s$ 
passes only through the base point $(0,(0:1))$ and not the 
other three. 
Also note that 
$f(t,(0:1))=(1:0)$ for all $t\neq 0$, so that the image of the marked point
is $(1:0)$ on all fibers other than the fiber over $t=0$.

To arrive at the family $C\pr$ of stable maps, we blow up at the four
base points.  Let $S$ be the resulting surface and 
$\r:S\rightarrow\PP^1\x\PP^1$ the blowup map.  We take for sections
$s_i:\PP^1\ra S$ the proper transforms of the original sections $s_i$.
By abuse of notation, keep the
labels $\pi$ and $f$ for the induced projection and map to the target as well.

It is clear that 
$\r$ is an isomorphism of curves on the fibers over points other than the four
base points, and that the fibers over 
the special values of $t$ now each have two
components. 
We will check that on such a special fiber, each the restriction of $f$ has
degree one on each component. Over $t=0$, the marked point is on the 
new component
coming from the exceptional divisor since its original image was blown up. 
In the other three special fibers, the marked points avoided the blowup, and
thus remain on the original component of those curves.

We start with the fiber over $t=0$.
By fixing $y=1$, we consider the affine piece of $\PP^1\x\PP^1$ with
coordinates $t$ and $x$. The blowup introduces additional coordinates 
$u$ and $v$ subject to the blowup equation $tv=ux$. Looking locally
on the blowup, we assume $v=1$ so that $t=ux$. The map $f$ becomes
\[((ux+1)x^2+ux:ux^3+(ux+1)x)=((ux+1)x+u:ux^2+ux+1)\text{.}\]
The exceptional divisor $E$ corresponds to $x=0$. Hence $f$ restricted to the
exceptional divisor is $(u:1)$.  Similarly, the other component in the
fiber $t=0$ corresponds to $u=0$, and $f$ restricted to that component is
given by $(x:1)$. Thus the map has degree 1 on each component of the fiber.

Second, we consider the blowup at the point $(-\frac{1}{2},1)$.
Here we use the coordinate system $((w,z),(u:v))$, where $w=t+\frac{1}{2}$,
$z=x-1$, and $wv=uz$. Take $v=1$ to arrive in local
coordinates, where now $w=uz$.
We find

\begin{eqnarray*}
 & & ((t+1)x^2+t:tx^2+(t+1)x) \\ 
& = & ((uz+\frac{1}{2})(z+1)^2+uz-\frac{1}{2}:
(uz-\frac{1}{2})(z+1)^2+(uz+\frac{1}{2})(z+1)) \\
& = & {\scriptstyle(uz^3+\frac{1}{2}z^2+2uz^2+uz+z+\frac{1}{2}+uz-\frac{1}{2}:
uz^3+2uz^2+uz-\frac{1}{2}z^2-z-\frac{1}{2}+uz^2+\frac{1}{2}z+uz+\frac{1}{2})}
\\
& = & (uz^2+\frac{1}{2}z+2uz+2u+1:uz^2+3uz+2u-\frac{1}{2}z-\frac{1}{2})
\text{.}
\end{eqnarray*}
The fiber over $t=-\frac{1}{2}$ has two components.
On the component corresponding to the exceptional divisor, $z=0$, so $f$
restricts to $(2u+1:2u-\frac{1}{2})$. On the other component, $u=0$, so
$f$ restricts to $(\frac{1}{2}z+1:-\frac{1}{2}z-\frac{1}{2})$. We see
that the map has degree 1 on each component. The computations on the other
two special fibers are similar.

Let us return to considering the fiber over $t=0$.
The value of $s$ at $t=0$ is determined by its values at nearby points
in $\PP^1\x\PP^1$. In local coordinates $s(t)=(t,0)$.  Switching
to coordinates $((t,x),(u:v))$ 
on the blowup, we find 
\[\lim_{t\ra 0}s(t)=\lim_{t\ra 0}((t,0),(u:0))
=\lim_{t\ra 0}((t,0),(1:0))=((0,0),(1:0))\text{,}\]
where $v=0$ because of the blowup equation $tv=ux$.
Thus $s(0)$ has value $(1:0)$ on the exceptional
divisor $E$ above $t=0$. 
The restriction of $f$ to $E$ described in local coordinates above extends
to $f((0,0),(1:0))=(1:0)$. This completes the verification that the marked
point of every stable map in this family has image $(1:0)$.
We have shown that $\int C\pr\cdot D=4$, and
we have also constructed $C\pr$ so that $\int C\pr\cdot H_1=0$.  
Furthermore, by standard intersection theory on surfaces (see 
\cite[Ch.V,\S 3]{H}), $E.(s(\PP^1))=1$, and the image of the section does
not intersect the other exceptional divisors. Finally, if $E_1$, $E_2$, and
$E_3$ are the other exceptional divisors, then
\[s^*(c_1(\w_{\pi}))=s^*(-2H_2+E+E_1+E_2+E_3)=s^*E\]
by standard results about dualizing sheaves of products and blowups and
about intersection theory on blowups (\cite{H}).
Therefore
\[4a=\int C\pr\cdot aD=\int C\pr\cdot\psi=\int i_{C\pr}^*\psi=
\int s^*(c_1(\w_{\pi}))=\int s^*(E)=1\text{,}\]
so that $a=\frac{1}{4}$.
We can now conclude
that
\[\psi=\frac{1}{4}D-H_1\]
in $\mm_{0,1}(\PP^1,2)$.

\chapter{Relations}
\label{sec:rel}

Relations come from three different sources: Pullbacks of relations
in $A^*(\mm_{0,1}(\PP^r,2))$ via the contractions forgetting a marked
point, relations given by the geometry of $\psi$-classes
on boundary divisors, and one linear
relation that so far has no {\em proven} geometric explanation. Instead,
we prove this last relation using the method of localization and linear
algebra, which will be explained in Section \ref{sec:lla}.

\section{Relations from pullbacks}
\label{sec:pb}

In Section \ref{sec:0112}, we found the relations $H_1^2$ and $D^3$ in
$A^*(\mm_{0,1}(\PP^1,2))$. Additionally, Section \ref{sec:exppsi} gives
the relation $\psi-\frac{1}{4}D+H_1$.

There are two contraction morphisms $\pi_1$ and $\pi_2$
from $\mm_{0,2}(\PP^1,2)$ to
$\mm_{0,1}(\PP^1,2)$. Recall that $\pi_i$ forgets the $i$'th marked point. 
Formulas for the pullbacks of the standard divisor classes under these
maps were given in Section \ref{sec:exppsi}.

First, pulling back the relation $H_1^2$ under $\pi_1$ and $\pi_2$
gives the relations $H_i^2$ for $i\in\und{2}$. (When considering the 
contraction $\pi_1$, either a monotonic relabeling of the marked points 
gives an index shift on pullback, or $H_1$ needs to be labeled as $H_2$.)
We can also see these relations via pullback under the evaluation morphisms.
Since $H^2=0$ in $A^*(\PP^1)$,
we have $H_i^2=\ev_i^*(H^2)=0$.

Second, applying Equation \ref{pulldiv} to the case at hand, 
we have $\pi_i^*(D)=D_1+D_2$ for $i\in\und{2}$.
Pulling back the relation $D^3$ in $A^*(\mm_{0,1}(\PP^1,2))$ by either one of 
these gives the cubic relation $(D_1+D_2)^3$ in $A^*(\mm_{0,2}(\PP^1,2))$.

Finally, 
linear relations expressing the $\psi$-classes in terms of the other divisor
classes are obtained by pulling back the relation $\psi-\frac{1}{4}D+H_1$.
We find
\[\psi_i=\pi_j^*(\psi_i)+D_0=\pi_j^*(\frac{1}{4}D-H_1)+D_0
=\frac{1}{4}D_1+\frac{1}{4}D_2+D_0-H_i\text{,}\]
where $i\neq j$. Thus we have relations
$\psi_i-\frac{1}{4}D_1-\frac{1}{4}D_2-D_0+H_i$.

\renewcommand{\baselinestretch}{1}
\section{Relations from the geometry of $\psi$-classes on boundary 
divisors}
\label{sec:geomrel}

We can interpret a product with
a factor of $D_i$ as a
restriction to that divisor and, via the gluing morphisms described
in Section \ref{sec:modintro}, ultimately as a class on
a fiber product of simpler 
moduli spaces.  This simplifies the computations of such products
once we know the pullbacks of the remaining factors.
We will use the gluing morphisms
\[j_0:\mm_{0,3}(\PP^1,0)\x_{\PP^1}\mm_{0,1}(\PP^1,2)
\rightarrow\mm_{0,2}(\PP^1,2)\]
and
\[j_1:\mm_{0,3}(\PP^1,1)\x_{\PP^1}\mm_{0,1}(\PP^1,1)
\rightarrow\mm_{0,2}(\PP^1,2)\text{,}\]
whose images are $D_0$ and $D_1$ respectively. By convention, 
both $j_0$ and $j_1$ glue the third
marked point of the first factor to the lone marked point of the second
factor. It is not hard to see that
$j_0$ is an isomorphism onto $D_0$, and we note further that
$\mm_{0,3}(\PP^1,0)\x_{\PP^1}\mm_{0,1}(\PP^1,2)\iso\mm_{0,1}(\PP^1,2)$.
Similarly, $\mm_{0,3}(\PP^1,1)\x_{\PP^1}\mm_{0,1}(\PP^1,1))
\iso\mm_{0,3}(\PP^1,1)$. However, $j_1$ is only an isomorphism away from
the divisor $D\x_{\PP^1}\mm_{0,1}(\PP^1,1)$, where $D$ is the boundary
divisor of $\mm_{0,3}(\PP^1,1)$ as described in Section \ref{sec:0n11}.
The image of this divisor is isomorphic to the global quotient stack 
\[[\mm_{0,1}(\PP^1,1)\x_{\PP^1}\mm_{0,4}(\PP^1,0)\x_{\PP^1}\mm_{0,1}(\PP^1,1)
/S_2]\text{,}\]
where the $S_2$-action switches the factors on the ends.
Thus the restriction of $j_0$ to $D$ has degree two.
Note that this last fiber product is isomorphic to 
$[\mm_{0,4}(\PP^1,0)/S_2]$,
where the $S_2$ action switches the third and fourth marked points.
 
For the $D_0$ case, the universal property of the moduli space shows that
$\psi_1$ and $\psi_2$ pull back to what may be considered as the first and
second $\psi$-classes on $\mm_{0,3}(\PP^1,0)$.
A similar statement holds for $D_1$ case,  
this time with the resulting $\psi$-classes on $\mm_{0,3}(\PP^1,1)$.

In the first case using this technique, we will show that the $\psi$-classes
vanish on $D_0$. This is because each 
$\psi$-class pulls back to zero under $j_0$ since the marked points lie
on the rigid component corresponding to $\mm_{0,3}(\PP^1,0)$. 
More rigorously, $j_0$ induces a family 
$(\C,\tilde{\s}_1,\tilde{\s}_2,\ev_3\circ\tilde{\jmath}_0)$ of
stable maps via the fiber diagram

\vspace{0.3in}

\psset{arrows=->}
\begin{center}
\begin{psmatrix}
$\C$  & $\mm_{0,3}(\PP^1,2)$ & $\PP^1$\\ 
$\mm_{0,3}(\PP^1,0)\x_{\PP^1}\mm_{0,1}(\PP^1,2)$ & 
$\mm_{0,2}(\PP^1,2)\text{,}$
\ncline{1,1}{2,1}\naput{$\tilde{\pi}$}
\ncline{1,1}{1,2}\naput{$\tilde{\jmath}_0$}
\ncline{1,2}{2,2}\naput{$\pi$}
\ncline{2,1}{2,2}\naput{$j_0$}
\ncline{1,2}{1,3}\naput{$\ev_3$}
\ncarc[arcangle=45]{2,1}{1,1}\naput{$\tilde{\s}_i$}
\ncarc[arcangle=45]{2,2}{1,2}\naput{$\s_i$}
\end{psmatrix}
\end{center}
\psset{arrows=-}

\vspace{0.2in}

\noindent
It is easy
to check that this family can be identified with a universal family of
stable maps over $D_0$. We can obtain another family over 
$\mm_{0,3}(\PP^1,0)\x_{\PP^1}\mm_{0,1}(\PP^1,2)$ by fiber product:

\vspace{0.3in}

\psset{arrows=->}
\begin{center}
\begin{psmatrix}
$\mm_{0,4}(\PP^1,0)\x_{\PP^1}\mm_{0,1}(\PP^1,2)$  & 
$\mm_{0,4}(\PP^1,0)$ & $\PP^1$ \\ 
$\mm_{0,3}(\PP^1,0)\x_{\PP^1}\mm_{0,1}(\PP^1,2)$ & 
$\mm_{0,3}(\PP^1,0)\text{,}$
\ncline{1,1}{2,1}\naput{$\tilde{\pi}_4$}
\ncline{1,1}{1,2}\naput{$\tilde{\pj}_1$}
\ncline{1,2}{2,2}\naput{$\pi_4$}
\ncline{2,1}{2,2}\naput{$\pj_1$}
\ncline{1,2}{1,3}\naput{$\ev_4$}
\ncarc[arcangle=45]{2,1}{1,1}\naput{$\tilde{s}_i$}
\ncarc[arcangle=45]{2,2}{1,2}\naput{$s_i$}
\end{psmatrix}
\end{center}
\psset{arrows=-}

\vspace{0.2in}

\noindent
where the right-hand part is the universal stable map.
Recall from Section \ref{sec:0n11} that the marked points in
$\mm_{0,3}(\PP^1,0)$ cannot vary.
It follows that $s_i^*(\w_{\pi_4})=\psi_i=0$ for 
$i\in\und{3}$. Note also that $\w_{\pi_4}$
pulls back to the relative dualizing sheaf of the left-hand
column.
There is an inclusion 
\[k_0:\mm_{0,4}(\PP^1,0)\x_{\PP^1}\mm_{0,1}(\PP^1,2)\ra
\mm_{0,3}(\PP^1,2)\]
very similar to the morphism $j_0$ described above. This morphism is part
of a 2-commutative diagram that guarantees the existence of a stack morphism
\[\iota:\mm_{0,4}(\PP^1,0)\x_{\PP^1}\mm_{0,1}(\PP^1,2)\ra\C\]
over $\mm_{0,3}(\PP^1,0)\x_{\PP^1}\mm_{0,1}(\PP^1,2)$. Furthermore,
$\iota$ is injective since $k_0$ is, and 
it is also compatible with the sections.
The images
of the sections $\tilde{\s_i}$ are contained in the image of $\iota$.
Therefore we can compute the $\psi$-classes of
$\mm_{0,3}(\PP^1,0)\x_{\PP^1}\mm_{0,1}(\PP^1,2)$ on this subfamily:
\begin{eqnarray*}
\psi_i & = & c_1(\tilde{\s}_i^*\w_{\tilde{\pi}})\\
 & = & c_1(\tilde{s}_i^*\w_{\tilde{\pi}_4})\\
 & = & \tilde{s}_i^*\tilde{\pj}_1^*(c_1(\w_{\pi_4}))\\
 & = & \pj_1^*s_i^*(c_1(\w_{\pi_4}))\\
 & = & 0\text{.}
\end{eqnarray*}
It follows that the pullbacks of $\psi_1$ and $\psi_2$
to $D_0$ vanish. This gives relations $D_0\psi_i$.

Now we will show that the product $D_1\psi_1\psi_2$ vanishes by computing the
pullback of $\psi_1\psi_2$ under $j_1$. 
Identifying $A^*(\mm_{0,3}(\PP^1,1)\x_{\PP^1}\mm_{0,1}(\PP^1,1))$ with
$A^*(\mm_{0,3}(\PP^1,1))$, the pullback is
\[\psi_1\psi_2=(H_2+H_3-D)(H_1+H_3-D)
=(H_1+H_2-D)(H_2+H_3-D)=0\]
according to the presentation given in Proposition \ref{m0311}. 
This shows that $\psi_1\psi_2$ is in the kernel of $j_1^*$.
It follows from \cite[Theorem 3.1]{Mu} that there is an group isomorphism
$A^2(D_1)\ra A^2(\mm_{0,3}(\PP^1,1))^{S_2}$, where the target is the subgroup
of $S_2$-invariants. This isomorphism is naturally identified with a
monomorphism into $A^2(\mm_{0,3}(\PP^1,1))$.
Furthermore, Mumford's proof of this theorem 
shows that, up to sign, this monomorphism
is the same as $\tilde{\jmath}_1^*$, where $\tilde{\jmath}_1$ is 
$j_1$ with the target
changed to $D_1$ (and considered as a map of these 
Chow groups). This is enough
to show that $\psi_1\psi_2$ restricts to zero on $D_1$.
Thus $D_1\psi_1\psi_2$ is a relation in $A^*(\mm_{0,2}(\PP^1,2))$.

\section{A relation from localization and linear algebra}
\label{sec:lla}

Ideally, the proof that (\ref{geomprez}) is a presentation for 
$A^*(\mm_{0,2}(\PP^1,2))$ should be uniformly geometric. Such a proof
would be easy to understand and generalize to higher-dimensional targets
$\PP^r$ (and other moduli spaces). 
Unfortunately, a sort of ``brute force" algebraic computation is
so far still required for two aspects of the proof. The first is the
derivation of the linear relation $D_2-\psi_1-\psi_2$. The second, closely
related, is the verification that the set of generators and relations is
complete in each degree.

The relation $D_2-\psi_1-\psi_2$ should be explained geometrically by
the existence of a section $s$ of the tensor product $L_1\* L_2$ 
of the cotangent line bundles whose zero
stack is $D_2$. Then
\[D_2=Z(s)=c_1(L_1\* L_2)=c_1(L_1)+c_1(L_2)=\psi_1+\psi_2\text{.}\]
The completeness of the set of generators and relations could be demonstrated
by giving an additive basis for $A^*(\mm_{0,2}(\PP^1,2))$. Such a basis
should be attainable 
using the stratification of the moduli space of stable maps
described in Section \ref{sec:ser}.

We call the computational algebraic method developed in this subsection
{\em localization and linear algebra}. In theory, it could be used to compute
relations in any moduli space of stable maps to projective space. (In 
practice, the method quickly becomes too tedious as the parameters increase.
See Section \ref{sec:direct} for an example.)
Its first step consists of using
localization (see Section \ref{sec:eq}) to
find the integrals of all degree four monomials in the boundary divisors and
hyperplane pullbacks in $A^*(\mm_{0,2}(\PP^1,2))$.
Then an arbitrary relation of given degree in these generating classes is
considered (with variable coefficients). 
This relation is multiplied by various monomials of complementary
dimension, and the result is integrated. The resulting rational polynomials 
in the coefficients place restrictions on these coefficients. Ultimately,
we can use these restrictions to describe precisely the form of 
relations that can occur in each degree in this Chow ring.

We begin by computing the integrals of degree four monomials. Notice first
that
the relations $D_0\psi_i$ and $\psi_i-\frac{1}{4}D_1-\frac{1}{4}D_2-D_0+H_i$
imply 
\begin{equation}\label{d0hi}
D_0H_1=D_0H_2\text{,}
\end{equation} 
and hence $H_1H_2D_0=H_1^2D_0=0$.
The integrals of the following types of monomials will be zero since the 
monomials themselves are zero:

\begin{enumerate}

\item Any monomial with a factor of $H_i^2$ for $i\in\underline{2}$.

\item Any monomial with a factor of $H_1H_2D_0$.

\end{enumerate}

These together take care of $2\left(\left({5 \atop 1}\right)+
\left({5 \atop 2}\right)\right)-1+3=32$ of the $\left({8\atop 4}\right)=70$ 
degree 4 monomials. We will use localization to compute the remaining 38 
integrals.   In particular, we apply Corollary \ref{stackloc} to
the smooth stack $\mm=\mm_{0,2}(\PP^1,2)$. Let $T=(\CC^*)^2$. 
Consider the usual $T$-action on $\mm$.  First we have to find the 
$T$-equivariant 
Euler classes of the normal bundles of the fixed point components.
 We do this using Theorem \ref{norm}. We will label the graphs 
of the 14 fixed components as follows.

\[\G_1=\text{
\psset{labelsep=2pt, tnpos=b,radius=2pt}
\pstree[treemode=R]{\TC*~{0}~[tnpos=a]{\{1\}}}
{
\TC*~{1}~[tnpos=a]{\{2\}}\taput{2} 
}
}\]

\[\G_2=\text{
\psset{labelsep=2pt, tnpos=b,radius=2pt}
\pstree[treemode=R]{\TC*~{0}~[tnpos=a]{\{2\}}}
{
\TC*~{1}~[tnpos=a]{\{1\}}\taput{2} 
}
}\]

\[\G_3=\text{
\psset{labelsep=2pt, tnpos=b,radius=2pt}
\pstree[treemode=R]{\TC*~{0}~[tnpos=a]{\underline{2}}}
{
\TC*~{1}\taput{2} 
}
}\]

\[\G_4=\text{
\psset{labelsep=2pt, tnpos=b,radius=2pt}
\pstree[treemode=R]{\TC*~{0}}
{
\TC*~{1}~[tnpos=a]{\underline{2}}\taput{2} 
}
}\]

\[\G_5=\text{
\psset{labelsep=2pt, tnpos=b,radius=2pt}
\pstree[treemode=R]{\TC*~{1}~[tnpos=a]{$\underline{2}$}}
{
\pstree{\TC*~{0}\taput{1}}
{
\TC*~{1}\taput{1} 
}
}
}\]

\[\G_6=\text{
\psset{labelsep=2pt, tnpos=b,radius=2pt}
\pstree[treemode=R]{\TC*~{0}~[tnpos=a]{$\underline{2}$}}
{
\pstree{\TC*~{1}\taput{1}}
{
\TC*~{0}\taput{1} 
}
}
}\]

\[\G_7=\text{
\psset{labelsep=2pt, tnpos=b,radius=2pt}
\pstree[treemode=R]{\TC*~{1}~[tnpos=a]{\{1\}}}
{
\pstree{\TC*~{0}~[tnpos=a]{\{2\}}\taput{1}}
{
\TC*~{1}\taput{1} 
}
}
}\]

\[\G_8=\text{
\psset{labelsep=2pt, tnpos=b,radius=2pt}
\pstree[treemode=R]{\TC*~{0}~[tnpos=a]{\{1\}}}
{
\pstree{\TC*~{1}~[tnpos=a]{\{2\}}\taput{1}}
{
\TC*~{0}\taput{1} 
}
}
}\]

\[\G_9=\text{
\psset{labelsep=2pt, tnpos=b,radius=2pt}
\pstree[treemode=R]{\TC*~{1}~[tnpos=a]{\{2\}}}
{
\pstree{\TC*~{0}~[tnpos=a]{\{1\}}\taput{1}}
{
\TC*~{1}\taput{1} 
}
}
}\]

\[\G_{10}=\text{
\psset{labelsep=2pt, tnpos=b,radius=2pt}
\pstree[treemode=R]{\TC*~{0}~[tnpos=a]{\{2\}}}
{
\pstree{\TC*~{1}~[tnpos=a]{\{1\}}\taput{1}}
{
\TC*~{0}\taput{1} 
}
}
}\]

\[\G_{11}=\text{
\psset{labelsep=2pt, tnpos=b,radius=2pt}
\pstree[treemode=R]{\TC*~{1}}
{
\pstree{\TC*~{0}~[tnpos=a]{$\underline{2}$}\taput{1}}
{
\TC*~{1}\taput{1} 
}
}
}\]

\[\G_{12}=\text{
\psset{labelsep=2pt, tnpos=b,radius=2pt}
\pstree[treemode=R]{\TC*~{0}}
{
\pstree{\TC*~{1}~[tnpos=a]{$\underline{2}$}\taput{1}}
{
\TC*~{0}\taput{1} 
}
}
}\]

\[\G_{13}=\text{
\psset{labelsep=2pt, tnpos=b,radius=2pt}
\pstree[treemode=R]{\TC*~{1}~[tnpos=a]{\{1\}}}
{
\pstree{\TC*~{0}\taput{1}}
{
\TC*~{1}~[tnpos=a]{\{2\}}\taput{1} 
}
}
}\]

\[\G_{14}=\text{
\psset{labelsep=2pt, tnpos=b,radius=2pt}
\pstree[treemode=R]{\TC*~{0}~[tnpos=a]{\{1\}}}
{
\pstree{\TC*~{1}\taput{1}}
{
\TC*~{0}~[tnpos=a]{\{2\}}\taput{1} 
}
}
}\]

\vspace{0.1in}

\noindent
Let $Z_i$ denote the fixed component corresponding to $\G_i$ for all $i$. 
All of the fixed components are points except for $Z_{11}$ and $Z_{12}$,
which are each isomorphic to the quotient $[\PP^1/S_2]$. Fixed components
$Z_1$, $Z_2$, $Z_3$, and $Z_4$ also have automorphism group $S_2$.

All
the degree 4 classes that are not automatically zero for the reasons above 
have factors $D_i$. Therefore they are supported on the boundary of
$\mm_{0,2}(\PP^1,2)$. Hence representatives of these classes 
will not intersect $Z_1$ or $Z_2$, which lie in the locus
$\nobarmm_{0,2}(\PP^1,2)$ where the domain curves are smooth.
Since restrictions of the relevant classes to $Z_1$ and $Z_2$ are thus zero,
their equivariant Euler classes are not needed for our computation. 

Recall from Section \ref{sec:0112} that
the term $e_F$ that appears in the formula for $e_{\G}^{\text{F}}$ in
Theorem \ref{norm} is zero on any fixed component which is a point.  
This was called a $\tilde{\psi}$-class in
Section \ref{sec:psi}. From the discussion there, we see that,
for the components that are isomorphic to $\MM_{0,4}\iso\PP^1$,
we can identify $e_F$ with a $\psi$-class on $\MM_{0,4}$.
We will argue that
the $\psi_i$ are all equal to the class of a point for $i\in\underline{4}$. We
denote this class by $\psi$, so that on these components $e_F=\psi$.
Note that $\psi^2=0$.

\begin{lem}\label{allpsiequal}
Let $H$ be the class of a point in $\MM_{0,4}\iso\PP^1$. Then $\psi_4=H$
in $A^*(\MM_{0,4})$.
\end{lem}

\noindent
{\bf Proof.} 
We may fix the first three marked points to be 0, 1, and $\infty$.
The universal curve for $\MM_{0,4}$ is $\MM_{0,5}$, which is
isomorphic to the blowup $\Bl_3(\PP^1\x\PP^1)$ of $\PP^1\x\PP^1$ at the
points 0, 1, and $\infty$ on its diagonal. (This is a simple case of the
construction of Keel in \cite{K}, for example.) We abuse notation by 
identifying the forgetful morphism 
$\pi_5$ of the universal curve with
the projection $\r_2$. Similarly, 
the universal section $s$ corresponds to the diagonal
morphism. Let $E_1$, $E_2$, and $E_3$ be the exceptional divisors. Then, using
arguments like those at the end of Section \ref{sec:exppsi},

\begin{eqnarray*}
\psi_4 & = & s^*c_1(\w_{\pi_5})\\
	& = & s^*c_1(\r_1^*(\O(-2))\*\O(E_1)\*\O(E_2)\*\O(E_3))\\
	& = & s^*(-2H_1+E_1+E_2+E_3)\\
	& = & -2H+H+H+H \\
	& = & H\text{.} 
\end{eqnarray*}

\begin{flushright}
$\Box$\end{flushright}

\noindent
The same follows for the other $\psi_i$ by symmetry.

For $\G_3$ we label the edges and vertices as follows.

\begin{center}
\psset{labelsep=2pt, tnpos=b,radius=2pt}
\pstree[treemode=R]{\TC*~{A}~[tnpos=a]{\underline{2}}}
{
\TC*~{B}\taput{a} 
}
\end{center}

We calculate
\[e_{\G_{3}}^{\text{F}}=\frac{\w_{A,a}-e_{A,a}}{(\la_0-\la_1)(\la_1-\la_0)}
=\frac{(\la_0-\la_1)/2}{(\la_0-\la_1)(\la_1-\la_0)}=\frac{1}{2(\la_1-\la_0)}
\text{,}\]
\[e_{\G_{3}}^{\text{v}}=\frac{(\la_0-\la_1)(\la_1-\la_0)}
{(\la_1-\la_0)/2}=2(\la_0-\la_1)
\text{,}\]
and
\[e_{\G_{3}}^{\text{e}}=\frac{2^2(\la_1-\la_0)^4}{2^4}=\frac{(\la_1-\la_0)^4}{4}
\text{.}\]
Thus
\[\eqeul(N_{\G_{3}})=\frac{2(\la_0-\la_1)}{2(\la_1-\la_0)}\frac{(\la_1-\la_0)^4}{4}
=\frac{-(\la_1-\la_0)^4}{4}
\text{.}\]

A similar calculation shows that $\eqeul(N_{\G_{4}})=\frac{-(\la_1-\la_0)^4}{4}$
 as 
well. For the remaining graphs we label the vertices and edges as follows.

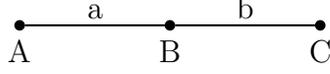
\begin{center}
\psset{labelsep=2pt, tnpos=b,radius=2pt}
\pstree[treemode=R]{\TC*~{A}}
{
\pstree{\TC*~{B}\taput{a}}
{
\TC*~{C}\taput{b} 
}
}
\end{center}

Then
\[e_{\G_{5}}^{\text{F}}=\frac{\w_{A,a}}{(\la_0-\la_1)^2(\la_1-\la_0)^2}
=\frac{\la_1-\la_0}{(\la_0-\la_1)^2(\la_1-\la_0)^2}=\frac{1}{(\la_1-\la_0)^3}
\text{,}\]
\[e_{\G_{5}}^{\text{v}}=(\la_1-\la_0)^2(\la_0-\la_1)\frac{\w_{B,a}+\w_{B,b}}
{\w_{C,b}}=\frac{(\la_1-\la_0)^2(\la_0-\la_1)2(\la_0-\la_1)}{\la_1-\la_0}
=2(\la_1-\la_0)^3
\text{,}\]
and
\[e_{\G_{5}}^{\text{e}}=(\la_1-\la_0)^4
\text{.}\]
Thus
\[\eqeul(N_{\G_{5}})=\frac{2(\la_1-\la_0)^3(\la_1-\la_0)^4}{(\la_1-\la_0)^3}
=2(\la_1-\la_0)^4
\text{.}\]

A similar calculation shows that $\eqeul(N_{\G_{6}})=2(\la_1-\la_0)^4$.
Next, for $Z_7$,
\[e_{\G_{7}}^{\text{F}}=\frac{\w_{B,a}\w_{B,b}}{(\la_1-\la_0)^2(\la_0-\la_1)^2}
=\frac{1}{(\la_1-\la_0)^2}
\text{,}\]
\[e_{\G_{7}}^{\text{v}}=\frac{(\la_1-\la_0)^2(\la_0-\la_1)}{\w_{C,b}}
=(\la_0-\la_1)(\la_1-\la_0)
\text{,}\]
and
\[e_{\G_{7}}^{\text{e}}=(\la_1-\la_0)^4
\text{.}\]
Thus
\[\eqeul(N_{\G_{7}})=\frac{(\la_0-\la_1)(\la_1-\la_0)(\la_1-\la_0)^4}{(\la_1-\la_0)^2}
=-(\la_0-\la_1)^4
\text{.}\]

Similar calculations show that $\eqeul(N_{\G_{8}})=\eqeul(N_{\G_{9}})=
\eqeul(N_{\G_{10}})=-(\la_0-\la_1)^4$ as well.

 For
$Z_{11}$ we get
\[e_{\G_{11}}^{\text{F}}=\frac{(\la_0-\la_1-\psi)^2}{(\la_0-\la_1)^4}
\text{,}\]

\[e_{\G_{11}}^{\text{v}}=\frac{(\la_1-\la_0)^2(\la_0-\la_1)}{(\la_1-\la_0)^2}
=\la_0-\la_1
\text{,}\]
and
\[e_{\G_{11}}^{\text{e}}=(\la_1-\la_0)^4
\text{.}\]
Thus
\[\eqeul(N_{\G_{11}})=(\la_0-\la_1)(\la_0-\la_1-\psi)^2
=(\la_0-\la_1)^2(\la_0-\la_1-2\psi)
\text{.}\]
A completely symmetric calculation shows that
\[\eqeul(N_{\G_{12}})=(\la_1-\la_0)^2(\la_1-\la_0-2\psi)
\text{.}\]

Finally,
\[e_{\G_{13}}^{\text{F}}=\frac{1}{(\la_0-\la_1)^2(\la_1-\la_0)^2}
\text{,}\]
\[e_{\G_{13}}^{\text{v}}=(\la_0-\la_1)^2(\la_1-\la_0)(\w_{B,a}+\w_{B,b})
=2(\la_0-\la_1)^2(\la_1-\la_0)^2
\text{,}\]
and
\[e_{\G_{13}}^{\text{e}}=(\la_1-\la_0)^4
\text{.}\]
Thus
\[\eqeul(N_{\G_{13}})=2(\la_0-\la_1)^4
\text{.}\]
A similar calculation shows that $\eqeul(N_{\G_{14}})=2(\la_0-\la_1)^4$ as well.

Next we need to know the restriction of each degree 4 monomial in the 
generating classes to each fixed component.  These come immediately from the
restrictions of the generating classes themselves to the fixed components,
which we now compute.
Results are given in the Table \ref{rest0212}.

The same computations were carried out for $\mm_{0,1}(\PP^1,2)$ in
Section \ref{sec:0112}, and much of the groundwork for the present case
was laid there. However, since the current space is more complicated,
we will need to extend the those remarks.  For example, there are
several boundary divisors to consider instead of just one. Perhaps most
strikingly, two of the fixed components are one-dimensional rather than
just points.
We again turn to the methods of \cite[\S 9.2]{CK}. As mentioned, some of the
fixed components are quotients of smooth varieties by finite groups.
We must be 
careful to perform these calculations on the variety before quotienting; the
integral formula of Corollary 
\ref{stackloc} takes this quotienting into account later.

The restrictions of the $H_i$ are computed exactly as in Section 
\ref{sec:0112}, except now there are two evaluation morphisms. Just use
$\ev_i$ for computing restrictions of $H_i$, so that the image of the $i$'th
marked point determines the weight of the restriction.

Many restrictions of boundary classes are zero simply because the 
corresponding boundary components are disjoint from the fixed components.
 We have $i_1^*(D_j)=i_2^*(D_j)=0$ for 
$j\in\{0,1,2\}$, because $Z_1$ and $Z_2$ have smooth domain curves and thus 
lie in the complement of the boundary.

The domain of every map 
in $D_1$ has two degree one components, as can be observed from
the degeneration diagrams in Section \ref{sec:bound}. The same holds for 
$D_2$.  Thus $i_3^*(D_1)=i_4^*(D_1)=i_3^*(D_2)=i_4^*(D_2)=0$ since the maps
of $Z_3$ and $Z_4$ do not have two degree one components.

Recall that the boundary divisor $D_0$ is the closure of the
locus consisting of stable
maps whose domains have a collapsed component and a degree two component. 
By stability, such a stable map must have both marked points on the collapsed
component.  Since we have seen that a collapsed component with three special
points is a rigid object, the property of possessing such a component is
preserved under limits. 
Thus the domain curve of every stable map lying in
$D_0$ must have both marked points on the same degree zero component.
On the other hand, the 
maps of $Z_7$, $Z_8$, $Z_9$, and $Z_{10}$ have their marked points on
different components.  
Thus $i_7^*(D_0)=i_8^*(D_0)=i_9^*(D_0)=i_{10}^*(D_0)=0$.
The above argument also shows that
the domain of every map in
$D_0$ has a component without any marked points. The same
holds for $D_1$.  Thus 
$i_{13}^*(D_0)=i_{14}^*(D_0)=i_{13}^*(D_1)=i_{14}^*(D_1)=0$ since the maps
of $Z_{13}$ and $Z_{14}$ each have a marked point on both components.

Similarly, the
domain curve of every
stable map in $D_2$ must have the marked points on distinct
components. Briefly, a generic stable map in $D_2$ has this property, and
in the limit marked points may only move to newly sprouted components.
(If both marked points collide with the node, {\em two} new components
will result.)
On the other hand, the 
maps of $Z_5$ and $Z_6$ have their marked points on the same component.
Thus $i_5^*(D_2)=i_6^*(D_2)=0$.

The simplest of the nonzero entries are the restrictions of $D_0$ and $D_2$
to $Z_{11}$ and $Z_{12}$. Each divisor intersects these fixed components
transversely before quotienting. The intersections consist of one point in
$\PP^1$ for $D_0$ and two points for $D_2$. Recall from Lemma \ref{allpsiequal}
that $\psi$ is the class
of a point in $\MM_{0,4}\iso\PP^1$. Thus we get $\psi$ for the restrictions
of $D_0$ and $2\psi$ for the restrictions of $D_2$.

The remaining entries are nonzero and are computed using Expression
(\ref{eq:sumtan}). This requires a bit of extra care. As before,
$i_j^*(D_i)=c_1(i_j^*\O(D_i))=c_1(i_j^*\N_{D_i/\mm})$ since we are working
on the smooth atlases over the fixed components. Equation (\ref{eq:sumtan}) 
gives
the results from smoothing all the nodes. The first Chern class of the 
restriction of the 
normal cone $\N_{D_i/\mm}$ to a fixed component
corresponds to smoothing a node in that degeneration locus that comes from
$D_i$. In particular,
this involves taking the sum of the weights of the tangent 
directions to each component at the ``node to be smoothed.'' Intuitively, this
is the node that, when smoothed, takes one outside the divisor being 
intersected with.  This makes sense because it gives the restriction of the 
normal cone of the 
divisor to its intersection with the fixed component. 

This is straightforward for $Z_3$ and $Z_4$ since these components have
only one node. Since the non-collapsed component of these stable maps have
degree two, the weights from the tangent space of $\PP^1$ are divided by two.
The weights of the restriction of $D_2$ to $Z_{13}$ and $Z_{14}$ are easy
to determine for the same reason.

Fixed components $Z_5$ and $Z_6$ have two nodes. When restricting $D_0$,
smooth the node contained in the collapsed component. When restricting $D_1$,
smooth the node that is the intersection of two degree one components.

Similarly, the curves associated to $Z_7$, $Z_8$, $Z_9$, and $Z_{10}$ have
two nodes.  When restricting $D_2$, smooth the node for which both irreducible
components have a marked point.  When restricting $D_1$, smooth the node
contained in a component without a marked point. At any rate, smoothing either
node gives the same weight for these components since both nodes map to the
same fixed point and lie at the same type of intersection.

Finally, we restrict $D_1$ to $Z_{11}$ and $Z_{12}$. Both of these components
lie completely inside of $D_1$, so that we still need to use the node 
splitting. Smoothing either node corresponds to the restriction of $\O(D_1)$.
For the first time, the collapsed components make a non-trivial contribution
since they contain four special points. Since $\psi$ is the first Chern
class of the cotangent bundle to any marked point in $\MM_{0,4}$, $-\psi$
is the first Chern class of the tangent bundle. We get contributions of the
form $\la_i-\la_j-\psi$. As in Section \ref{sec:0112}, the $S_2$-action
causes these contributions to be doubled on the overlying variety.


\begin{table}
\renewcommand{\baselinestretch}{1}
\small\normalsize
\begin{tabular*}{5.9375in}{|c|c|c|c|c|c|c|c|c|@{\extracolsep{\fill}}c|}
\hline
& $Z_1$ & $Z_2$ & $Z_3$ & $Z_4$ & $Z_5$ & $Z_6$ & $Z_7$
& $Z_8$  & $Z_9$\\ \hline 

$H_1$ & $\la_0$ & $\la_1$ &  $\la_0$ &  $\la_1$ &  $\la_1$ &  $\la_0$ &  $\la_1$
& $\la_0$ & $\la_0$  \\ \hline

$H_2$ &  $\la_1$ & $\la_0$ & $\la_0$ &  $\la_1$ & $\la_1$ & $\la_0$ & $\la_0$ 
 & $\la_1$ & $\la_1$ \\ \hline

$D_0$ & 0 & 0 & $\frac{\la_0-\la_1}{2}$ & $\frac{\la_1-\la_0}{2}$ 
& $\la_1-\la_0$ & $\la_0-\la_1$ & 0 & 0  & 0 \\ \hline

$D_1$ & 0 & 0 & 0 & 0 & $2(\la_0-\la_1)$ & $2(\la_1-\la_0)$ & 
$\la_0-\la_1$ & $\la_1-\la_0$  & $\la_0-\la_1$ \\ \hline

$D_2$ & 0 & 0 & 0 & 0 & 0 & 0 & $\la_0-\la_1$ & $\la_1-\la_0$ & $\la_0-\la_1$ 
\\ \hline

\end{tabular*}

\begin{tabular}{|c|c|c|c|c|c|c|}
\hline
& $Z_{10}$ & $Z_{11}$ & $Z_{12}$ & 
$Z_{13}$ & $Z_{14}$
 \\ \hline

$H_1$ &  $\la_1$ &  $\la_0$ &  $\la_1$ &  $\la_1$ &  $\la_0$
  \\ \hline

$H_2$  & $\la_0$ & $\la_0$ & $\la_1$ & $\la_1$ & $\la_0$
  \\ \hline

$D_0$ & 0 & $\psi$ & $\psi$ & 0 & 0   \\ \hline

$D_1$ & $\la_1-\la_0$ &  $2\la_0-2\la_1-2\psi$ 
& $2\la_1-2\la_0-2\psi$ & 0 & 0  \\ \hline

$D_2$  & $\la_1-\la_0$ &  $2\psi$ &  $2\psi$ & $2(\la_0-\la_1)$ & 
$2(\la_1-\la_0)$   \\ \hline

\end{tabular}

\caption{Restrictions of divisor classes in $A_T^*(\mm_{0,2}(\PP^1,2))$
to fixed components
\label{rest0212}}
\end{table}


We are now ready to compute the integrals.  
Details are shown for two of the integrals below.
Computations of the remaining integrals are relegated to Appendix 
\ref{sec:ints}.
Note first that
\[(\la_0-\la_1-2\psi)^{-1}=(\la_0-\la_1)^{-1}(1-2\psi/(\la_0-\la_1))^{-1}
=\frac{1+2\psi/(\la_0-\la_1)}{(\la_0-\la_1)}\text{,}\]
and similarly with $\la_0$ and $\la_1$ switched. Also note the factors of
two appearing in the denominators of the integrands for $Z_{11}$ and $Z_{12}$,
which are due to the $S_2$ automorphism group of these components. 
Also, in later computations, a factor of 2 is inserted into
the denominators of integrands for 
$Z_3$ and $Z_4$, with the $S_2$ automorphisms here occurring because of the
involution switching the sheets of the double covers.
For $D_1^2H_1H_2$, we get
\begin{eqnarray*}
\int_{\mm_T}D_1^2H_1H_2
& = & \int_{(Z_{5})_T}\frac{4\la_1^2(\la_1-\la_0)^2}{2(\la_0-\la_1)^4}+
\int_{(Z_{6})_T}\frac{4\la_0^2(\la_0-\la_1)^2}{2(\la_0-\la_1)^4}\\
&   & +\int_{(Z_{7})_T}\frac{\la_0\la_1(\la_0-\la_1)^2}{-(\la_0-\la_1)^4}
 +\int_{(Z_{8})_T}\frac{\la_0\la_1(\la_1-\la_0)^2}{-(\la_0-\la_1)^4}\\
&   & + \int_{(Z_{9})_T}\frac{\la_0\la_1(\la_0-\la_1)^2}{-(\la_0-\la_1)^4}+
\int_{(Z_{10})_T}\frac{\la_0\la_1(\la_1-\la_0)^2}{-(\la_0-\la_1)^4} \\
&   & +\int_{(Z_{11})_T}\frac{\la_0^2(2\la_0-2\la_1-2\psi)^2}{2(\la_0-\la_1)^2(\la_0-\la_1-2\psi)}\\
&   & +\int_{(Z_{12})_T}\frac{\la_1^2(2\la_1-2\la_0-2\psi)^2}{2(\la_1-\la_0)^2(\la_1-\la_0-2\psi)} \\
& = & \frac{2\la_1^2}{(\la_0-\la_1)^2}+\frac{2\la_0^2}{(\la_0-\la_1)^2}
-\frac{\la_0\la_1}{(\la_0-\la_1)^2}-\frac{\la_0\la_1}{(\la_0-\la_1)^2} \\
&  & -\frac{\la_0\la_1}{(\la_0-\la_1)^2}
-\frac{\la_0\la_1}{(\la_0-\la_1)^2}+\int_{(Z_{11})_T}
\frac{2\la_0^2(\la_0-\la_1-2\psi)}{(\la_0-\la_1)(\la_0-\la_1-2\psi)} \\
&  & 
+\int_{(Z_{12})_T}\frac{2\la_1^2(\la_1-\la_0-2\psi)}{(\la_1-\la_0)(\la_1-\la_0-2\psi)} \\
& = & \frac{2\la_0^2+2\la_1^2-4\la_0\la_1}{(\la_0-\la_1)^2}+0+0 \\
& = & 2
\text{.}\end{eqnarray*}

The integrals over $Z_{11}$ and $Z_{12}$ vanish because the integrands are
not of top codimension.
For $D_1D_2H_1H_2$, we get

\begin{eqnarray*}
\int_{\mm_T}D_1D_2H_1H_2
& = & \int_{(Z_{7})_T}\frac{\la_0\la_1(\la_0-\la_1)^2}{-(\la_0-\la_1)^4}
+  \int_{(Z_{8})_T}\frac{\la_0\la_1(\la_1-\la_0)^2}{-(\la_0-\la_1)^4} \\
&   & +\int_{(Z_{9})_T}\frac{\la_0\la_1(\la_0-\la_1)^2}{-(\la_0-\la_1)^4} 
+\int_{(Z_{10})_T}\frac{\la_0\la_1(\la_1-\la_0)^2}{-(\la_0-\la_1)^4} \\
&   & +\int_{(Z_{11})_T}\frac{\la_0^22\psi(2\la_0-2\la_1-2\psi)} 
{2(\la_0-\la_1)^2(\la_0-\la_1-2\psi)} \\
&   &
+\int_{(Z_{12})_T}\frac{\la_1^22\psi(2\la_1-2\la_0-2\psi)}
{2(\la_1-\la_0)^2(\la_1-\la_0-2\psi)} \\
& = & -4\frac{\la_0\la_1}{(\la_0-\la_1)^2}+
\int_{(Z_{11})_T}\frac{2\la_0^2\psi(2\la_0-2\la_1)(1+2\psi/(\la_0-\la_1))}
{2(\la_0-\la_1)^3} \\
&  & +\int_{(Z_{12})_T}\frac{2\la_1^2\psi(2\la_1-2\la_0)(1+2\psi/(\la_1-\la_0))}
{2(\la_1-\la_0)^3}\\
& = & \frac{-4\la_0\la_1}{(\la_0-\la_1)^2}+\int_{(Z_{11})_T}\frac{2\la_0^2\psi}
{(\la_0-\la_1)^2}+\int_{(Z_{12})_T}\frac{2\la_1^2\psi}{(\la_1-\la_0)^2} \\
& = & \frac{2\la_0^2-4\la_0\la_1+2\la_1^2}{(\la_0-\la_1)^2} \\
& = & 2
\text{.}\end{eqnarray*}

Computations of the other integrals (shown in Appendix \ref{sec:ints})
are quite similar and give the results 
shown in Table
\ref{deg4ints}.
Any integral of a degree four monomial not listed there is automatically
zero for one of the reasons given at the beginning of the section.

\renewcommand{\baselinestretch}{1}
\small\normalsize

\begin{table}
\begin{center}
\begin{tabular}{|p{2in}p{2.8in}|}
\hline
\rule[-3mm]{0mm}{8mm}$\int_{\mm}D_2^4=12$ & $\int_{\mm}D_2^3H_1=-4$ \\ 

\rule[-3mm]{0mm}{8mm}$\int_{\mm}D_2^3D_1=-4$ & $\int_{\mm}D_2^3H_2=-4$ \\
\rule[-3mm]{0mm}{8mm}
$\int_{\mm}D_2^3D_0=0$ & $\int_{\mm}D_2^2D_1H_1=0$\\
\rule[-3mm]{0mm}{8mm}
$\int_{\mm}D_2^2D_1^2=-4$ & $\int_{\mm}D_2^2D_1H_2=0$ \\
\rule[-3mm]{0mm}{8mm}
$\int_{\mm}D_2^2D_1D_0=0$ & $\int_{\mm}D_2^2D_0H_1=\int_{\mm}D_2^2D_0H_2=0$ \\
\rule[-3mm]{0mm}{8mm}
$\int_{\mm}D_2^2D_0^2=0$ & $\int_{\mm}D_2D_1^2H_1=4$\\
\rule[-3mm]{0mm}{8mm}
$\int_{\mm}D_2D_1^3=12$  & $\int_{\mm}D_2D_1^2H_2=4$\\
\rule[-3mm]{0mm}{8mm}
$\int_{\mm}D_2D_1^2D_0=0$ & $\int_{\mm}D_2D_1D_0H_1=\int_{\mm}D_2D_1D_0H_2=0$ 
\\
\rule[-3mm]{0mm}{8mm}
$\int_{\mm}D_2D_1D_0^2=0$ &  $\int_{\mm}D_2D_0^2H_1=\int_{\mm}D_2D_0^2H_2=0$\\
\rule[-3mm]{0mm}{8mm}
$\int_{\mm}D_2D_0^3=0$ &  $\int_{\mm}D_1^3H_1=-8$\\
\rule[-3mm]{0mm}{8mm}
$\int_{\mm}D_1^4=-20$ & $\int_{\mm}D_1^3H_2=-8$ \\
\rule[-3mm]{0mm}{8mm}
$\int_{\mm}D_1^3D_0=0$  & $\int_{\mm}D_1^2D_0H_1=\int_{\mm}D_1^2D_0H_2=4$\\
\rule[-3mm]{0mm}{8mm}
$\int_{\mm}D_1^2D_0^2=4$ & $\int_{\mm}D_1D_0^2H_1=\int_{\mm}D_1D_0^2H_2=-1$ \\
\rule[-3mm]{0mm}{8mm}
$\int_{\mm}D_1D_0^3=-2$ & $\int_{\mm}D_0^3H_1=\int_{\mm}D_0^3H_2=\frac{1}{4}$
\\
\rule[-3mm]{0mm}{8mm}
$\int_{\mm}D_0^4=\frac{3}{4}$ & $\int_{\mm}D_2D_1H_1H_2=2$ \\
\rule[-3mm]{0mm}{8mm}
$\int_{\mm}D_2^2H_1H_2=2$ & $\int_{\mm}D_1^2H_1H_2=2$ \\ \hline
\end{tabular}
\end{center}
\caption{Integrals of degree four classes on $\mm_{0,2}(\PP^1,2)$
\label{deg4ints}}
\end{table}

Expressions for the $\psi$-classes in $\mm_{0,2}(\PP^1,2)$ were given
in Section \ref{sec:pb}, so it suffices to consider the
basic divisor classes
$H_1$, $H_2$, $D_0$, $D_1$, and $D_2$.
Since the first Betti number is four, there must be a relation among these
five divisor classes.  Suppose we have a relation
\[aH_1+bH_2+cD_0+dD_1+eD_2=0\text{.}\]
We can place restrictions on the coefficients by multiplying the above 
equation by degree three monomials and then integrating. 
For example, multiplying by $D_2^2H_1$ gives
\[aD_2^2H_1^2+bD_2^2H_2H_1+cD_2^2D_0H_1+dD_2^2D_1H_1+eD_2^3H_1=0\text{.}\]
Now integration gives the equation
$2b-4e=0$, using the integral values from Table \ref{deg4ints}. 
Continuing with some other choices of monomial, we 
get the system of 
restrictions given in Table \ref{rest}.


\renewcommand{\baselinestretch}{1}
\small\normalsize

\begin{table}
\begin{center}
\begin{tabular}{|c|c|}
\hline
Monomial & Resulting relation on coefficients \\ \hline
 $D_2^2H_1$ & $ 2b-4e=0$ \\ 
 $D_2^2H_2$ & $ 2a-4e=0$ \\ 
 $D_1D_0H_1$ & $ -c+4d=0$ \\ 
$D_1H_1H_2$ & $ 2d+2e=0$ \\ \hline
\end{tabular}
\caption{Restrictions placed on coefficients of a linear relation by 
integration \label{rest}}
\end{center}
\end{table}



Together these restrictions show that up to a constant multiple
the only possible linear
relation among these five divisor classes is
\[2H_1+2H_2-4D_0-D_1+D_2=0\text{.}\]
Thus this must indeed be a relation, and the 
remaining four classes must be independent.  Hence, we have additionally
found that 
the classes $H_1$, $H_2$,
$D_0$, and $D_1$ generate the degree one piece of the graded ring
$A^*(\mm_{0,2}(\PP^1,2))$. The same method will be used in Section
\ref{sec:prez} to show monomials in the basic 
classes generate in degrees two and three also.
We can write the relation above as 
$D_2=4D_0+D_1-2H_1-2H_2$. Notice that this linear relation 
can also be written $D_2-\psi_1-\psi_2=0$ using the expressions for the 
$\psi$-classes from Section \ref{sec:pb}.

An interesting consequence of this relation, together with the relations
$D_0\psi_i$, is the relation $D_0D_2=0$. We can find this relation directly
by arguing that the divisors $D_0$ and $D_2$ are disjoint.
We have seen that
the domain curve of every stable map lying in
$D_0$ must have both marked points on the same degree zero component, and that
the
domain curve of every
stable map in $D_2$ must have the marked points on distinct
components.
These mutually exclusive properties of stable maps in $D_0$ and $D_2$
validate the claim of their disjointness.  Indeed, one may need to use such
a direct argument for the relation $D_0D_2$ when $r>1$, since it is not clear
whether the relation $D_2-\psi_1-\psi_2$ holds there.

\chapter{The Presentation}
\label{sec:prez}

\renewcommand{\baselinestretch}{1}
\section{
Completeness of the higher degree parts of the presentation}
\label{sec:comp}

We have already seen in Section \ref{sec:lla} that there can be no extra 
divisor classes independent of the ones given in Section \ref{sec:gen}, since
the first Betti number is four, and four of these classes are independent.
The same method used there also demonstrates that monomials in the
standard divisor classes generate in higher degrees.
As in Section \ref{sec:lla}, we need not consider the $\psi$-classes since
they can be expressed in terms of the other generators. 
We need not consider monomials involving $D_2$ for the same reason:
As found in Section \ref{sec:lla}, $D_2$ can be expressed in terms of
other generators. Furthermore,
relations involving the $\psi$-classes can be reformulated to reduce
the number of spanning monomials in the other divisor classes.

In degree 2, there are ten monomials in the remaining 
classes. However, the $H_i^2$
will not play a role since they vanish. Furthermore,
it is easy to check that $\psi_1-\psi_2=H_1-H_2$, so that
$D_0H_1-D_0H_2=D_0\psi_1-D_0\psi_2=0$. Hence we can also discount the
monomial $D_0H_2$, leaving seven monomials spanning the degree two part
of the Chow ring. Suppose we have a 
relation of the form 
\[aD_1^2+bD_1D_0+cD_1H_1+dD_1H_2+eD_0^2+fD_0H_1+gH_1H_2=0\text{.}\]
Then, multiplying this expression by each of these seven degree two monomials
and integrating the results, we obtain a system of seven linear equations
in seven variables.  Solving this system, we see that the only possible 
relation of this form (up to a constant multiple) is 
\[D_1D_0+4D_0^2-4D_0H_1=0\text{.}\]
Using the expression for $D_2$ derived in Section \ref{sec:lla}, this can be
rewritten as $D_0D_2$, a relation we have already discovered.
Thus six of the above degree two classes are independent.  Since the second
Betti number is six, these six classes generate the degree two part of 
$A^*(\mm_{0,2}(\PP^1,2))$, and so there are no other generators in degree two.

In degree three, any monomial with a factor of $D_1D_0$ can be expressed
in terms of other monomials.
Taking this into account together with the other
relations of lower degree, we have
six remaining degree three monomials in the divisor classes.  A generic 
relation among these has the form
\[aD_1^3+bD_1^2H_1+cD_1^2H_2+dD_1H_1H_2+eD_0^3+fD_0^2H_1=0\text{.}\]
Multiplying this expression by the four independent divisor classes and
integrating, we get a system of four linear equations in six variables.
The set of solutions is two-dimensional: any relation must have the form
\[aD_1^3+bD_1^2H_1+bD_1^2H_2+(-6a-4b)D_1H_1H_2+(32a-8b)D_0^3+(-96a-8b)D_0^2H_1
=0\text{.}\]
One can compute that the relations
\begin{equation}
\label{eqn:c1}
(D_1+D_2)^3=2^3(D_1^3-3D_1^2H_1-3D_1^2H_2+6D_1H_1H_2+56D_0^3-72D_0^2H_1)
\end{equation}
and
\begin{equation}
\label{eqn:c2} 
D_1\psi_1\psi_2=\frac{1}{4}D_1^3-D_1^2H_1-D_1^2H_2+
\frac{5}{2}D_1H_1H_2+16D_0^3-16D_0^2H_1
\end{equation}
found above satisfy these conditions. Moreover, they
are clearly independent.
So they span the space of relations.
Thus four of these six degree three classes must be independent.  Since the
third Betti number is four, they generate the degree three part, and there
are no additional generators.

The degree four part is one-dimensional, so since $D_1^2H_1H_2$ is nonzero,
it generates the degree four part.

\section{Two presentations for $A^*(\mm_{0,2}(\PP^1,2))$}

We have established the following result.

\begin{thm}\label{thm:prez}
With notation as established in Section \ref{sec:gen}, we have an isomorphism

\begin{equation}\label{geomprez}
A^*(\mm_{0,2}(\PP^1,2))\iso\frac{\QQ[D_0,D_1,D_2,H_1,H_2,\psi_1,\psi_2]}
{\left({H_1^2, H_2^2,D_0\psi_1,D_0\psi_2,D_2-\psi_1-\psi_2,
\psi_1-\frac{1}{4}D_1-\frac{1}{4}D_2-D_0+H_1, \atop
\psi_2-\frac{1}{4}D_1-\frac{1}{4}D_2-D_0+H_2, (D_1+D_2)^3,
D_1\psi_1\psi_2}\right)}\text{.}
\end{equation}
of graded rings.
\end{thm}

There are many different presentations for a given ring, and some are more
practical and satisfying than others. Here is another presentation for
$A^*(\mm_{0,2}(\PP^1,2))$.

\begin{prop}\label{altprez}
We also have an isomorphism
\[A^*(\mm_{0,2}(\PP^1,2))\iso\frac{\QQ[D_0,D_1,H_1,H_2]}
{\left({H_1^2, H_2^2, D_0H_1-D_0H_2, D_1D_0+4D_0^2-4D_0H_1, \atop
C_1(D_0,D_1,H_1,H_2), C_2(D_0,D_1,H_1,H_2)}\right)}\]
of graded rings,
where 
\[C_1(D_0,D_1,H_1,H_2)=D_1^3-3D_1^2H_1-3D_1^2H_2+6D_1H_1H_2+56D_0^3-72D_0^2H_1
\]
and
\[C_2(D_0,D_1,H_1,H_2)=D_1^2H_1+D_1^2H_2-4D_1H_1H_2-8D_0^3-8D_0^2H_1\]
are cubic relations. 
\end{prop}

This presentation is more
efficient than (\ref{geomprez}) in the sense that it has fewer
generators and relations. In fact, it has the minimum possible number of each.
However, 
it is not very geometric; one would be hard-pressed to give a geometric
explanation for some of the relations. A goal of efficiency also leads
to some complicated relations, and this type of presentation will be 
difficult to generalize.  Its derivation relies even more heavily on the
kind of brute force linear algebra techniques described and used in Sections
\ref{sec:lla} and \ref{sec:comp}. Including the $\psi$-classes as generators
leads to the geometric presentation (\ref{geomprez}), 
which is more beautiful and will
also be more useful.

\renewcommand{\baselinestretch}{1}
\section{
Directions for generalization}
\label{sec:direct}
\renewcommand{\baselinestretch}{2}

The most natural direction to extend Theorem \ref{thm:prez} would involve
giving 
presentations for all the rings $A^*(\mm_{0,2}(\PP^r,2))$, or at least for
those with some other small values of $r$. We already know the Poincar\'{e}
polynomials of these spaces from Chapter \ref{sec:ser}; indeed, we saw there
that the degeneration strata are the same for all $r$. So the computations 
involved and the resulting presentations
should bear many similarities to what we have seen
in the case of $A^*(\mm_{0,2}(\PP^1,2))$.

The biggest obstacle to such extension is the dependence on the localization
and linear algebra method of Section \ref{sec:lla}. 
This dependence must be removed in order to
obtain a general result. Localization and linear algebra can still be used
in obtaining presentations one dimension at a time, but even in the next
simplest case $r=2$, the magnitude of computations increases substantially.
Since $\dim(\mm_{0,2}(\PP^2,2))=7$, the method specifies finding integrals
of all degree seven monomials in the generating classes $H_1$, $H_2$, $D_0$,
$D_1$, and $D_2$. There are ${11 \choose 4}=330$ such monomials, though
as before many are zero or equivalent to other integrals
by inspection. Any monomial with a factor of $H_i^3$
will be zero since $H^3=0$ in $A^*(\PP^2)$. There are $2{8 \choose 4}-5=135$
monomials like this. The relation $D_0H_1=D_0H_2$ 
still holds by a direct argument.
(Roughly, both marked points must have the same image when they are on the
same collapsed component.) Thus we need not worry about any
monomial with a factor of $D_0H_2$. There are ${9 \choose 4}-{7 \choose 3}
-{6 \choose 2}+1=77$ of these monomials not already counted above. Finally,
$D_0$ and $D_2$ are still disjoint, so the 
${8 \choose 3}-{5 \choose 3}=46$ monomials (not already counted
above) with a factor of
$D_0D_2$ are zero as well. Taking these relations into account, 72 integrals
remain to be calculated. While the  effort necessary to perform these
calculations is not too unreasonable, it is clear that the workload will
quickly explode as $r$ increases. It could be reduced somewhat by calculating
only the integrals actually needed to find relations by linear algebra.
However, the benefit of such picking and choosing would be fleeting
as $r$ grows. Including the $\psi$-classes and making full use of 
the existing algorithms for computing
gravitational correlators may afford more significant efficiency. At any rate,
this approach cannot extend beyond some computational proofs for low $r$.

The steps needed to free our method from dependence on localization
and linear algebra were described at the beginning of Section \ref{sec:lla}.
Making rigorous the heuristic calculations of \cite{Wit} may give a way to
construct the desired section of $L_1\* L_2$. At least for low $r$, the
methods of Mumford in 
\cite{Mu} are appropriate for finding an additive basis for 
$A^*(\mm_{0,2}(\PP^r,2))$. Using the excision sequence (\ref{eksiz}), 
the author has already used these methods to
construct an additive basis for $A^*(\mm_{0,1}(\PP^1,2))$, and work on such
a construction for $A^*(\mm_{0,2}(\PP^1,2))$ is in progress. It should not
be too hard get an additive basis for general $A^*(\mm_{0,2}(\PP^r,2))$ after
seeing the pattern of the first few cases, although a different method of
proof will be required.

It should also be possible to find a presentation for all Chow rings
$A^*(\mm_{0,2}(\PP^r,2))$ by building from already existing general
presentations for $A^*(\mm_{0,0}(\PP^r,2))$ (\cite{BO}) or
$A^*(\mm_{0,1}(\PP^r,2))$ (\cite{M}, \cite{MM}). The contraction morphisms provide
a key connection between $\mm_{0,2}(\PP^r,2)$ and these moduli spaces, and
would play a central role in such a method.

Once a presentation for a Chow ring $A^*(\MM_{0,2}(\PP^r,2))$ is known, the
techniques in the previous paragraph should also apply to give a presentation
for $A^*(\MM_{0,3}(\PP^r,2))$. Ultimately, presentations of Chow rings
$A^*(\MM_{0,n}(\PP^r,2))$ for small $n$ and arbitrary $r$ should be computable
using these ideas. Since Musta\c{t}\v{a} has given a presentation for
$A^*(\MM_{0,1}(\PP^r,d))$ with $r$ and $d$ arbitrary, it may be possible
to let $d$ increase beyond two as well in this vision. The final goal in this 
realm is a description of presentations for all the Chow rings $A^*(\m0)$.
(Presentations for positive $g$ are also desirable, but this seems to be 
in a somewhat different realm.) Doing so requires more general methods.
The recent work of Oprea (\cite{O},\cite{O2}) may be a substantial step in
the right direction.

\renewcommand{\baselinestretch}{1}
\chapter{Computation of gravitational correlators using the presentation}
\label{sec:app}

Let $X$ be a smooth variety, $\b\in H_2(X)$, and $\g_1,\ldots,\g_n\in H^*(X)$.
Recall the $\psi$-classes $\psi_1,\ldots,\psi_n\in A^*(\mmm)$ defined in
Section \ref{sec:psi}. The gravitational correlator 
$\bra\t_{d_1}\g_1,\ldots,\t_{d_n}\g_n\ket$ is a rational number given by
\[\bra\t_{d_1}\g_1,\ldots,\t_{d_n}\g_n\ket_{g,\b}
=\int_{[\mmm]^\text{vir}}\prod_{i=1}^n
\left(\psi_i^{d_i}\ev_i^*(\g_i)\right)\text{.}\]
If $\b=0$, we must require either $g>0$ or $n\geq 3$. 
We usually suppress $\t_0$ from the notation, and we suppress the fundamental
class of $X$ or write it as $1$. We define any gravitational
correlator including an argument $\t_{-1}$ to be zero.

Physicists are actually
interested in correlation functions of operators called {\em gravitational
descendants}. A certain model of string theory associates to each
cohomology class $\g\in H^*(X)$ a {\em local operator} $\O_{\g}$, which is
roughly a function from a coordinate patch in $X$
to the space of
operators on a Hilbert
space of states. 
Associated to these operators are other operators $\O_{i,\g}$ called
gravitational descendants.
Moreover, all these operators act as distributions in their 
coordinates.
As such, we can take the correlation functions 
$\bra\O_{d_1,\g_1},\ldots,\O_{d_n,\g_n}\ket_{g,\b}$ of any number of
gravitational
descendants by integrating. 
In string theory they describe how various particles interact.
In fact, according to \cite{CK}, ``In a very real sense,
one can argue that the $n$-point functions contain all of the physical 
predictions of the theory."
(A correlation function with $n$ variables is also called
an $n$-point function. See \cite{CK} for more on correlation functions.)  
These are related to gravitational correlators as follows:
\[\bra\t_{d_1}\g_1,\ldots,\t_{d_n}\g_n\ket_{g,\b}
=\frac{\bra\O_{d_1,\g_1},\ldots,
\O_{d_n,\g_n}\ket_{g,\b}}{\prod_{i=1}^n d_i!}\text{.}\]

Generalizing gravitational correlators, 
we can define {\em gravitational classes} in \linebreak[4]
$H^*(\mm_{g,n},\QQ)$. Assume that $2g+n\geq 3$.
Let $\pi:\mmm\ra X^n\x \mm_{g,n}$ be defined using
the evaluation morphisms and the morphism forgetting the map. Denote the
projection from $X^n\x \mm_{g,n}$ onto its $i$'th factor by $p_i$. Let 
$PD:H^*(\mm_{g,n},\QQ)\ra H_{6g-6+2n-*}(\mm_{g,n},\QQ)$ 
be the Poincar\'{e} duality isomorphism on $\mm_{g,n}$. Then given
$\b$ and $\g_1,\ldots,\g_n$ as above, the gravitational class
$I_{g,n,\b}(\t_{d_1}\g_1,\ldots,\t_{d_n}\g_n)$ is defined to be
\[PD^{-1} {p_2}_*\left(\prod_{i=1}^n \psi_i^{d_i}\cup 
p_1^*(\g_1\*\cdots\*\g_n)\cap\pi_*([\mmm]^\text{virt})\right)\text{.}\]
The gravitational correlators and gravitational classes are related by
\[\bra\t_1\g_1,\ldots,\t_{d_n}\g_n\ket_{g,\b}
=\int_{\mm_{g,n}}I_{g,n,\b}(\t_{d_1}\g_1,\ldots,\t_{d_n}\g_n)\text{.}\]

Gromov-Witten invariants, defined in Section \ref{sec:intro}, are just
gravitational correlators with $d_i=0$ for all $i\in\und{n}$, so that there
are no $\psi$-classes in the corresponding integral.
Gromov-Witten classes can also be defined as a special case of gravitational
classes in the same way.
The above information and more can be found in \cite[Chapter 10]{CK}.

In this chapter, we will use the presentation given in Theorem \ref{thm:prez}
to compute all the genus zero, degree two, two-point gravitational correlators
of $\PP^1$. Algorithms for computing gravitational correlators have already
been constructed using indirect methods. We will show that our results agree
with the numbers computed by these existing methods. This provides a check
on the validity of the presentation.

Gravitational correlators are known to satisfy certain axioms, and the 
algorithms mentioned above make use of some of these axioms in computing the 
correlators. We now list the axioms that we will use. For the most part, these
are the same axioms listed in \cite[Chapter 10]{CK}.

\noindent
{\bf Degree Axiom.} Assume that the $\g_i$ are homogeneous.
A gravitational correlator 
$\bra\t_{d_1}\g_1,\ldots,\t_{d_n}\g_n\ket_{g,\b}$ 
can be nonzero only if the cohomological degrees of
the classes being integrated add up to twice 
the virtual dimension of the moduli
space, {\em i.e.},
\[\sum_{i=1}^n (\deg(\g_i)+2d_i)=2(1-g)(\dim_{\CC}(X)-3)
-2\int_\b K_X+2n\text{.}\]
When $X=\PP^r$, which is the case of interest, all the $\g_i$ have even
degrees. So we can divide the above equation by two and use algebraic
degrees for the $\g_i$.

\noindent
{\bf Equivariance Axiom.}
Assume the $\g_i$ are homogeneous. Then for $i\in\und{n-1}$
\begin{eqnarray*}
& & \bra\t_{d_1}\g_1,\ldots,\t_{d_{i+1}}\g_{i+1},
\t_{d_i}\g_{i},\ldots,\t_{d_n}\g_n\ket_{g,\b}\\
& = & (-1)^{\deg\g_i\cdot\deg\g_{i+1}}\bra\t_{d_1}\g_1,
\ldots,\t_{d_i}\g_{i},\t_{d_{i+1}}\g_{i+1},
\ldots,
\t_{d_n}\g_n\ket_{g,\b}\text{.}
\end{eqnarray*}
This has an obvious extension to any permutation of the entries, so that
``equivariance" refers to $S_n$-equivariance. Again, we consider only
cases where the cohomology lives
in even degrees, so that the gravitational correlators are in fact
{\em invariant} under permutation of the entries.

\noindent
{\bf Fundamental Class Axiom.}
Assume that either  $\b\neq 0$ and $n\geq 1$ or $n+2g\geq 4$. Recall that 
$1\in H^0(X,\QQ)$ denotes the fundamental class of $X$. Then
\begin{eqnarray*}
&   & \bra\t_{d_1}\g_1,\ldots,\t_{d_{n-1}}\g_{n-1},1\ket_{g,\b}\\
& = & \sum_{i=1}^{n-1}\bra\t_{d_1}\g_1,\ldots,\t_{d_{i-1}}\g_{i-1},
\t_{d_i-1}\g_{i},\t_{d_{i+1}}\g_{i+1},\ldots,\t_{d_{n-1}}\g_{n-1}
\ket_{g,\b}\text{.}
\end{eqnarray*}

\noindent
{\bf Divisor Axiom.}
Let $D\in H^2(X,\QQ)$ be a divisor class. Again assume that either  
$\b\neq 0$ and $n\geq 1$ or $n+2g\geq 4$. Then
\begin{eqnarray*}
& & \bra\t_{d_1}\g_1,\ldots,\t_{d_{n-1}}\g_{n-1},D\ket_{g,\b}\\
& = & (\int_\b D)\bra\t_{d_1}\g_1,\ldots,\t_{d_{n-1}}\g_{n-1}\ket_{g,\b}\\
 & & +\sum_{i=1}^{n-1}\bra\t_{d_1}\g_1,\ldots,\t_{d_{i-1}}\g_{i-1},
\t_{d_i-1}\g_{i}\cup D,\t_{d_{i+1}}\g_{i+1},\ldots,\t_{d_{n-1}}\g_{n-1}
\ket_{g,\b}\text{.}
\end{eqnarray*}

\noindent
{\bf Splitting Axiom.}
This axiom is easier to state in terms of gravitational classes.
Recall the gluing morphisms
\[\phi:\mm_{g_1,n_1+1}\x\mm_{g_2,n_2+1}\ra\mm_{g_1+g_2,n_1+n_2}\]
defined in Section \ref{sec:modintro}. Let $T_i$ be a homogeneous
basis for $H^*(X,\QQ)$. For each $i$ and $j$, let $g_{ij}=\int_X T_i\cup T_j$
and $(g^{ij})$ be the inverse of the matrix $(g_{ij})$. Then
\begin{eqnarray*}
& & \phi^*I_{g,n,\b}(\t_{d_1}\g_1,\ldots,\t_{d_n}\g_n)\\
& = & \sum_{\b=\b_1+\b_2}\sum_{i,j}(g^{ij}
I_{g_1,n_1+1,\b_1}(\t_{d_1}\g_1,\ldots,\t_{d_{n_1}}\g_{n_1},T_i)\\
&   &   
\text{\makebox[1in]{}}
\* I_{g_2,n_2+1,\b_2}(T_j,\t_{d_{n_1+1}}\g_{n_1+1},\ldots,
\t_{d_{n_1+n_2}}\g_{n_1+n_2}))
\text{.}
\end{eqnarray*}

\noindent
{\bf Dilaton Axiom.}
For $n\geq 1$,
\[\bra\t_1,\t_{d_1}\g_1,\ldots,\t_{d_{n}}\g_{n}\ket_{g,\b}
=(2g-2+n)\bra\t_{d_1}\g_1,\ldots,\t_{d_{n}}\g_{n}\ket_{g,\b}\text{.}\]

There are sixteen genus zero, degree two, two-point gravitational correlators
of $\PP^1$. First, we compute them using the presentation given in Section
\ref{sec:prez}. We invoke the Equivariance Axiom in
order to reduce the number of calculations to nine. This is purely a
matter of convenience; the other integrals can be calculated just as easily
and have the expected values.
The computations were carried out using the algebraic geometry
software system Macaulay 2 (\cite{GS}) by entering the presentation for the
Chow ring and instructing the program to perform multiplications in this ring.
The Macaulay 2 input and output for these calculations are shown in
Appendix \ref{sec:m2}.
The top codimension ({\em i.e.} degree) 
piece of the Chow ring is generated by any non-trivial class
of that codimension.
Thus, once we know the degree of one such class, all integrals can be computed
in terms of it. The degrees of many such classes are given in Table
\ref{deg4ints}.
For our computations in Macaulay 2, it is convenient to use
the value
\[\int_{\mm}D_1^4=-20\text{.}\]
One may wish to avoid dependence on localization by using instead a
``geometrically obvious" degree. The value
\[\int_{\mm}D_2D_1H_1H_2=2\]
is appropriate for such a purpose. Indeed, it is rather clear that there
are two points in the moduli space that satisfy the conditions these
classes impose. They correspond to the following stable maps.

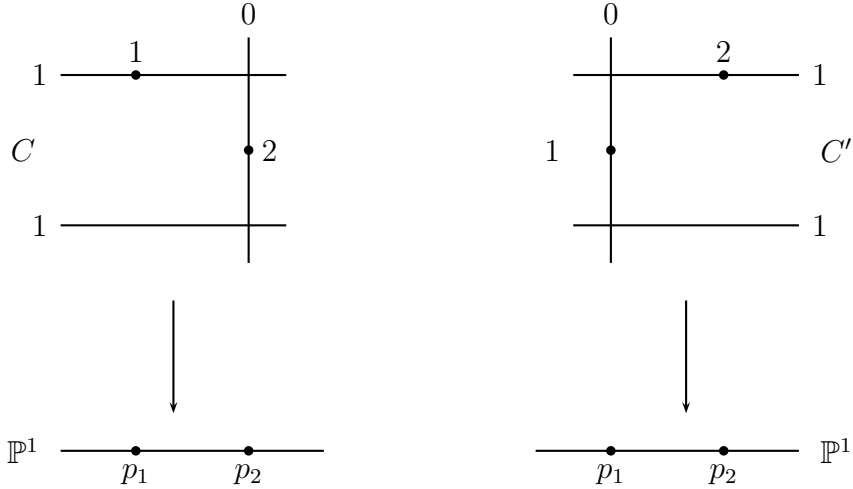
\begin{center}
\begin{pspicture}(0,-2.5)(4,4)
\rput(0,2){$C$}
\pnode(0.5,1){a}
\pnode(3.5,1){b}
\pnode(3,0.5){d}
\pnode(3,3.5){c}
\pnode(.5,3){f}
\pnode(3.5,3){e}
\dotnode(1.5,3){z}
\dotnode(3,2){y}
\ncline{a}{b}
\ncline{c}{d}
\ncline{e}{f}
\uput{5pt}[l](.5,1){1}
\uput{5pt}[l](.5,3){1}
\uput{5pt}[u](3,3.5){0}
\uput{5pt}[u](1.5,3){1}
\uput{5pt}[r](3,2){2}
\rput(0,-2){$\PP^1$}
\pnode(0.5,-2){g}
\pnode(4,-2){h}
\dotnode(1.5,-2){x}
\dotnode(3,-2){w}
\ncline{g}{h}
\uput{5pt}[d](3,-2){$p_2$}
\uput{5pt}[d](1.5,-2){$p_1$}
\pnode(2,0){i}
\pnode(2,-1.5){j}
\ncline{->}{i}{j}
\end{pspicture}
\hspace{1in}
\begin{pspicture}(0,-2.5)(4,4)
\rput(4,2){$C\pr$}
\pnode(0.5,1){a}
\pnode(3.5,1){b}
\pnode(1,0.5){d}
\pnode(1,3.5){c}
\pnode(.5,3){f}
\pnode(3.5,3){e}
\dotnode(2.5,3){z}
\dotnode(1,2){y}
\ncline{a}{b}
\ncline{c}{d}
\ncline{e}{f}
\uput{5pt}[r](3.5,1){1}
\uput{5pt}[r](3.5,3){1}
\uput{5pt}[u](1,3.5){0}
\uput{5pt}[u](2.5,3){2}
\uput{5pt}[l](0.5,2){1}
\rput(4,-2){$\PP^1$}
\pnode(0,-2){g}
\pnode(3.5,-2){h}
\dotnode(1,-2){x}
\dotnode(2.5,-2){w}
\ncline{g}{h}
\uput{5pt}[d](2.5,-2){$p_2$}
\uput{5pt}[d](1,-2){$p_1$}
\pnode(2,0){i}
\pnode(2,-1.5){j}
\ncline{->}{i}{j}
\end{pspicture}
\end{center}
Note also that these stable maps have no non-trivial automorphisms.

We show details for one example, the gravitational correlator
$\bra\t_2H,\t_1\ket_{0,2}$. By definition,
\[\bra\t_2H,\t_1\ket_{0,2}=\int_\mm \psi_1^2 H_1 \psi_2\text{.}\]
Macaulay 2 reduces the integrand to $\frac{1}{80}D^4$. Since $D^4$ has
degree $-20$, $\bra\t_2H_1,\t_1\ket=-\frac{1}{4}$. Similar computations result
in the values found in Table \ref{gravcor}. (See Appendix \ref{sec:m2}.)

\begin{table}
\begin{center}
\begin{tabular}{|cccclcr|}
\hline
$\bra\t_4,1\ket_{0,2}$ & = & $\bra 1,\t_4\ket_{0,2}$ & = & $-20\cdot\frac{3}{80}$ 
& = & $-\frac{3}{4}$\\
$\bra\t_3H,1\ket_{0,2}$ & = & $\bra 1,\t_3H\ket_{0,2}$ & = & $-20\cdot-\frac{1}{80}$ 
& = & $\frac{1}{4}$\\
$\bra\t_3,H\ket_{0,2}$ & = & $\bra H,\t_3\ket_{0,2}$ & = & $-20\cdot\frac{1}{16}$ 
& = & $-\frac{5}{4}$ \\
$\bra\t_3,\t_1\ket_{0,2}$ & = & $\bra\t_1,\t_3\ket_{0,2}$ & = & $-20\cdot-\frac{3}{80}$ 
& = & $\frac{3}{4}$\\
$\bra\t_2H,\t_1\ket_{0,2}$ & = & $\bra\t_1,\t_2H\ket_{0,2}$ & = & $-20\cdot\frac{1}{80}$ 
& = & $-\frac{1}{4}$\\
$\bra\t_2H,H\ket_{0,2}$ & = & $\bra H,\t_2H\ket_{0,2}$ & = & $-20\cdot-\frac{1}{40}$ 
& = & $\frac{1}{2}$\\
& & $\bra\t_2,\t_2\ket_{0,2}$ & = & $-20\cdot-\frac{1}{16}$ 
& = & $\frac{5}{4}$\\
$\bra\t_2,\t_1H\ket_{0,2}$ & = & $\bra\t_1H,\t_2\ket_{0,2}$ & = & $-20\cdot\frac{3}{80}$ 
& = & $-\frac{3}{4}$\\ 
& & $\bra\t_1H,\t_1H\ket_{0,2}$ & = & $-20\cdot-\frac{1}{40}$ & = & $\frac{1}{2}$\\
\hline
\end{tabular}
\end{center}
\caption{Gravitational correlators via the presentation for 
$A^*(\mm_{0,2}(\PP^1,2))$
\label{gravcor}}
\end{table}

Now we will verify the values of these gravitational correlators
by computing them using previously established methods. We will
use the Equivariance Axiom again to reduce the number of computations.
First, the following identities are derived in \cite[Chapter 10]{CK} using
the axioms and a method attributed to R. Pandharipande:
\[\bra\t_{2d-1}H,1\ket_{0,d}=\frac{1}{(d!)^2}\]
and
\[\bra\t_{2d},1\ket_{0,d}=\frac{-2}{(d!)^2}\left(1+\frac{1}{2}+\cdots+
\frac{1}{d}\right)\text{.}\]
For $d=2$ these give $\bra\t_3H,1\ket_{0,2}=\frac{1}{4}$ and 
$\bra\t_4,1\ket_{0,2}=-\frac{3}{4}$. From these, the Fundamental
Class Axiom gives $\bra\t_2H\ket_{0,2}=\frac{1}{4}$ and
$\bra\t_3\ket_{0,2}=-\frac{3}{4}$ as well. By the Divisor Axiom,
\[\bra\t_3,H\ket_{0,2}=2\bra\t_3\ket_{0,2}+\bra\t_2H\ket_{0,2}
=-\frac{5}{4}\text{.}\]
Similarly,
\[\bra\t_2H,H\ket_{0,2}=2\bra\t_2H\ket_{0,2}+\bra\t_1H^2\ket_{0,2}
=\frac{1}{2}\text{,}\]
where the second term in the sum is zero because $H^2=0$. By the
Dilaton Axiom,
\[\bra\t_1,\t_3\ket_{0,2}=-\bra\t_3\ket_{0,2}
=\frac{3}{4}\text{.}\]
Similarly,
\[\bra\t_1,\t_2H\ket_{0,2}=-\bra\t_2H\ket_{0,2}
=-\frac{1}{4}\text{.}\]

For the remaining calculations, we use results of Kontsevich and Manin
in \cite{KM2}. Let $X$ be a smooth projective manifold (variety), and let
$\bra\D_a\ket$ and $\bra\D^a\ket$ be Poincar\'{e} dual bases of $H^*(X)$.

\begin{prop}[Kontsevich-Manin]\label{wellknownid}
For $g=0$, $n=3$, and $d_1\geq 1$, we have
\[\bra\t_{d_1}\g_1,\t_{d_2}\g_2,\t_{d_3}\g_3\ket_{0,\b}
=\sum_{\b_1+\b_2=\b,a} \bra\t_{d_1-1}\g_1,\D_a\ket_{0,\b_1}
\bra\D^a,\t_{d_2}\g_2,\t_{d_3}\g_3\ket_{0,\b_2}\text{.}\]
\end{prop}

Choose some divisor $\g_0$
such that $(\g_0,\b)=\int_\b \g_0\neq 0$. 
Using the Divisor Axiom, they derive the
following identity for  genus zero, two-point correlators with $d_1>0$:

\begin{eqnarray*}
\bra\t_{d_1}\g_1,\t_{d_2}\g_2\ket_{0,\b}
& = & \frac{1}{(\g_0,\b)}(\bra\g_0,\t_{d_1}\g_1,\t_{d_2}\g_2\ket_{0,\b}\\
&   & -\bra\t_{d_1-1}(\g_0\cup\g_1),\t_{d_2}\g_2\ket_{0,\b}
-\bra\t_{d_1}\g_1,\t_{d_2-1}(\g_0\cup\g_2)\ket_{0,\b})\text{.}
\end{eqnarray*}

A gravitational correlator is called {\em primary} if all the $d_i$ are
zero, {\em i.e.}, if it is a Gromov-Witten invariant. The identity above 
can be used repeatedly to reduce to an expression in terms of primary
correlators, whose calculations are relatively straightforward. In the last
two terms, the sums of the $\t$ subscripts are already smaller. To reduce
the first term, we apply Proposition \ref{wellknownid} together
with the Equivariance Axiom to get
\[\bra\g_0,\t_{d_1}\g_1,\t_{d_2}\g_2\ket_{0,\b}
=\sum_{\b_1+\b_2=\b,a}\bra\t_{d_1-1}\g_1,\D_a\ket_{0,\b_1}
\bra\D^a,\g_0,\t_{d_2}\g_2\ket_{0,\b_2}\text{.}\]

We take the dual bases $\bra 1,H\ket$ and $\bra H,1\ket$ for $H^*(\PP^1)$, and
$\g_0=H$ is the obvious choice. Applying the above procedure gives
\begin{eqnarray*}
\bra\t_1H,\t_1H\ket_{0,2} & = & \frac{1}{2}(\bra H,\t_1H,\t_1H\ket_{0,2}
-\bra H^2,\t_1H\ket_{0,2}-\bra \t_1H,H^2\ket_{0,2})\\
 & = & \frac{1}{2}\sum_{d_1+d_2=2,a}\bra H,\D_a\ket_{0,d_1}\bra\D^a,H,\t_1 H
\ket_{0,d_2}-0-0\\
 & = & \frac{1}{2}(\bra H,H\ket_{0,1}\bra 1,H,\t_1H\ket_{0,1})\\
 & = & \frac{1}{2}\cdot 1\cdot\bra H,H\ket_{0,1}\\
 & = & \frac{1}{2}\text{,}
\end{eqnarray*}
where we used the Fundamental Class Axiom in going from the third line to the
fourth.  Notice also that the
Degree Axiom substantially limits the number of terms in the sum that can be 
nonzero. Finally, $\bra H,H\ket_{0,1}=1$ is the degree of the class of a
point in $\PP^1\x\PP^1$.
We will use the same procedure to calculate 
$\bra\t_1H,\t_2\ket_{0,2}$ and $\bra\t_2,\t_2\ket_{0,2}$. We get
\begin{eqnarray*}
\bra\t_1H,\t_2\ket_{0,2} & = & \frac{1}{2}(\bra H,\t_1H,\t_2\ket_{0,2}
-\bra H^2,\t_2\ket_{0,2}-\bra \t_1H,\t_1H\ket_{0,2})\\
 & = & \frac{1}{2}\sum_{d_1+d_2=2,a}\bra H,\D_a\ket_{0,d_1}\bra\D^a,H,\t_2
\ket_{0,d_2}-0-\frac{1}{4}\\
 & = & \frac{1}{2}\bra H,H\ket_{0,1}\bra 1,H,\t_2\ket_{0,1}-\frac{1}{4}\\
 & = & \frac{1}{2}\cdot 1\cdot\bra H,\t_1\ket_{0,1}-\frac{1}{4}\\
 & = & -\frac{3}{4}
\end{eqnarray*}
and
\begin{eqnarray*}
\bra\t_2,\t_2\ket_{0,2} & = & \frac{1}{2}(\bra H,\t_2,\t_2\ket_{0,2}
-\bra \t_1H,\t_2\ket_{0,2}-\bra \t_2,\t_1H\ket_{0,2})\\
 & = & \frac{1}{2}\sum_{d_1+d_2=2,a}\bra \t_1,\D_a\ket_{0,d_1}\bra\D^a,H,\t_2
\ket_{0,d_2}+\frac{3}{8}+\frac{3}{8}\\
 & = & \frac{1}{2}\bra \t_1,H\ket_{0,1}\bra 1,H,\t_2\ket_{0,1}+\frac{3}{4}\\
 & = & \frac{1}{2}\cdot(-1)\cdot(-1)+\frac{3}{4}\\
 & = & \frac{5}{4}\text{,}
\end{eqnarray*}
where the last steps are obtained using the Fundamental Class Axiom, the
Dilaton Axiom, and $\bra H\ket_{0,1}=1$ in each case. (The gravitational
correlator $\bra H\ket_{0,1}$ is just the degree of the class of a point
in $\PP^1$.) Observe
that all of the values computed by these standard methods agree with those
in Table \ref{gravcor}.

\appendix



\renewcommand{\baselinestretch}{1}
\chapter{Integral computations in $A^*(\MM_{0,2}(\PP^1,2))$}
\label{sec:ints}

In Section \ref{sec:lla}, we describe the ``localization and linear algebra"
method for finding relations and additive bases of Chow rings of moduli spaces
of stable maps to projective spaces. We then apply this method to learn
some things about the structure of the ring $A^*(\mm_{0,2}(\PP^1,2))$.
The key step involves computing the integrals of all degree four monomials
in the boundary divisor and hyperplane pullback classes. Section \ref{sec:eq}
develops the general theory of localization required for these computations,
while Section \ref{sec:lla} gives the specific ingredients needed to apply
localization in $A^*(\mm_{0,2}(\PP^1,2))$. Recall that there are 70 such
integrals, and 32 are immediately seen to be zero by previously known 
relations. Calculations for $D_1^2H_1H_2$ and $D_1D_2H_1H_2$ are performed
in Section \ref{sec:lla}. In this Appendix, we will compute the remaining 36
integrals. We can reduce the number of computations slightly using the fact
that $D_0H_1=D_0H_2$ (\ref{d0hi}). Because of this, we automatically know
that the integral of any multiple of $D_0H_2$ will have the same value as
the integral of the corresponding multiple of $D_0H_1$.
The results are summarized in Table \ref{deg4ints}.
Let $\mm=\mm_{0,2}(\PP^1,2)$. We now mention a couple of details that are
important to keep in mind. Recall that $\psi^2=0$ on $Z_{11}$ and $Z_{12}$.
Remember also that extra factors of 2 occur in the denominators of the
integrands over $Z_3$, $Z_4$, $Z_{11}$, and $Z_{12}$ because these fixed 
components have automorphism group $S_2$. Lastly, note that 
the $\la_i$ are treated as constants by the equivariant integrals.

For $D_2^4$, we compute

\renewcommand{\baselinestretch}{1}
\small\normalsize

\begin{eqnarray*}
& & \int_{\mm_T}D_2^4\\
& = & \int_{(Z_{7})_T}\frac{(\la_0-\la_1)^4}{-(\la_0-\la_1)^4}
+\int_{(Z_{8})_T}\frac{(\la_0-\la_1)^4}{-(\la_0-\la_1)^4}\\
& & +\int_{(Z_{9})_T}\frac{(\la_0-\la_1)^4}{-(\la_0-\la_1)^4}
+\int_{(Z_{10})_T}\frac{(\la_0-\la_1)^4}{-(\la_0-\la_1)^4}\\
& & +\int_{(Z_{11})_T}\frac{16\psi^4}{2(\la_0-\la_1)^2(\la_0-\la_1-2\psi)}
+\int_{(Z_{12})_T}\frac{16\psi^4}{2(\la_1-\la_0)^2(\la_1-\la_0-2\psi)}\\
& & +\int_{(Z_{13})_T}\frac{16(\la_0-\la_1)^4}{2(\la_0-\la_1)^4}
+\int_{(Z_{14})_T}\frac{16(\la_0-\la_1)^4}{2(\la_0-\la_1)^4}  \\
& = & -1-1-1-1+0+0+8+8  \\
& = & 12
\text{.}\end{eqnarray*}

For $D_2^3D_1$, we compute

\begin{eqnarray*}
& & \int_{\mm_T}D_2^3D_1 \\
& = & \int_{(Z_{7})_T}\frac{(\la_0-\la_1)^4}{-(\la_0-\la_1)^4}
+\int_{(Z_{8})_T}\frac{(\la_0-\la_1)^4}{-(\la_0-\la_1)^4}\\
&   & +\int_{(Z_{9})_T}\frac{(\la_0-\la_1)^4}{-(\la_0-\la_1)^4}
+\int_{(Z_{10})_T}\frac{(\la_0-\la_1)^4}{-(\la_0-\la_1)^4}\\
&   & +\int_{(Z_{11})_T}\frac{8\psi^32(\la_0-\la_1-\psi)}{2(\la_0-\la_1)^2(\la_0-\la_1-2\psi)}
+\int_{(Z_{12})_T}\frac{8\psi^32(\la_1-\la_0-\psi)}{2(\la_1-\la_0)^2(\la_1-\la_0-2\psi)}\\
& = & -1-1-1-1+0+0\\
& = & -4
\text{.}\end{eqnarray*}

For $D_2^3D_0$, we compute

\begin{eqnarray*}
& & \int_{\mm_T}D_2^3D_0 \\
& = & 
\int_{(Z_{11})_T}\frac{8\psi^3\psi}{2(\la_0-\la_1)^2(\la_0-\la_1-2\psi)}
+\int_{(Z_{12})_T}\frac{8\psi^3\psi}{2(\la_1-\la_0)^2(\la_1-\la_0-2\psi)}\\
& = & 0+0\\
& = & 0
\text{.}\end{eqnarray*}

For $D_2^2D_1^2$, we compute

\begin{eqnarray*}
& & \int_{\mm_T}D_2^2D_1^2 \\
& = & \int_{(Z_{7})_T}\frac{(\la_0-\la_1)^4}{-(\la_0-\la_1)^4}
+\int_{(Z_{8})_T}\frac{(\la_0-\la_1)^4}{-(\la_0-\la_1)^4}\\
&   & +\int_{(Z_{9})_T}\frac{(\la_0-\la_1)^4}{-(\la_0-\la_1)^4}
+\int_{(Z_{10})_T}\frac{(\la_0-\la_1)^4}{-(\la_0-\la_1)^4}\\
&   & +\int_{(Z_{11})_T}\frac{4\psi^24(\la_0-\la_1-\psi)^2}{2(\la_0-\la_1)^2(\la_0-\la_1-2\psi)}
+\int_{(Z_{12})_T}\frac{4\psi^24(\la_1-\la_0-\psi)^2}{2(\la_1-\la_0)^2(\la_1-\la_0-2\psi)}\\
& = & -1-1-1-1+0+0\\
& = & -4
\text{.}\end{eqnarray*}

For $D_2^2D_1D_0$, we compute

\begin{eqnarray*}
& & \int_{\mm_T}D_2^2D_1D_0 \\
& = & \int_{(Z_{11})_T}\frac{4\psi^2}{2(\la_0-\la_1)^2(\la_0-\la_1-2\psi)}
+\int_{(Z_{12})_T}\frac{4\psi^2}{2(\la_1-\la_0)^2(\la_1-\la_0-2\psi)}\\
& = & 0+0\\
& = & 0
\text{.}\end{eqnarray*}

For $D_2^2D_0^2$, we compute

\begin{eqnarray*}
& & \int_{\mm_T}D_2^2D_0^2 \\
& = & \int_{(Z_{11})_T}\frac{4\psi^2\psi^2}{2(\la_0-\la_1)^2(\la_0-\la_1-2\psi)}
+\int_{(Z_{12})_T}\frac{4\psi^2\psi^2}{2(\la_1-\la_0)^2(\la_1-\la_0-2\psi)}\\
& = & 0+0\\
& = & 0
\text{.}\end{eqnarray*}

For $D_2D_1^3$, we compute

\begin{eqnarray*}
& & \int_{\mm_T}D_2D_1^3 \\
 & = & \int_{(Z_{7})_T}\frac{(\la_0-\la_1)^4}{-(\la_0-\la_1)^4}
+\int_{(Z_{8})_T}\frac{(\la_0-\la_1)^4}{-(\la_0-\la_1)^4}\\
&   & +\int_{(Z_{9})_T}\frac{(\la_0-\la_1)^4}{-(\la_0-\la_1)^4}
+\int_{(Z_{10})_T}\frac{(\la_0-\la_1)^4}{-(\la_0-\la_1)^4}\\
&   & +\int_{(Z_{11})_T}\frac{2\psi 8(\la_0-\la_1-\psi)^3}{2(\la_0-\la_1)^2(\la_0-\la_1-2\psi)}
+\int_{(Z_{12})_T}\frac{2\psi 8(\la_1-\la_0-\psi)^3}{2(\la_1-\la_0)^2(\la_1-\la_0-2\psi)}\\
& = & -1-1-1-1+\int_{(Z_{11})_T}\frac{8\psi(\la_0-\la_1-\psi)}
{\la_0-\la_1-2\psi}\\
&   & +\int_{(Z_{12})_T}\frac{8\psi(\la_1-\la_0-\psi)}{\la_1-\la_0-2\psi}\\
& = & -4+\int_{(Z_{11})_T}\frac{8\psi(\la_0-\la_1-\psi)(1+2\psi/(\la_0-\la_1))}
{\la_0-\la_1}\\
&   & +\int_{(Z_{12})_T}\frac{8\psi(\la_1-\la_0-\psi)(1+2\psi/(\la_1-\la_0))}
{\la_1-\la_0}\\
& = & -4+\int_{(Z_{11})_T}8\psi(1+2\psi/(\la_0-\la_1))
+\int_{(Z_{12})_T}8\psi(1+2\psi/(\la_1-\la_0))\\
& = & -4+8+8\\
& = & 12
\text{.}\end{eqnarray*}

For $D_2D_1^2D_0$, we compute

\begin{eqnarray*}
& & \int_{\mm_T}D_2D_1^2D_0 \\
& = & \int_{(Z_{11})_T}\frac{2\psi 4(\la_0-\la_1-\psi)^2\psi}{2(\la_0-\la_1)^2(\la_0-\la_1-2\psi)}
+\int_{(Z_{12})_T}\frac{2\psi 4(\la_1-\la_0-\psi)^2\psi}{2(\la_1-\la_0)^2(\la_1-\la_0-2\psi)}\\
& = & 0+0\\
& = & 0
\text{.}\end{eqnarray*}

For $D_2D_1D_0^2$, we compute

\begin{eqnarray*}
& & \int_{\mm_T}D_2D_1D_0^2 \\
& = & \int_{(Z_{11})_T}\frac{2\psi 2(\la_0-\la_1-\psi)\psi^2}{2(\la_0-\la_1)^2(\la_0-\la_1-2\psi)}
+\int_{(Z_{12})_T}\frac{2\psi 2(\la_1-\la_0-\psi)\psi^2}{2(\la_1-\la_0)^2(\la_1-\la_0-2\psi)}\\
& = & 0+0\\
& = & 0
\text{.}\end{eqnarray*}

For $D_2D_0^3$, we compute

\begin{eqnarray*}
& & \int_{\mm_T}D_2D_0^3 \\
& = & \int_{(Z_{11})_T}\frac{2\psi\psi^3}{2(\la_0-\la_1)^2(\la_0-\la_1-2\psi)}
+\int_{(Z_{12})_T}\frac{2\psi\psi^3}{2(\la_1-\la_0)^2(\la_1-\la_0-2\psi)}\\
& = & 0+0\\
& = & 0
\text{.}\end{eqnarray*}

For $D_1^4$, we compute

\begin{eqnarray*}
&   & \int_{\mm_T}D_1^4 \\
& = & \int_{(Z_{5})_T}\frac{16(\la_0-\la_1)^4}{2(\la_1-\la_0)^4}
+\int_{(Z_{6})_T}\frac{16(\la_0-\la_1)^4}{2(\la_1-\la_0)^4}\\
&   & +\int_{(Z_{7})_T}\frac{(\la_0-\la_1)^4}{-(\la_0-\la_1)^4}
+\int_{(Z_{8})_T}\frac{(\la_0-\la_1)^4}{-(\la_0-\la_1)^4}\\
&   & +\int_{(Z_{9})_T}\frac{(\la_0-\la_1)^4}{-(\la_0-\la_1)^4}
+\int_{(Z_{10})_T}\frac{(\la_0-\la_1)^4}{-(\la_0-\la_1)^4}\\
&   & +\int_{(Z_{11})_T}\frac{16(\la_0-\la_1-\psi)^4}{2(\la_0-\la_1)^2(\la_0-\la_1-2\psi)}
+\int_{(Z_{12})_T}\frac{16(\la_1-\la_0-\psi)^4}{2(\la_1-\la_0)^2(\la_1-\la_0-2\psi)}\\
& = & 8+8-1-1-1-1+\int_{(Z_{11})_T}\frac{8(\la_0-\la_1)(\la_0-\la_1-4\psi)}
{\la_0-\la_1-2\psi}\\
&   & +\int_{(Z_{12})_T}\frac{8(\la_1-\la_0)(\la_1-\la_0-4\psi)}
{\la_1-\la_0-2\psi}\\
& = & 12+\int_{(Z_{11})_T}8(\la_0-\la_1-4\psi)(1+2\psi/(\la_0-\la_1))\\
&   & +\int_{(Z_{12})_T}8(\la_1-\la_0-4\psi)(1+2\psi/(\la_1-\la_0))\\
& = & 12-16-16\\
& = & -20
\text{.}\end{eqnarray*}

For $D_1^3D_0$, we compute

\begin{eqnarray*}
& & \int_{\mm_T}D_1^3D_0 \\
& = & \int_{(Z_{5})_T}\frac{8(\la_0-\la_1)^3(\la_1-\la_0)}{2(\la_1-\la_0)^4}
+\int_{(Z_{6})_T}\frac{8(\la_1-\la_0)^3(\la_0-\la_1)}{2(\la_1-\la_0)^4}\\
&   & +\int_{(Z_{11})_T}\frac{8(\la_0-\la_1-\psi)^3\psi}{2(\la_0-\la_1)^2(\la_0-\la_1-2\psi)}
+\int_{(Z_{12})_T}\frac{8(\la_1-\la_0-\psi)^3\psi}{2(\la_1-\la_0)^2(\la_1-\la_0-2\psi)}\\
& = & -4-4+\int_{(Z_{11})_T}\frac{4(\la_0-\la_1-3\psi)\psi(1+2\psi/(\la_0-\la_1))}
{\la_0-\la_1}\\
&   & +\int_{(Z_{12})_T}\frac{4(\la_1-\la_0-\psi)\psi(1+2\psi/(\la_1-\la_0))}
{\la_1-\la_0}\\
& = & -8+4+4\\
& = & 0
\text{.}\end{eqnarray*}

For $D_1^2D_0^2$, we compute

\begin{eqnarray*}
& & \int_{\mm_T}D_1^2D_0^2 \\
& = & \int_{(Z_{5})_T}\frac{4(\la_0-\la_1)^2(\la_1-\la_0)^2}{2(\la_1-\la_0)^4}
+\int_{(Z_{6})_T}\frac{4(\la_1-\la_0)^2(\la_0-\la_1)^2}{2(\la_1-\la_0)^4}\\
&   & +\int_{(Z_{11})_T}\frac{4(\la_0-\la_1-\psi)^2\psi^2}{2(\la_0-\la_1)^2(\la_0-\la_1-2\psi)}
+\int_{(Z_{12})_T}\frac{4(\la_1-\la_0-\psi)^2\psi^2}{2(\la_1-\la_0)^2(\la_1-\la_0-2\psi)}\\
& = & 2+2+0+0\\
& = & 4
\text{.}\end{eqnarray*}

For $D_1D_0^3$, we compute

\begin{eqnarray*}
& & \int_{\mm_T}D_1D_0^3 \\
& = & \int_{(Z_{5})_T}\frac{2(\la_0-\la_1)(\la_1-\la_0)^3}{2(\la_1-\la_0)^4}
+\int_{(Z_{6})_T}\frac{2(\la_1-\la_0)(\la_0-\la_1)^3}{2(\la_1-\la_0)^4}\\
&   & +\int_{(Z_{11})_T}\frac{2(\la_0-\la_1-\psi)\psi^3}{2(\la_0-\la_1)^2(\la_0-\la_1-2\psi)}
+\int_{(Z_{12})_T}\frac{2(\la_1-\la_0-\psi)\psi^3}{2(\la_1-\la_0)^2(\la_1-\la_0-2\psi)}\\
& = & -1-1+0+0\\
& = & -2
\text{.}\end{eqnarray*}

For $D_0^4$, we compute

\begin{eqnarray*}
& & \int_{\mm_T}D_0^4 \\
& = & \int_{(Z_{3})_T}\frac{(\la_1-\la_0)^4/16}{-2(\la_1-\la_0)^4/4}
+\int_{(Z_{4})_T}\frac{(\la_1-\la_0)^4/16}{-2(\la_1-\la_0)^4/4}\\
&   & +\int_{(Z_{5})_T}\frac{(\la_1-\la_0)^4}{2(\la_1-\la_0)^4}
+\int_{(Z_{6})_T}\frac{(\la_1-\la_0)^4}{2(\la_1-\la_0)^4}\\
&   & +\int_{(Z_{11})_T}\frac{\psi^4}{2(\la_0-\la_1)^2(\la_0-\la_1-2\psi)}
+\int_{(Z_{12})_T}\frac{\psi^4}{2(\la_1-\la_0)^2(\la_1-\la_0-2\psi)}\\
& = & -\frac{1}{8}-\frac{1}{8}+\frac{1}{2}+\frac{1}{2}+0+0\\
& = & \frac{3}{4}
\text{.}\end{eqnarray*}

For $D_2^3H_1$, we compute

\begin{eqnarray*}
& & \int_{\mm_T}D_2^3H_1 \\
& = & \int_{(Z_{7})_T}\frac{(\la_0-\la_1)^3\la_1}{-(\la_0-\la_1)^4}
+\int_{(Z_{8})_T}\frac{(\la_1-\la_0)^3\la_0}{-(\la_0-\la_1)^4}\\
&   & +\int_{(Z_{9})_T}\frac{(\la_0-\la_1)^3\la_0}{-(\la_0-\la_1)^4}
+\int_{(Z_{10})_T}\frac{(\la_1-\la_0)^3\la_1}{-(\la_0-\la_1)^4}\\
&   & +\int_{(Z_{11})_T}\frac{8\psi^3\la_0}{2(\la_0-\la_1)^2(\la_0-\la_1-2\psi)}
+\int_{(Z_{12})_T}\frac{8\psi^3\la_1}{2(\la_1-\la_0)^2(\la_1-\la_0-2\psi)}\\
&   & +\int_{(Z_{13})_T}\frac{8(\la_0-\la_1)^3\la_1}{2(\la_0-\la_1)^4}
+\int_{(Z_{14})_T}\frac{8(\la_1-\la_0)^3\la_0}{2(\la_0-\la_1)^4}  \\
& = & -\frac{\la_1}{\la_0-\la_1}-\frac{\la_0}{\la_1-\la_0}
-\frac{\la_0}{\la_0-\la_1}-\frac{\la_1}{\la_1-\la_0}\\
&   & +0+0+\frac{4\la_1}{\la_0-\la_1}+\frac{4\la_0}{\la_1-\la_0}\\
& = & 1-1-4\\
& = & -4
\text{.}\end{eqnarray*}

For $D_2^3H_2$, we compute

\begin{eqnarray*}
& & \int_{\mm_T}D_2^3H_2 \\
& = & \int_{(Z_{7})_T}\frac{(\la_0-\la_1)^3\la_0}{-(\la_0-\la_1)^4}
+\int_{(Z_{8})_T}\frac{(\la_1-\la_0)^3\la_1}{-(\la_0-\la_1)^4}\\
&   & +\int_{(Z_{9})_T}\frac{(\la_0-\la_1)^3\la_1}{-(\la_0-\la_1)^4}
+\int_{(Z_{10})_T}\frac{(\la_1-\la_0)^3\la_0}{-(\la_0-\la_1)^4}\\
&   & +\int_{(Z_{11})_T}\frac{8\psi^3\la_0}{2(\la_0-\la_1)^2(\la_0-\la_1-2\psi)}
+\int_{(Z_{12})_T}\frac{8\psi^3\la_1}{(2\la_1-\la_0)^2(\la_1-\la_0-2\psi)}\\
&   & +\int_{(Z_{13})_T}\frac{8(\la_0-\la_1)^3\la_1}{2(\la_0-\la_1)^4}
+\int_{(Z_{14})_T}\frac{8(\la_1-\la_0)^3\la_0}{2(\la_0-\la_1)^4}  \\
& = & -\frac{\la_0}{\la_0-\la_1}-\frac{\la_1}{\la_1-\la_0}
-\frac{\la_1}{\la_0-\la_1}-\frac{\la_0}{\la_1-\la_0}\\
&   & +0+0+\frac{4\la_1}{\la_0-\la_1}+\frac{4\la_0}{\la_1-\la_0}\\
& = & -1+1-4\\
& = & -4
\text{.}\end{eqnarray*}

For $D_2^2D_1H_1$, we compute

\begin{eqnarray*}
& & \int_{\mm_T}D_2^2D_1H_1 \\
& = & \int_{(Z_{7})_T}\frac{(\la_0-\la_1)^2(\la_0-\la_1)\la_1}{-(\la_0-\la_1)^4}
+\int_{(Z_{8})_T}\frac{(\la_0-\la_1)^2(\la_1-\la_0)\la_0}{-(\la_0-\la_1)^4}\\
&   & +\int_{(Z_{9})_T}\frac{(\la_0-\la_1)^2(\la_0-\la_1)\la_0}{-(\la_0-\la_1)^4}
+\int_{(Z_{10})_T}\frac{(\la_0-\la_1)^2(\la_1-\la_0)\la_1}{-(\la_0-\la_1)^4}\\
&   & +\int_{(Z_{11})_T}\frac{4\psi^22(\la_0-\la_1-\psi)\la_0}{2(\la_0-\la_1)^2(\la_0-\la_1-2\psi)}
+\int_{(Z_{12})_T}\frac{4\psi^22(\la_1-\la_0-\psi)\la_1}{2(\la_1-\la_0)^2(\la_1-\la_0-2\psi)}\\
& = & -\frac{\la_1}{\la_0-\la_1}-\frac{\la_0}{\la_1-\la_0}\\
&   & -\frac{\la_0}{\la_0-\la_1}-\frac{\la_1}{\la_1-\la_0}+0+0\\
& = & 1-1\\
& = & 0
\text{.}\end{eqnarray*}

For $D_2^2D_1H_2$, we compute

\begin{eqnarray*}
& & \int_{\mm_T}D_2^2D_1H_2 \\
& = & \int_{(Z_{7})_T}\frac{(\la_0-\la_1)^2(\la_0-\la_1)\la_0}{-(\la_0-\la_1)^4}
+\int_{(Z_{8})_T}\frac{(\la_0-\la_1)^2(\la_1-\la_0)\la_1}{-(\la_0-\la_1)^4}\\
&   & +\int_{(Z_{9})_T}\frac{(\la_0-\la_1)^2(\la_0-\la_1)\la_1}{-(\la_0-\la_1)^4}
+\int_{(Z_{10})_T}\frac{(\la_0-\la_1)^2(\la_1-\la_0)\la_0}{-(\la_0-\la_1)^4}\\
&   & +\int_{(Z_{11})_T}\frac{4\psi^22(\la_0-\la_1-\psi)\la_0}{2(\la_0-\la_1)^2(\la_0-\la_1-2\psi)}
+\int_{(Z_{12})_T}\frac{4\psi^22(\la_1-\la_0-\psi)\la_1}{2(\la_1-\la_0)^2(\la_1-\la_0-2\psi)}\\
& = & -\frac{\la_0}{\la_0-\la_1}-\frac{\la_1}{\la_1-\la_0}\\
&   & -\frac{\la_1}{\la_0-\la_1}-\frac{\la_0}{\la_1-\la_0}+0+0\\
& = & -1+1\\
& = & 0
\text{.}\end{eqnarray*}

For $D_2^2D_0H_1$, we compute

\begin{eqnarray*}
& & \int_{\mm_T}D_2^2D_0H_1 \\
& = & \int_{(Z_{11})_T}\frac{4\psi^2\psi\la_0}{2(\la_0-\la_1)^2(\la_0-\la_1-2\psi)}
+\int_{(Z_{12})_T}\frac{4\psi^2\psi\la_1}{2(\la_1-\la_0)^2(\la_1-\la_0-2\psi)}\\
& = & 0+0\\
& = & 0
\text{.}\end{eqnarray*}

Thus $\int_{\mm}D_2^2D_0H_2=0$ as well.
For $D_2D_1^2H_1$, we compute

\begin{eqnarray*}
& & \int_{\mm_T}D_2D_1^2H_1 \\
& = & \int_{(Z_{7})_T}\frac{(\la_0-\la_1)(\la_0-\la_1)^2\la_1}{-(\la_0-\la_1)^4}
+\int_{(Z_{8})_T}\frac{(\la_1-\la_0)(\la_0-\la_1)^2\la_0}{-(\la_0-\la_1)^4}\\
&   & +\int_{(Z_{9})_T}\frac{(\la_0-\la_1)(\la_0-\la_1)^2\la_0}{-(\la_0-\la_1)^4}
+\int_{(Z_{10})_T}\frac{(\la_1-\la_0)(\la_0-\la_1)^2\la_1}{-(\la_0-\la_1)^4}\\
&   & +\int_{(Z_{11})_T}\frac{2\psi 4(\la_0-\la_1-\psi)^2\la_0}{2(\la_0-\la_1)^2(\la_0-\la_1-2\psi)}
+\int_{(Z_{12})_T}\frac{2\psi 4(\la_1-\la_0-\psi)^2\la_1}{2(\la_1-\la_0)^2(\la_1-\la_0-2\psi)}\\
& = & -\frac{\la_1}{\la_0-\la_1}-\frac{\la_0}{\la_1-\la_0}
 -\frac{\la_0}{\la_0-\la_1}-\frac{\la_1}{\la_1-\la_0}\\
&   & +\int_{(Z_{11})_T}\frac{4\psi(\la_0-\la_1-2\psi)\la_0
}{(\la_0-\la_1)(\la_0-\la_1-2\psi)}\\
&  & +\int_{(Z_{12})_T}\frac{4\psi (\la_1-\la_0-2\psi)\la_1
}{(\la_1-\la_0)(\la_1-\la_0-2\psi)} \\
& = & 1-1+\frac{4\la_0}{\la_0-\la_1}+\frac{4\la_1}{\la_1-\la_0}\\
& = & 4
\text{.}\end{eqnarray*}

For $D_2D_1^2H_2$, we compute

\begin{eqnarray*}
& & \int_{\mm_T}D_2D_1^2H_2 \\
& = & \int_{(Z_{7})_T}\frac{(\la_0-\la_1)(\la_0-\la_1)^2\la_0}{-(\la_0-\la_1)^4}
+\int_{(Z_{8})_T}\frac{(\la_1-\la_0)(\la_0-\la_1)^2\la_1}{-(\la_0-\la_1)^4}\\
&   & +\int_{(Z_{9})_T}\frac{(\la_0-\la_1)(\la_0-\la_1)^2\la_1}{-(\la_0-\la_1)^4}
+\int_{(Z_{10})_T}\frac{(\la_1-\la_0)(\la_0-\la_1)^2\la_0}{-(\la_0-\la_1)^4}\\
&   & +\int_{(Z_{11})_T}\frac{2\psi 4(\la_0-\la_1-\psi)^2\la_0}{2(\la_0-\la_1)^2(\la_0-\la_1-2\psi)}
+\int_{(Z_{12})_T}\frac{2\psi 4(\la_1-\la_0-\psi)^2\la_1}{2(\la_1-\la_0)^2(\la_1-\la_0-2\psi)}\\
& = & -\frac{\la_0}{\la_0-\la_1}-\frac{\la_1}{\la_1-\la_0}
 -\frac{\la_1}{\la_0-\la_1}-\frac{\la_0}{\la_1-\la_0}\\
&   & +\int_{(Z_{11})_T}\frac{4\psi(\la_0-\la_1-2\psi)\la_0
}{(\la_0-\la_1)(\la_0-\la_1-2\psi)}\\
&   & +\int_{(Z_{12})_T}\frac{4\psi (\la_1-\la_0-2\psi)\la_1
}{(\la_1-\la_0)(\la_1-\la_0-2\psi)}\\
& = & -1+1+\frac{4\la_0}{\la_0-\la_1}+\frac{4\la_1}{\la_1-\la_0}\\
& = & 4
\text{.}\end{eqnarray*}

For $D_2D_1D_0H_1$, we compute

\begin{eqnarray*}
& & \int_{\mm_T}D_2D_1D_0H_1 \\
& = & 
\int_{(Z_{11})_T}\frac{2\psi 2(\la_0-\la_1-\psi)\psi\la_0}{2(\la_0-\la_1)^2(\la_0-\la_1-2\psi)}
+\int_{(Z_{12})_T}\frac{2\psi 2(\la_1-\la_0-\psi)\psi\la_1}{2(\la_1-\la_0)^2(\la_1-\la_0-2\psi)}\\
& = & 0+0\\
& = & 0
\text{.}\end{eqnarray*}

Thus $\int_{\mm}D_2D_1D_0H_2=0$ as well.
For $D_2D_0^2H_1$, we compute

\begin{eqnarray*}
& & \int_{\mm_T}D_2D_0^2H_1 \\
& = & \int_{(Z_{11})_T}\frac{2\psi\psi^2\la_0}{2(\la_0-\la_1)^2(\la_0-\la_1-2\psi)}
+\int_{(Z_{12})_T}\frac{2\psi\psi^2\la_1}{2(\la_1-\la_0)^2(\la_1-\la_0-2\psi)}\\
& = & 0+0\\
& = & 0
\text{.}\end{eqnarray*}

Thus $\int_{\mm}D_2D_0^2H_2=0$ as well.
For $D_1^3H_1$, we compute

\begin{eqnarray*}
& & \int_{\mm_T}D_1^3H_1 \\
& = & \int_{(Z_{5})_T}\frac{8(\la_0-\la_1)^3\la_1}{2(\la_1-\la_0)^4}
+\int_{(Z_{6})_T}\frac{8(\la_1-\la_0)^3\la_0}{2(\la_1-\la_0)^4}\\
&   & +\int_{(Z_{7})_T}\frac{(\la_0-\la_1)^3\la_1}{-(\la_0-\la_1)^4}
+\int_{(Z_{8})_T}\frac{(\la_1-\la_0)^3\la_0}{-(\la_0-\la_1)^4}\\
&   & +\int_{(Z_{9})_T}\frac{(\la_0-\la_1)^3\la_0}{-(\la_0-\la_1)^4}
+\int_{(Z_{10})_T}\frac{(\la_1-\la_0)^3\la_1}{-(\la_0-\la_1)^4}\\
&   & +\int_{(Z_{11})_T}\frac{8(\la_0-\la_1-\psi)^3\la_0}{2(\la_0-\la_1)^2(\la_0-\la_1-2\psi)}
+\int_{(Z_{12})_T}\frac{8(\la_1-\la_0-\psi)^3\la_1}{2(\la_1-\la_0)^2(\la_1-\la_0-2\psi)}\\
& = & \frac{4\la_1}{\la_0-\la_1}+\frac{4\la_0}{\la_1-\la_0}
-\frac{\la_1}{\la_0-\la_1}\\
&   & -\frac{\la_0}{\la_1-\la_0}-\frac{\la_0}{\la_0-\la_1}
-\frac{\la_1}{\la_1-\la_0}\\
&   & +\int_{(Z_{11})_T}\frac{4(\la_0-\la_1-3\psi)\la_0(1+2\psi/(\la_0-\la_1))}
{\la_0-\la_1}\\
&   & +\int_{(Z_{12})_T}\frac{4(\la_1-\la_0-3\psi)\la_1(1+2\psi/(\la_1-\la_0))}
{\la_1-\la_0}\\
& = & -4+1-1-\frac{4\la_0}{\la_0-\la_1}-\frac{4\la_1}{\la_1-\la_0}\\
& = & -4-4\\
& = & -8
\text{.}\end{eqnarray*}

For $D_1^3H_2$, we compute

\begin{eqnarray*}
& & \int_{\mm_T}D_1^3H_2 \\
& = & \int_{(Z_{5})_T}\frac{8(\la_0-\la_1)^3\la_1}{2(\la_1-\la_0)^4}
+\int_{(Z_{6})_T}\frac{8(\la_1-\la_0)^3\la_0}{2(\la_1-\la_0)^4}\\
&   & +\int_{(Z_{7})_T}\frac{(\la_0-\la_1)^3\la_0}{-(\la_0-\la_1)^4}
+\int_{(Z_{8})_T}\frac{(\la_1-\la_0)^3\la_1}{-(\la_0-\la_1)^4}\\
&   & +\int_{(Z_{9})_T}\frac{(\la_0-\la_1)^3\la_1}{-(\la_0-\la_1)^4}
+\int_{(Z_{10})_T}\frac{(\la_1-\la_0)^3\la_0}{-(\la_0-\la_1)^4}\\
&   & +\int_{(Z_{11})_T}\frac{8(\la_0-\la_1-\psi)^3\la_0}
{2(\la_0-\la_1)^2(\la_0-\la_1-2\psi)}
+\int_{(Z_{12})_T}\frac{8(\la_1-\la_0-\psi)^3\la_1}
{2(\la_1-\la_0)^2(\la_1-\la_0-2\psi)}\\
& = & \frac{4\la_1}{\la_0-\la_1}+\frac{4\la_0}{\la_1-\la_0}
-\frac{\la_0}{\la_0-\la_1}
-\frac{\la_1}{\la_1-\la_0}-\frac{\la_1}{\la_0-\la_1}
-\frac{\la_0}{\la_1-\la_0}\\
&   & +\int_{(Z_{11})_T}\frac{4(\la_0-\la_1-3\psi)\la_0(1+2\psi/(\la_0-\la_1))}
{\la_0-\la_1}\\
&   & +\int_{(Z_{12})_T}\frac{4(\la_1-\la_0-3\psi)\la_1(1+2\psi/(\la_1-\la_0))}
{\la_1-\la_0}\\
& = & -4-1+1-\frac{4\la_0}{\la_0-\la_1}-\frac{4\la_1}{\la_1-\la_0}\\
& = & -4-4\\
& = & -8
\text{.}\end{eqnarray*}

For $D_1^2D_0H_1$, we compute

\begin{eqnarray*}
& & \int_{\mm_T}D_1^2D_0H_1 \\
& = & \int_{(Z_{5})_T}\frac{4(\la_0-\la_1)^2(\la_1-\la_0)\la_1}{2(\la_1-\la_0)^4}
+\int_{(Z_{6})_T}\frac{4(\la_0-\la_1)^2(\la_0-\la_1)\la_0}{2(\la_1-\la_0)^4}\\
&   & +\int_{(Z_{11})_T}\frac{4(\la_0-\la_1-\psi)^2\psi\la_0}{2(\la_0-\la_1)^2(\la_0-\la_1-2\psi)}
+\int_{(Z_{12})_T}\frac{4(\la_1-\la_0-\psi)^2\psi\la_1}{2(\la_1-\la_0)^2(\la_1-\la_0-2\psi)}\\
& = & \frac{2\la_1}{\la_1-\la_0}+\frac{2\la_0}{\la_0-\la_1}\\
&   & +\int_{(Z_{11})_T}\frac{2(\la_0-\la_1-2\psi)\psi\la_0}
{(\la_0-\la_1)(\la_0-\la_1-2\psi)}\\
&   & +\int_{(Z_{12})_T}\frac{2(\la_1-\la_0-2\psi)\psi\la_1}
{(\la_1-\la_0)(\la_1-\la_0-2\psi)}\\
& = & 2+\frac{2\la_0}{\la_0-\la_1}+\frac{2\la_1}{\la_1-\la_0}\\
& = & 4
\text{.}\end{eqnarray*}

Thus $D_1^2D_0H_2=4$ as well.
For $D_1D_0^2H_1$, we compute

\begin{eqnarray*}
& & \int_{\mm_T}D_1D_0^2H_1 \\
& = & \int_{(Z_{5})_T}\frac{2(\la_0-\la_1)(\la_1-\la_0)^2\la_1}{2(\la_1-\la_0)^4}
+\int_{(Z_{6})_T}\frac{2(\la_1-\la_0)(\la_0-\la_1)^2\la_0}{2(\la_1-\la_0)^4}\\
&   & +\int_{(Z_{11})_T}\frac{2(\la_0-\la_1-\psi)\psi^2\la_0}
{2(\la_0-\la_1)^2(\la_0-\la_1-2\psi)}
+\int_{(Z_{12})_T}\frac{2(\la_1-\la_0-\psi)\psi^2\la_1}
{2(\la_1-\la_0)^2(\la_1-\la_0-2\psi)}\\
& = & \frac{\la_1}{\la_0-\la_1}+\frac{\la_0}{\la_1-\la_0}+0+0\\
& = & -1
\text{.}\end{eqnarray*}

Thus $\int_{\mm}D_1D_0^2H_2=-1$ as well.
For $D_0^3H_1$, we compute

\begin{eqnarray*}
& & \int_{\mm_T}D_0^3H_1 \\
& = & \int_{(Z_{3})_T}\frac{\la_0(\la_0-\la_1)^3/8}{-2(\la_1-\la_0)^4/4}
+\int_{(Z_{4})_T}\frac{\la_1(\la_1-\la_0)^3/8}{-2(\la_1-\la_0)^4/4}\\
&   & +\int_{(Z_{5})_T}\frac{(\la_1-\la_0)^3\la_1}{2(\la_1-\la_0)^4}
+\int_{(Z_{6})_T}\frac{(\la_0-\la_1)^3\la_0}{2(\la_1-\la_0)^4}\\
&   & +\int_{(Z_{11})_T}\frac{\psi^3\la_0}{2(\la_0-\la_1)^2(\la_0-\la_1-2\psi)}
+\int_{(Z_{12})_T}\frac{\psi^3\la_1}{2(\la_1-\la_0)^2(\la_1-\la_0-2\psi)}\\
& = & -\frac{\la_0}{4(\la_0-\la_1)}-\frac{\la_1}{4(\la_1-\la_0)}\\
&   & +\frac{\la_1}{2(\la_1-\la_0)}+\frac{\la_0}{2(\la_0-\la_1)}+0+0\\
& = & -\frac{1}{4}+\frac{1}{2}\\
& = & \frac{1}{4}
\text{.}\end{eqnarray*}

Thus $\int_{\mm}D_0^3H_2=\frac{1}{4}$ as well.
Finally, for $D_2^2H_1H_2$, we compute

\begin{eqnarray*}
& & \int_{\mm_T}D_2^2H_1H_2 \\
& = & \int_{(Z_{7})_T}\frac{(\la_0-\la_1)^2\la_0\la_1}{-(\la_0-\la_1)^4}
+\int_{(Z_{8})_T}\frac{(\la_0-\la_1)^2\la_0\la_1}{-(\la_0-\la_1)^4}\\
&   & +\int_{(Z_{9})_T}\frac{(\la_0-\la_1)^2\la_0\la_1}{-(\la_0-\la_1)^4}
+\int_{(Z_{10})_T}\frac{(\la_0-\la_1)^2\la_0\la_1}{-(\la_0-\la_1)^4}\\
&   & +\int_{(Z_{11})_T}\frac{4\psi^2\la_0^2}{2(\la_0-\la_1)^2(\la_0-\la_1-2\psi)}
+\int_{(Z_{12})_T}\frac{4\psi^2\la_1^2}{2(\la_1-\la_0)^2(\la_1-\la_0-2\psi)}\\
&   & +\int_{(Z_{13})_T}\frac{4(\la_0-\la_1)^2\la_1^2}{2(\la_0-\la_1)^4}
+\int_{(Z_{14})_T}\frac{4(\la_1-\la_0)^2\la_0^2}{2(\la_0-\la_1)^4}\\
& = & -\frac{\la_0\la_1}{(\la_0-\la_1)^2}-\frac{\la_0\la_1}{(\la_0-\la_1)^2}
 -\frac{\la_0\la_1}{(\la_0-\la_1)^2} \\
&   &-\frac{\la_0\la_1}{(\la_0-\la_1)^2}+0+0+\frac{2\la_1^2}{(\la_0-\la_1)^2}
+\frac{2\la_0^2}{(\la_0-\la_1)^2}\\
& = & \frac{2(\la_0-\la_1)^2}{(\la_0-\la_1)^2}\\
& = & 2
\text{.}\end{eqnarray*}

\renewcommand{\baselinestretch}{1}
\chapter{Macaulay 2 code and output for computing gravitational correlators}
\label{sec:m2}

In order to find the genus zero, degree two,
two-point gravitational correlators of $\PP^1$ using the presentation
(\ref{geomprez}), we need to compute integrals of products of
 the $\psi$-classes and
hyperplane pullback classes in $A^*(\mm_{0,2}(\PP^1,2))$. Since
$\dim(A^4(\mm_{0,2}(\PP^1,2)))=1$, such a project reduces to two steps.

\begin{enumerate}
\item \label{first}
Find the integral of one class in $A^4(\mm_{0,2}(\PP^1,2))$.
\item \label{second}
Express the relevant products as multiples of the class above.
\end{enumerate}

The integral of such a product has value equal to the product of the 
numbers from (\ref{first}) and (\ref{second}) above. Part (\ref{first})
was carried out in Section \ref{sec:lla},
where, for example, the degree of $D_1^4$ was noted to be $-20$. Part
(\ref{second}) uses relations from the presentation to identify any degree
four monomial with a multiple of $D_1^4$. This
could be done by hand, but it is more efficient to use a computer algebra
system.  Macaulay 2 (\cite{GS}) is especially well-suited for this purpose.
The program, as well as helpful tutorials and documentation, is available
at {\tt http://www.math.uiuc.edu/Macaulay2/}.

This appendix records both the input code and output values involved in
calculating the necessary products in Macaulay 2. Note that $D_1$ is 
listed last
among the generators. This causes the output from the products to be
expressed in terms of $D_1^4$.

Some lines have been wrapped manually to remain within margins. Otherwise
the interface is preserved.
Take special note of the minus signs on the output; they are easy to miss.

\renewcommand{\baselinestretch}{1}
\begin{verbatim}

Macaulay 2, version 0.9.2
--Copyright 1993-2001, D. R. Grayson and M. E. Stillman
--Singular-Factory 1.3b, copyright 1993-2001, G.-M. Greuel, et al.
--Singular-Libfac 0.3.2, copyright 1996-2001, M. Messollen

i1 : R=QQ[D0,D2,H1,H2,P1,P2,D1]

o1 = R

o1 : PolynomialRing

i2 : I=ideal(H1^2,H2^2,D0*P1,D0*P2,D2-P1-P2,1/4*D1+1/4*D2+D0-H1-P1,
1/4*D1+1/4*D2+D0-H2-P2,(D1+D2)^3,D1*P1*P2)

              2    2                                   1           
o2 = ideal (H1 , H2 , D0*P1, D0*P2, D2 - P1 - P2, D0 + -*D2 - H1 
                                                       4         
       1         1                1       3      2           2 
- P1 + -*D1,D0 + -*D2 - H2 - P2 + -*D1, D2  + 3D2 D1 + 3D2*D1
       4         4                4
     3
 + D1 , P1*P2*D1)

o2 : Ideal of R

i3 : S=R/I

o3 = S

o3 : QuotientRing

i4 : P1^4

      3   4
o4 = --*D1
     80

o4 : S

i5 : P1^3*H1

       1   4
o5 = ---*D1
      80

o5 : S

i6 : P1^3*H2

\end{verbatim}

\pagebreak

\begin{verbatim}


      1   4
o6 = --*D1
     16

o6 : S

i7 : P1^3*P2

       3   4
o7 = ---*D1
      80

o7 : S

i8 : P1^2*H1*P2

      1   4
o8 = --*D1
     80

o8 : S

i9 : P1^2*H1*H2

       1   4
o9 = ---*D1
      40

o9 : S

i10 : P1^2*P2^2

        1   4
o10 = ---*D1
       16

o10 : S

i11 : P1^2*P2*H2

       3   4
o11 = --*D1
      80

o11 : S

i12 : P1*H1*P2*H2

        1   4
o12 = ---*D1
       40

o12 : S
\end{verbatim}

\end{document}